\documentclass{article}

\usepackage{arxiv}

\usepackage[utf8]{inputenc} 
\usepackage[T1]{fontenc}    
\usepackage{hyperref}       
\usepackage{url}            
\usepackage{booktabs}       
\usepackage{amsfonts}       
\usepackage{nicefrac}       
\usepackage{microtype}      
\usepackage{amsmath}
\usepackage{graphicx}
\usepackage{amsthm,amsmath,amsfonts,amssymb}
\usepackage[numbers]{natbib}
\usepackage{amsmath}
\usepackage{color, graphicx}
\usepackage[colorinlistoftodos]{todonotes}
\usepackage{xcolor}
\usepackage{hyperref}
\usepackage{verbatim}
\usepackage{color,soul}
\usepackage{amsthm}
\usepackage{amsmath}
\usepackage{tikz}
\usepackage{float}
\usepackage{amsmath}
\usepackage{float}
\usepackage{siunitx}
\usepackage{amssymb}
\usepackage{comment}
\usepackage{graphicx}
\usepackage{svg}
\usepackage{marvosym}
\usepackage{wrapfig}

\usepackage[toc,page]{appendix}
\usepackage{bbm}
\usepackage{graphicx}
\usepackage{amsbsy}
\usepackage{lipsum}
\usepackage{multicol}
\usepackage{amsthm}
\usepackage{algorithm,algorithmic}
\usepackage{pifont}
\usepackage{bbm}
\usepackage{bm}
\usepackage{dsfont}
\usepackage{bigints}
\usepackage{mathtools}
\usepackage{graphicx}
\usepackage{svg}
\usepackage{calc}
\usepackage{scalerel}
\usepackage{accents}
\usepackage{mwe}
\usepackage{graphicx,wrapfig,lipsum}
\usepackage{mathtools}
\usepackage{mathrsfs}
\usetikzlibrary{mindmap,backgrounds}
\usepackage{blindtext}
\usepackage{hyperref}
\usepackage{pdfpages}
\allowdisplaybreaks

\usepackage{enumitem}
\usepackage{letltxmacro}
\usepackage{nameref,hyperref}
\usepackage[capitalize]{cleveref}

\newlist{thmlist}{enumerate}{1}
\setlist[thmlist]{label=\alph{thmlisti}., ref=\thetheorem.\alph{thmlisti}.,noitemsep}

\numberwithin{equation}{section}
\newcommand{\pb}[0]{\mathbb{P}}
\newcommand{\ex}[0]{\mathbb{E}}
\newcommand{\ind}[0]{\mathbf{1}}
\newcommand{\ubar}[1]{\underaccent{\bar}{#1}}

\DeclareMathOperator*{\argmax}{arg\,max}

\newlength\myindent
\newcommand\numberthis{\addtocounter{equation}{1}\tag{\theequation}}
\theoremstyle{plain}
\newtheorem{theorem}{Theorem}

\newtheorem{lemma}{Lemma}

\theoremstyle{remark}
\newtheorem{example}{Example}

\newtheorem{remark}{Remark}

\newtheorem*{definition*}{Definition}

\Crefname{theorem}{Theorem}{Theorems}

\title{Adaptive and Efficient Isotonic Estimation in Wicksell's Problem}

\author{ Francesco Gili \\
	Delft University of Technology,\\
	Mekelweg 4, Delft 2628CD, The Netherlands. \\
	\texttt{F.Gili@tudelft.nl} \\
	\And
	Geurt Jongbloed \\
	Delft University of Technology,\\
	Mekelweg 4, Delft 2628CD, The Netherlands. \\
	\texttt{G.Jongbloed@tudelft.nl} \\
    \And
    Aad van der Vaart \\
	Delft University of Technology,\\
	Mekelweg 4, Delft 2628CD, The Netherlands. \\
	\texttt{A.W.vanderVaart@tudelft.nl} \\
}

\hypersetup{
pdftitle={Adaptive and Efficient Isotonic Estimation in Wicksell's Problem - F. Gili, G. Jongbloed, A. van der Vaart},
pdfsubject={},
pdfauthor={Francesco ~Gili, Geurt ~Jongbloed, Aad ~van der Vaart},
pdfkeywords={Wicksell's problem, Isotonic estimation, Efficiency theory},
}

\begin{document}
\maketitle

\begin{abstract}
We consider nonparametric estimation in Wicksell's problem which has relevant applications in astronomy for estimating the distribution of the positions of the stars in a galaxy given projected stellar positions and in material sciences to determine the 3D microstructure of a material, using its 2D cross sections. In the classical setting, we study the isotonized version of the plug-in estimator (IIE) for the underlying cdf $F$ of the spheres' squared radii. This estimator is fully automatic, in the sense that it does not rely on tuning parameters, and we show it is adaptive to local smoothness properties of the distribution function $F$ to be estimated. Moreover, we prove a local asymptotic minimax lower bound in this non-standard setting, with $\sqrt{\log{n}/n}$-asymptotics and where the functional $F$ to be estimated is not regular. Combined, our results prove that the isotonic estimator (IIE) is an adaptive, easy-to-compute, and efficient estimator for estimating the underlying distribution function $F$. 
\end{abstract}


\section{Introduction}

\begin{figure}
  \centering
  \includegraphics[width=7.5cm]{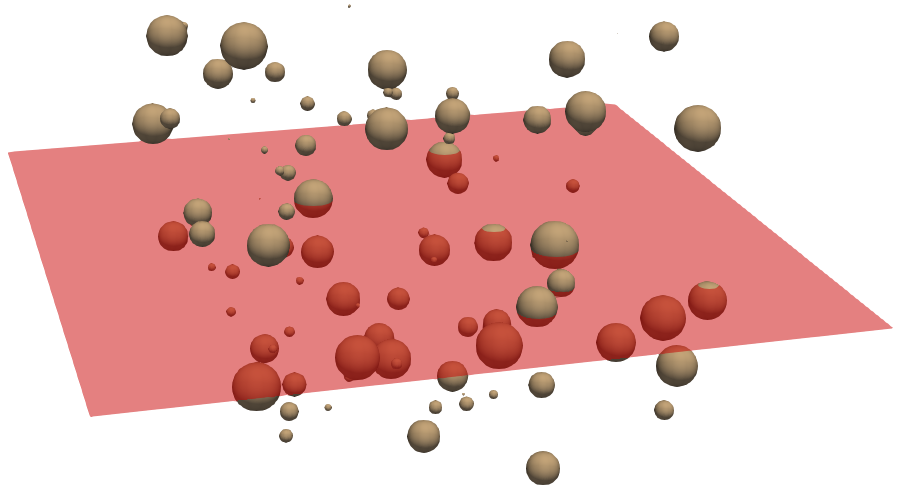}
  \includegraphics[width=3.8cm]{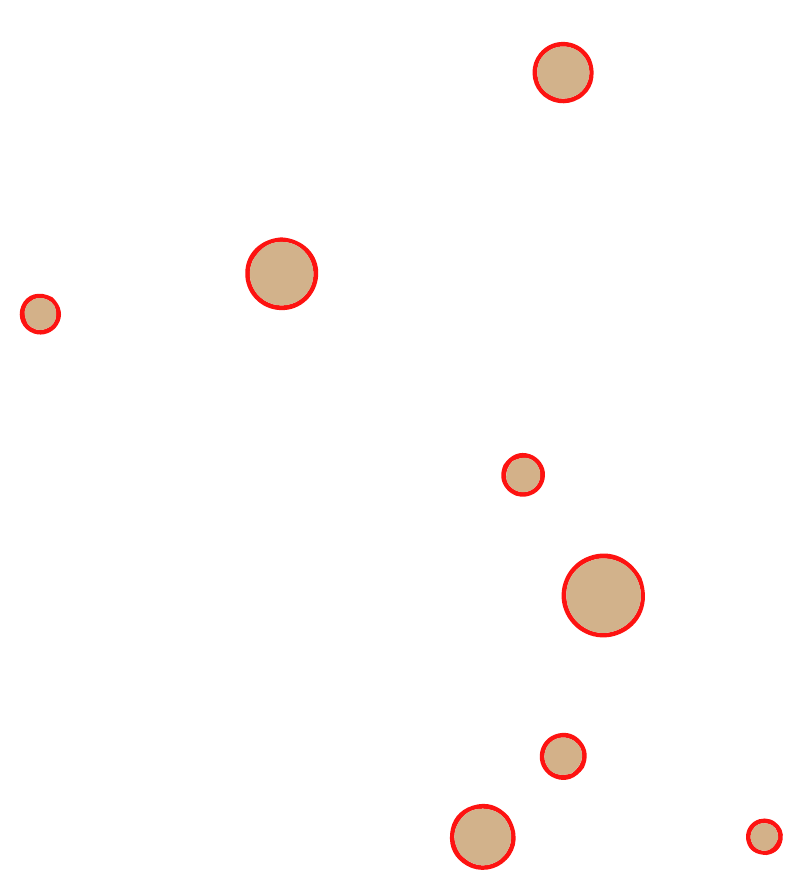}
  \setlength{\belowcaptionskip}{-4pt}
  \caption{Sampling in Wicksell's problem. Left panel: 3D configuration with spheres and the intersecting plane. Right panel: observable circular profiles.}\label{fig: 1 sampling}
\end{figure}

In the field of stereology, scientists are usually concerned with the study of three-dimensional properties of materials and objects by interpreting their two-dimensional cross-sections. This differs from three-dimensional reconstructions in that some 3D quantities can be estimated without expensive 3D reconstructions. In 1925, Sven Wicksell published an article \cite{14} that studied the following classical stereological inverse problem. Suppose that a number of spheres are embedded in an opaque three-dimensional medium. Because the medium is opaque, we are not able to observe the spheres directly; however, we can observe cross-sections of the medium (see Fig. \ref{fig: 1 sampling}), which show circular sections of the spheres that happen to be cut by the plane. This model has present-day relevant applications in astronomy for estimating the positions of stars that are spherically symmetrically distributed in galactic clusters, given projected stellar positions (see \cite{19}), or to determine the 3D microstructure of a material given its cross sections, useful to assess its structural properties like resistance and fracture toughness (see \cite{21,22} for applications of this setting). 

Assume the spheres' squared radii are i.i.d.\ realizations from a random variable $X$ with distribution function $F$, the object to be estimated. Following the same notation as in \cite{1}, we have that a version of the density $g$ of the observable squared circle radii \cite{15} is given by:
\begin{align}\label{eq: density of the observations}
g(z)=\frac{1}{2 m_0} \int_{z}^{\infty} \frac{d F(x)}{\sqrt{x-z}},    
\end{align}
where $ 0 < m_0=\int_0^{\infty} \sqrt{y} d F(y) < \infty$ is the expected sphere radius under $F$.
Wicksell \cite{14} inverted this equation, by recognizing an Abel-type integral, and found an expression for $F$ in terms of $g$:
\begin{align}\label{eq: inv}
F(x)=1-\frac{\int_x^{\infty}(z-x)^{-1 / 2} g(z) \: dz}{\int_0^{\infty} z^{-1 / 2} g(z) \: dz} = 1 - \frac{V(x)}{V(0)}, \quad x>0
\end{align}
where:
$$
V(x) \vcentcolon=\int_x^{\infty} \frac{g(z)}{\sqrt{z-x}} \: d z,
$$
and $F(x) = 0$ for all $x \leq 0$. Therefore, in order to estimate $F$ at a point $x>0$, the essential object to be estimated is the function $V$ at $x$ and at $0$. Given a sample $Z_1, \ldots, Z_n$ from the density $g$ and for $\mathbb{G}_n$ the empirical distribution function based on this sample, a naive (empirical plug-in) estimator of $V$ is:
\begin{align}\label{eq: naive estimator of V}
    V_n(x) \vcentcolon=\int_x^{\infty} \frac{d\mathbb{G}_n(z)}{\sqrt{z-x}}.
\end{align}
The estimator for $F$ is then defined by $F_n(x) = 1- V_n(x)/V_n(0)$. As a function of $x$, this estimator has some undesirable properties, being non-monotone and having infinite discontinuities at a set of points becoming dense in the support of the density $g$. This led in \cite{1} to the definition of the so-called "isotonic inverse estimator" (IIE) $\hat{F}_n$ for $F$, which is obtained by replacing $V_n$ by its isotonic projection. This can be calculated as the right continuous derivative of the least concave majorant of the primitive function of $V_n$, which is:
$$
U_n(x)=\int_0^x V_n(y) \, dy=2 \int_0^{\infty} \sqrt{z} \, d\mathbb{G}_n(z)-2 \int_x^{\infty} \sqrt{z-x} \, d \mathbb{G}_n(z).
$$
Note that $U(x) = \int_{0}^x V(y) \, dy$ is a concave and increasing function on $[0,\infty)$, whereas $U_n$ is not. Taking the least concave majorant of $U_n$ on $[0,\infty)$ yields an increasing concave function. Denote by $\hat{V}_n$ the right derivative of this least concave majorant of $U_n$, and note that this will be a decreasing function on $[0,\infty)$. Then $\hat{F}_n$ is given by: 
\begin{align}\label{eq: isotonic estimator of F}
\hat{F}_n(x) := 1 - \hat{V}_n(x)/\hat{V}_n(0).
\end{align}
This estimator has the important properties of being monotone, right-continuous and does not depend on tuning parameters, unlike other estimators in the literature (cf.\ \cite{10}). 

In this paper, we prove that the IIE $\hat{F}_n$ has the further properties of automatically adapting to the local smoothness of $F$ and of being asymptotically efficient. No other estimator in the literature simultaneously has these properties. Consequently, the IIE emerges as the unequivocal choice, surpassing alternative estimators previously proposed in Wicksell's problem. First, we prove that the IIE is asymptotically normal and the smoother the cdf $F$ is at $x$ and $0$, the smaller the asymptotic variance automatically attained by the IIE (Theorem  \ref{thm: main thm iso 1}). Here, the smoothness of $F$ at $x$ and $0$ is measured by considering constants $K >0$ and $\gamma > \frac{1}{2}$ such that:
\begin{align*}
    & \quad H_{x}(\delta ) \vcentcolon = \int_{0}^{1} \left( F(x+u\delta)-F(x) \right) \, du \sim \operatorname{sgn}(\delta) |\delta|^{\gamma} K, \quad \quad \text{as} \: \delta \rightarrow 0. \numberthis \label{eq: condition roughness intro} 
\end{align*}
The adaptivity extends also to the situation where $F$ is constant around $x$ and $0$. Then a faster $\sqrt{n}$-rate of convergence is attained (Theorem \ref{thm: main thm iso 2}). Note that previously there were no results in the literature for estimators that attained $\sqrt{n}$ rate in Wicksell's problem. In the case in which $F$ is constant only around $x$ but not around $0$, or vice versa, the rate will be again $\sqrt{n} \: (\log{n})^{-1/2}$. Finally, we prove a lower bound for the minimax risk of any estimator sequence, also depending on local smoothness properties of $F$. Combining these results shows that the IIE is an easy-to-compute, adaptive, and efficient estimator for $F$. We now state the first main result and its \hyperlink{proof of main thm}{proof} is given in Appendix A. 

\begin{theorem}\label{thm: main thm iso}
Let $F$ be the distribution function of the squared sphere radii and $g$ the corresponding density of the squared circle radii according to \eqref{eq: density of the observations}. Assume $\int_{0}^{\infty} y^{\frac{3}{2}} \: dF(y)  < \infty$ and choose $x >0$. 
\begin{thmlist}
    \item If $F$ satisfies condition \eqref{eq: condition roughness intro} both at $0$ and at $x$ for $\gamma_0 > \frac{1}{2}$ and $\gamma_x > \frac{1}{2}$ respectively and $g(x)<\infty$ holds true, then, as $n \rightarrow \infty$: \label{thm: main thm iso 1}
    \begin{align*}
    & \sqrt{\frac{n}{\log n}}\left( \hat{F}_n(x)-F(x)\right) \rightsquigarrow N\left(0, \frac{4m^2_0}{\pi^2} \left(\frac{g(x)}{2\gamma_x} + (1-F(x))^2 \frac{g(0)}{2\gamma_0} \right)\right).
    \end{align*} 
    \item If $F$ is constant around $x$ and $0$, with $K_x$ and $K_0$ the biggest compact intervals in $[0,\infty)$ that contain $x$ and $0$, respectively, so that $F$ is constant, on $K_x$ and $K_0$, with $K_x \cap K_0 = \emptyset$, then, as $n \rightarrow \infty$: \label{thm: main thm iso 2}
    \begin{align*}
        \sqrt{n} \left( \hat{F}_n(x) - F(x) \right) \rightsquigarrow \frac{2m_0}{\pi} \Big(L_x+ (1-F(x))L_0 \Big), 
    \end{align*}
    where, for any $a,b \in \mathbb{R}$, the joint distribution of $L_x$ and $L_0$ is given by:
    \begin{align*}
        \pb & \left( L_x \leq a, L_0 \leq b \right) \\
        &= \pb \bigg( \argmax_{s \in K_x} \left\{ \mathbb{Z}_x(s) - a s\right\} \leq x, \argmax_{s \in K_0} \left\{ \mathbb{Z}_0(s) - b s\right\} \leq 0 \bigg). 
    \end{align*}
    Here $\mathbb{Z}_p(s) = \mathbb{Z}(s) - \mathbb{Z}(p),$ for $s \geq 0$ and $p \in \{0,x\}$. $\mathbb{Z}$ is a continuous centered Gaussian Process with covariance structure given by:
    \begin{align}\label{eq: covariance gauss process}
        \mathrm{Cov} \left( \mathbb{Z}(t),\mathbb{Z}(s) \right) = \mathrm{Cov} \Big( \sqrt{\smash[b]{(Z)_{+}}} - \sqrt{\smash[b]{(Z_t)_{+}}}, \sqrt{\smash[b]{(Z)_{+}}} - \sqrt{\smash[b]{(Z_s)_{+}}} \Big),
    \end{align}
    where $Z_t = Z - t$ and $Z_s = Z - s$.
\end{thmlist}
\end{theorem}

The second main result of this paper consists of showing that the IIE is asymptotically efficient in this setting, i.e.\ the attained asymptotic variance is the smallest attainable. This kind of result is obtained by proving a minimax lower bound for any estimator sequence estimating $F(x)$. Here we give a weaker (but perhaps more insightful) version of this result and its \hyperlink{proof of thm 2}{proof} is given in Appendix B. 

\begin{theorem}\label{thm: LAM weaker}
    Let the conditions of Theorem \ref{thm: main thm iso 1} hold, assume that $F$ for such $\gamma_0, \gamma_x$ satisfies \eqref{eq: condition roughness infimum} and that $g$ is continuous at $x$ and $0$. Let $\ell$ be a bowl-shaped loss function, i.e., $\ell$ is non-negative, symmetric and subconvex, then for every estimator sequence $\left( F_n(x)\right)_{n \in \mathbb{N}}$:
    \begin{align}\label{eq: LAM simplified}
        \lim_{\delta \downarrow 0} \: \: \liminf_{n \rightarrow \infty} \, \sup_{\| F^{\prime } - F \|_{\infty} < \delta} \, \ex_{\scriptscriptstyle{G_{{\scaleto{F^{\prime}}{3.7pt}}}}} \, \ell\left(\sqrt{\frac{n}{\log{n}}}\left(F_n(x)- F^{\prime}(x)\right)\right)  \geq \ex \ell(L),
    \end{align}
where $L$ possesses the normal distribution in Theorem \ref{thm: main thm iso 1}.
\end{theorem}

\subsubsection*{Motivation, connections with the literature and results.} 

There is a vast literature on Wicksell's problem, see for instance \cite{1,2,3,9,10, 15, 16, 17, 18, 19, 20}. For a detailed introduction to Wicksell's problem we refer the reader to \cite{9}. The setting of Wicksell's problem has several present-day applications as the estimation of the distribution of stars in a galactic cluster (cf.\ \cite{19}) or the estimation of the 3D microstructure of materials (cf.\ \cite{21,22}). In applications, it is common to see the use of the so-called Saltykov method or methods based on numerical discretization. These methods are far from being efficient and in some cases not even consistent, making our work of practical relevance. The estimation of $F$ given data from $g$ is of interest not only for stereological procedures, but also from the mathematical point of view. Consider the class of inverse problems that entails the estimation of a functional: $\Psi_\alpha(x ; g)=\int_x^{\infty}(z-x)^{-\alpha} g(z) \, dz, \: \: 0<\alpha<1$ based on samples from $g$. This type of functional estimation emerges in a vast variety of applications. For $\alpha < \frac{1}{2}$, the above gives a regular statistical functional (cf.\ Theorem 4 in \cite{24}), whereas for $\alpha>\frac{1}{2}$ gives instances of irregular functionals. Wicksell's problem entails the boundary case for $\alpha = 1/2$, which makes the problem mathematically intriguing, both because of the unusual rate of convergence and because of the non-standard efficiency theory that in this paper is developed.

This work has a special connection with the work of \cite{1} and \cite{10}. In \cite{1}, the authors observed a peculiar behavior of the isotonic estimator (IIE) when compared to the natural naive (empirical plug-in) estimator for $F$, given by $F_n(x) = 1- V_n(x)/V_n(0)$. Asymptotic properties for fixed $x$ like consistency and pointwise asymptotic normality, after rescaling by the rate $(\log{n})^{-1/2}\sqrt{n} $ (cf.\ \cite{1, 15}), have been derived for this naive estimator. However, this estimator, viewed as a function of $x$, is ill-behaved, as we already pointed out. The curious behavior observed in \cite{1} is that the rate of estimation is the same as for the isotonic estimator and, in both cases, a normal limiting behavior in distribution appears. However, the asymptotic variance $\sigma^2(F) = \frac{4m^2_0}{\pi^2}(g(x) + g(0)(1-F(x))^2)$ for the naive estimator is exactly twice the one for the isotonic estimator. The authors of \cite{10} propose an explanation for this phenomenon. They show, if $F$ belongs to the Hölder class $\Lambda_\gamma$ of functions having finite norm \eqref{eq: holder norm definition}, for some $\gamma > \frac{1}{2}$, then an efficient estimator attains the asymptotic variance: $\sigma^2(F)/(2\gamma)$. In \cite{1}, the authors assume that $F$ possesses a bounded density $f$ and therefore this corresponds to the case $\gamma = 1$, leading indeed to a variance cut by 2. 

In this paper, we show that it is not needed to assume that $F$ belongs to some global Hölder class to obtain such asymptotic variances. It is enough to consider the local smoothness of the distribution function $F$ at the points $x$ and $0$. We will measure the local smoothness in some technical conditions, like \eqref{eq: condition roughness} and \eqref{eq: condition roughness infimum}, where some $\gamma_x > \frac{1}{2}$ and $\gamma_0 > \frac{1}{2}$ appear. These can be interpreted as the tight degree of Hölder smoothness at the points $x$ and $0$. Indeed, also in this case, for $F$ with density $f$ (i.e.\ differentiable at $x$ and $0$) we retrieve the $\gamma_x = \gamma_0 = 1$ that explains the peculiar behavior observed in \cite{1}. 

Moreover, the authors of \cite{10} propose an estimator based on Kernel Density Estimation which attains the variance $\sigma^2(F)/(2\gamma)$ (where $\gamma$ is as in their interpretation) and which is, to the best of our knowledge, the only other efficient estimator in the literature. However, the estimator they propose depends directly on the degree $\gamma$ of smoothness, which should be guesstimated before using the estimator. This means that the KDE estimator is not adaptive, in the sense that it does not recognize automatically the smoothness of $F$ at the estimation point. Furthermore, they derive an asymptotic minimax lower bound assuming $F$ is in the Hölder class $\Lambda_\gamma$, for $\gamma > \frac{1}{2}$. Interestingly, also in \cite{10} two local conditions at $x$ and $0$ are needed to prove the minimax theorem (cf. (24) and (25)). 

In this paper, we do two main things. On the one hand, we show that the isotonic estimator is adaptive in the sense specified above. The attained asymptotic variance is adaptive with respect to the local smoothness of $F$ at $x$ and $0$. On the other hand, we prove a lower bound on the minimax risk for any estimator sequence. The minimax theorem holds as well under local conditions of smoothness of $F$ around $x$ and $0$, as for the asymptotic normality of the isotonic estimator. If all conditions are met, this shows that the isotonic estimator is asymptotically efficient. Thus the IIE emerges as the unequivocal choice, surpassing alternative estimators previously proposed in Wicksell's problem.

Before continuing, we want to mention some work related to the present paper. The first time the naive plug-in estimator was studied is in \cite{15}. Later in \cite{1} the proof of the asymptotic normality of this estimator was simplified considerably. However, viewed as a function, the naive plug-in estimator is severely ill-behaved as we already pointed out. Naturally, one may wonder why not use the Nonparametric Maximum Likelihood Estimator. For Wicksell's problem, the usual definition of the NPMLE is also ill-posed, as the log-likelihood can attain value infinity. However, work has been done, for instance in \cite{3} and \cite{18}, to make the NPMLE an employable estimator. For results about confidence intervals in Wicksell's problem, we mention: \cite{19, 20}. On the other hand for KDE-based estimators, we refer to \cite{16, 10}. For extensions to other convex bodies other than spheres see \cite{23}. Finally, regarding efficient estimation, it is worth mentioning that one of the authors of \cite{10} did some further work on Wicksell's problem in \cite{17} where it was shown how to attain efficient smoothing for the plug-in estimator under square losses.

\section{Isotonic Estimator}

In this section, we prove two main theorems that will imply Theorem \ref{thm: main thm iso}. In order to prove Theorem \ref{thm: main thm iso 1}, it is enough to look at the behavior of $\hat{V}_n(x)$ and then apply the same result also to $\hat{V}_n(0)$. We now state the conditions needed for Theorem \ref{thm: main thm iso} and \ref{thm: normality Holder}. First, a natural assumption is a moment constraint for $Z$, which can be translated into a moment condition on $X$ (for $0 < m_0 < \infty$): 
\begin{align*}
    \quad \int_{0}^{\infty} z \, g(z) \, dz = \frac{2}{3m_0} \int_{0}^{\infty} y^{\frac{3}{2}} \: dF(y)  < \infty. \numberthis \label{eq: finite first moment}
\end{align*}
Additionally, for $x \geq 0$, we need our version of $g$ to be well defined, i.e.:
\begin{align*}
    \frac{1}{2m_0}\int_{x}^{\infty} \frac{dF(s)}{\sqrt{s - x}}  < \infty. \numberthis \label{eq: condition finite g at x}
\end{align*}
Furthermore, we measure the smoothness of $F$ at $x$; $\exists$ $K >0$ and $\gamma > \frac{1}{2}$ such that:
\begin{align*}
    & \quad H_{x}(\delta ) \vcentcolon = \int_{0}^{1} \left( F(x+u\delta)-F(x) \right) \, du \sim \operatorname{sgn}(\delta) |\delta|^{\gamma} K, \quad \quad \text{as} \: \delta \rightarrow 0. \numberthis \label{eq: condition roughness} 
\end{align*}
Before proving Theorem \ref{thm: main thm iso 1}, we show that condition \eqref{eq: condition roughness} entails the local smoothness of the distribution function $F$ by looking at the following examples \ref{example: cdf holder} and \ref{example: cdf discrete holder} below.
\begin{example}\label{example: cdf holder}
For $\gamma > \frac{1}{2}$ consider the class of distribution functions locally behaved as follows. Let $F$ be so that there exists a neighborhood $U_x$ of $x$ and some $K>0$:
\begin{align*}
    F(u) = \begin{cases}
        F(x) + K(u-x)^{\gamma}, \quad \quad \text{for} \: u \geq x, \: u \in U_x \\
        F(x) - K(x-u)^{\gamma}, \quad \quad \text{for} \: u < x, \: u \in U_x
    \end{cases}
\end{align*}
One can think of such distribution functions as tightly Hölder smooth of degree $\gamma$ at the point $x$. This class is of particular importance because it illustrates what type of functions we are thinking of when we derive results under smoothness conditions at the point $x$ and 0. Such functions satisfy both conditions \eqref{eq: condition roughness} and \eqref{eq: condition roughness 2}. Condition \eqref{eq: condition roughness 2} follows immediately. For condition \eqref{eq: condition roughness} note that as $|\delta| \downarrow 0$, small enough so that $[x-|\delta|, x + |\delta|] \subset U_x$:
\begin{align*}
    H_x(\delta) = \int_{0}^{1} \left( F(x+u\delta )-F(x) \right) \, du =  \operatorname{sgn}(\delta) |\delta|^\gamma K \int_{0}^{1} u^{\gamma} \, du + o(\delta).
\end{align*}
\end{example}

\begin{example}\label{example: cdf discrete holder}
For any measurable set $A$, consider the discrete probability measure $\mu_x$:
\begin{align*}
    \mu_x(A) = \frac{1}{2} \sum_{i=1}^{\infty} p_i \delta_{x+t_i}(A) + \frac{1}{2} \sum_{i=1}^{\infty} p_i \delta_{x-t_i} (A),
\end{align*}
for sequences $p_i$ of probability weights and $t_i \downarrow 0$. In particular set: $p_i = \frac{6}{\pi^2} i^{-2}$ and $t_i = i^{-\frac{1}{\gamma}}$ where $\gamma > 1/2$ and $\delta_x(\cdot)$ indicates the Dirac measure at $x$. Denote by $F$ the associated cdf, so that, because $\gamma > 1/2$, we have:
\begin{align*}
    \int_{x}^{\infty} \frac{dF(s)}{\sqrt{s-x}} = \frac{3}{\pi^2} \sum_{i=1}^{\infty} \frac{1}{i^2}\frac{1}{\sqrt{t_i}} = \frac{3}{\pi^2} \sum_{i=1}^{\infty} \frac{1}{i^{2- \frac{1}{2\gamma}}} < \infty.
\end{align*}
Moreover, we have for $\delta \downarrow 0$ (similar reasonings hold also for a generic sequence $|\delta| \downarrow 0$, we consider $\delta >0$ for simplicity):
\begin{align*}
    H_x(\delta) = \int_{0}^{1} \left( F(x+u\delta) - F(x) \right) \, du = \sum_i \frac{p_i}{2} \int_{0}^{1} \ind_{t_i \leq u \delta} \, du = \sum_i \frac{p_i (\delta - t_i)_{+}}{2\delta},
\end{align*}
\begin{figure}
  \vspace{-0.2cm}
  \centering
  \includegraphics[width = 6cm]{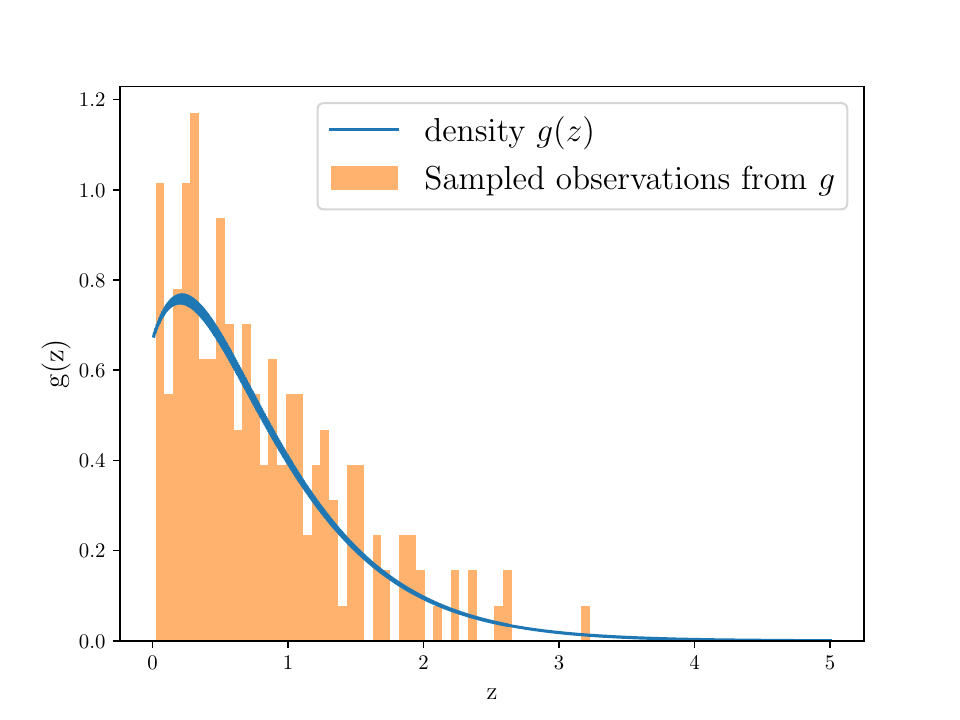}
  \includegraphics[width = 6cm]{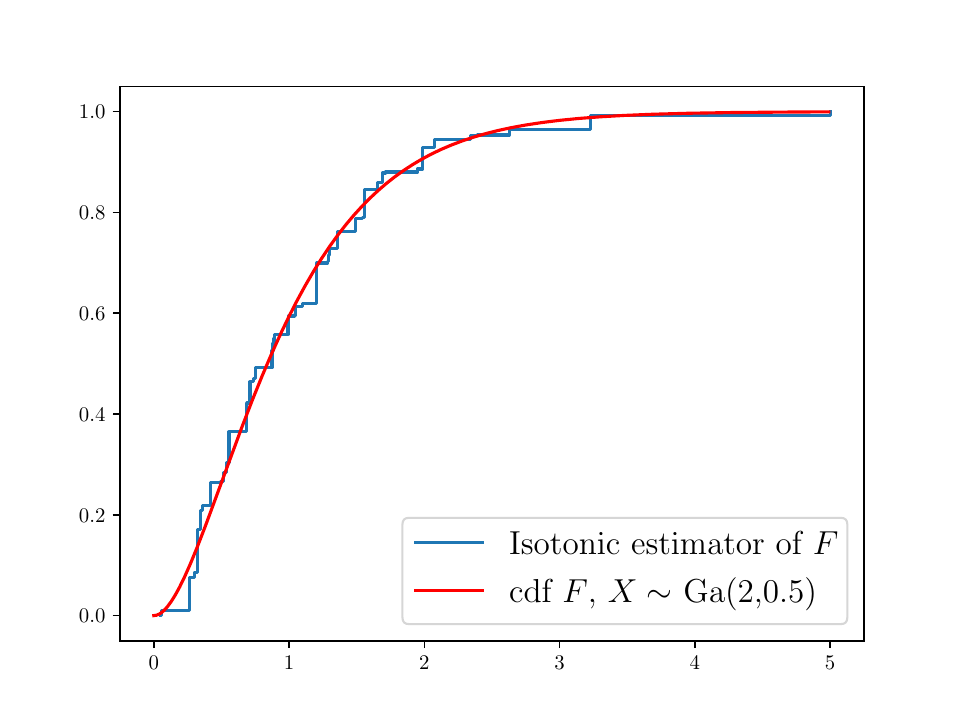}
  \setlength{\abovecaptionskip}{-4pt}
  \setlength{\belowcaptionskip}{-8
  pt}
  \caption{200 observed samples from $g$ and isotonic estimator for $F \sim \operatorname{Ga}(2,0.5)$.}\label{fig: estimator performance}
\end{figure}
which means $t_i< \delta$, thus $i > \frac{1}{\delta^{\gamma}}$. Therefore $H_x(\delta) = \sum_{i > 1/\delta^{\gamma}} \frac{p_i(\delta - i^{-1/\gamma})}{2\delta}$. Now we prove the following fact, as $\delta \downarrow 0$ and $k \geq 2$:
\begin{align*}
    \frac{1}{(\delta^{\gamma})^{k-1}} \sum_{i> \frac{1}{\delta^{\gamma}}} \frac{1}{i^k} \stackrel{\delta \downarrow 0}{\longrightarrow} \frac{1}{k-1}
\end{align*}
This follows directly from the fact that for $\alpha \rightarrow \infty$ (e.g.\ $\alpha = 1/\delta^{\gamma}$) we have:
\begin{align*}
    \frac{1}{k-1} \frac{1}{(\alpha+1)^{k-1}} = \int_{\alpha + 1}^{\infty} \frac{1}{x^k} dx \leq \sum_{i > \alpha} \frac{1}{i^k} \leq \int_{\alpha-1}^{\infty} \frac{1}{x^k} dx = \frac{1}{k-1} \frac{1}{(\alpha - 1)^{k-1}}
\end{align*}
Moreover $(\alpha/(\alpha \pm 1))^{k-1} \rightarrow 1$. Therefore as $\delta \downarrow 0$ we have:
\begin{align*}
    H_x(\delta) = \frac{6}{\pi^2} \left( \sum_{i > 1/ \delta^{\gamma}} \frac{1}{i^2} - \frac{1}{\delta} \sum_{i > 1/ \delta^{\gamma}} \frac{1}{i^{2+ \frac{1}{\gamma}}}\right) \sim  \frac{6}{\pi^2 (1+\gamma)} \delta^{\gamma}.
\end{align*}
This proves the fact that also for a discrete cdf the required conditions for Theorem \ref{thm: normality Holder} can be satisfied.
\end{example}

In Figure \ref{fig: estimator performance} the reader can observe the performance of the isotonic estimator (IIE).

\begin{theorem}\label{thm: normality Holder}
Let $x \geq 0$, suppose that \eqref{eq: finite first moment}, \eqref{eq: condition finite g at x} and \eqref{eq: condition roughness} hold true and $g(x)$ is given by \eqref{eq: density of the observations}. Then:
\begin{align*}
    \sqrt{\frac{n}{\log n}}\left(\hat{V}_n(x)-V(x)\right) \rightsquigarrow N\left(0, \frac{1}{2 \gamma} g(x)\right) \quad \text { as } n \rightarrow \infty \text {.}
\end{align*}
\end{theorem}
\begin{proof}
We generalize the proof of \cite{1}, for $\gamma \neq 1$.

First note that \eqref{eq: condition roughness}, implies that $F$ is not constant in a neighborhood of $x$ (c.f.\ Lemma \ref{lemma: finite g if rough 2}). 
Fix $x \geq 0$ such that $x$ satisfies the conditions stated above. Define, for $a>0$, the process $T_n$ as follows:
$$
T_n(a)=\inf \left\{t \geq 0: U_n(t)-a t \text { is maximal }\right\}.
$$

The relation between $T_n$ and $\hat{V}_n$ is given by the switch relation (see \cite{2}):
\begin{align}\label{eq: switch relation}
T_n(a) \leq t \quad \Leftrightarrow \quad \hat{V}_n(t) \leq a .
\end{align}
From this relation (Lemma 3.2, \cite{2}), it follows that for each vanishing sequence $\left(\delta_n\right)$ of positive numbers, for $a\in \mathbb{R}$, for $\delta^*_n = \delta_n^{1/\gamma}$ ($\gamma$ as in \eqref{eq: condition roughness}), and all $n$ sufficiently large:
\begin{align*}
\delta_n^{-1}\left(\hat{V}_n(x)-V(x)\right) \leq a \: \: \: \: \: &\Leftrightarrow \: \: \: \: \: \delta_n^{-1}\left(T_n\left(a_0+\delta_n a\right)-x\right) \leq 0 \: \: \: \: \: \\
&\Leftrightarrow \: \: \: \: \: (\delta^*_n)^{-1}\left(T_n\left(a_0+\delta_n a\right)-x\right) \leq 0. 
\end{align*}
We aim to study the probability of the event on the left-hand side. We write $a_0$ for $V(x)$, which is fixed. The introduction of $\delta^*_n = \delta_n^{1/\gamma}$ is one of the main steps that differentiates the current proof from the one in \cite{1}. Now denote by:
$$
I_x := \begin{cases}{[0, \infty),} & \text { if } x=0 \\ \mathbb{R}, & \text { if } x>0\end{cases}
$$
\color{black} Define the localized process $\tilde{Z}_n$ in $C\left(I_x\right)$, the space of continuous functions on $I_x$, as follows (recall $a_0 = V(x)$):
$$
\tilde{Z}_n(t)=\delta_n^{-1} (\delta^*_n)^{-1} \left[U_n\left(x+\delta^*_n t\right)-U_n(x)-a_0 \delta^*_n t\right] .
$$
For each fixed $a \in \mathbb{R}$ we can, for all $n$ sufficiently large, write 
$$
\begin{aligned}
& (\delta^*_n)^{-1}\left(T_n\left(a_0+\delta_n a\right)-x\right) \\
&\quad = (\delta^*_n)^{-1} \big(\operatorname { i n f } \big\{x+t \geq 0: U_n(x+t)-a_0 x +\\
&\quad \quad \quad \quad \quad \quad \quad \quad \quad \quad \quad \quad -a_0 t-\delta_n a x-\delta_n a t \text { is maximal }\big\}-x\big) \\
&\quad = (\delta_n^*)^{-1} \inf \left\{t \geq-x: U_n(x+t)-U_n(x)-a_0 t-\delta_n a t \text { is maximal }\right\} \\
&\quad = \inf \left\{t \geq-(\delta^*_n)^{-1} x: U_n\left(x+\delta^*_n t\right)-U_n(x)-a_0 \delta^*_n t-\delta_n \delta^*_n a t \text { is maximal }\right\} \\
&\quad = \inf \left\{t \geq-(\delta^*_n)^{-1} x: \tilde{Z}_n(t)-a t \text { is maximal }\right\} ,
\end{aligned}
$$
where we use that the location of a maximum of a function is equivariant under translations and invariant under multiplication by a positive number and addition of a constant. Thus we analyze the asymptotic behavior of the process:
\begin{align*}
    \tilde{Z}_n(t)-a t =& \: \delta_n^{-1} ( \delta^*_n)^{-1} \left\{U_n\left(x+\delta^*_n t\right)-U_n(x)- U\left(x+\delta^*_n t\right)+U(x)\right\} \color{black}+ \numberthis \label{eq: rewrite Z tilde 1}\\ 
    &+ \delta_n^{-1} ( \delta^*_n)^{-1}  \left\{ U\left(x+\delta^*_n t\right)-U(x)-V(x)\delta^*_n t \right\} - at. \numberthis \label{eq: rewrite Z tilde}
\end{align*}
We will now analyze the behavior of the stochastic and deterministic terms \eqref{eq: rewrite Z tilde 1} and \eqref{eq: rewrite Z tilde} separately. Let us start with \eqref{eq: rewrite Z tilde}. The Taylor expansion of $U$ with integral form of the remainder yields:
\begin{align*}
    U(x + \delta^*_n t) - U(x) - V(x) \delta^*_n t &= \delta^*_n t \int_{0}^{1} V(x + u \delta^*_n t) - V(x) \, du \\
    &= -\frac{\delta^*_n t \pi}{2 m_0} \int_{0}^{1} F(x + u \delta^*_n t) - F(x) \, du.
\end{align*}
When we look at the limit of $\tilde{Z}_n(t) - at$ we multiply the above expression by $ \delta_n^{-1} (\delta^*_n)^{-1}= (\delta^*_n)^{-\gamma - 1}$. Therefore we have to look at the limit behavior of:
\begin{align}\label{eq: reasoning asympt holder}
 - (\delta^*_n)^{-\gamma - 1} \frac{\delta^*_n t \pi}{2 m_0} \int_{0}^{1} F(x + u \delta^*_n t) - F(x) \, du = -\frac{ |t|^{\gamma + 1} \pi}{2 m_0} (K + o(1)),
\end{align}
where we use assumption \eqref{eq: condition roughness}. 

Now we consider \eqref{eq: rewrite Z tilde 1}. Defining the processes $\left(Z_n\right)$ in $C\left(I_x\right)$ as
\begin{align*}
&Z_n(t) = \delta_n^{-1} ( \delta^*_n)^{-1} \left\{U_n\left(x+\delta^*_n t\right)-U_n(x)- U\left(x+\delta^*_n t\right)+U(x)\right\} \numberthis \label{eq: Z_n(1)}= \\
&2 \delta_n^{-1} (\delta^*_n)^{-1} \int\left\{\sqrt{z-x} \ind_{[x, \infty)}(z)-\sqrt{z-x-\delta^*_n t} \ind_{\left[x+\delta^*_n t, \infty\right)}(z)\right\} \, d\left(\mathbb{G}_n-G\right)(z)
\end{align*}
we will now establish convergence in distribution of the process $Z_n$. In order to do so, we need that $g(x) < \infty$. This follows directly from \eqref{eq: condition finite g at x}. Moreover, we need the asymptotic covariance structure of the process $Z_n$, which is given by:
\begin{align}\label{eq: covariance structure final}
\operatorname{Cov}\left(Z_n(s), Z_n(t)\right)=\frac{1}{2 \gamma} g(x) s t\left(1-\frac{\log \log n}{\log n}\right)+O\left(\frac{1}{\log n}\right) \text {, }
\end{align}
for each fixed $s, t \in I_x$, provided that $\delta_n=n^{-1 / 2} \sqrt{\log n}$.
Therefore, using \eqref{eq: covariance structure final} and the Lindeberg-Feller Central Limit theorem we have that: 
\begin{align}\label{eq: marginal convergence}
    Z_n(1) \rightsquigarrow N(0,g(x)/(2\gamma)),
\end{align} 
(the \hyperlink{proof of lemma 1}{proof} of \eqref{eq: covariance structure final} and the \hyperlink{proof of lemma 2}{proof} of \eqref{eq: marginal convergence} are provided in Appendix A). On the other hand, by Chebyshev's inequality and \eqref{eq: covariance structure final}, we have, $\forall \: \varepsilon >0$:
\begin{align*}
    \mathbb{P}\left( | Z_n (t) - t Z_n (1) | > \varepsilon \right) \leq &\frac{\mathbb{V} ( Z_n (t) - t Z_n (1))}{\varepsilon^2} \rightarrow 0.
\end{align*}
Therefore, the finite-dimensional distributions of $Z_n$ converge weakly to the finite-dimensional distributions of the process: $ Z(t) = tX,$ in $C(I_x)$, where $X \sim N(0,g(x)/(2\gamma))$. We now need a stochastic equicontinuity condition to apply Theorem 2.3 from \cite{12}. Take the class of functions for $\varepsilon >0, \: M>0$:  
\begin{align*}
\mathscr{F}_{n, \varepsilon}^M=\big\{\sqrt{z-x-\delta^*_n t} & \ind_{\left[x+\delta^*_n t, \infty\right)}(z)-\sqrt{z-x-\delta^*_n s} \ind_{\left[x+\delta^*_n s, \infty\right)}(z) :  \numberthis \label{eq: class of functions} \\
&: \: s, t \in I_x,|s-t|<\varepsilon, \max (|s|,|t|)<M\big\}.
\end{align*}
This class has square integrable natural envelope $H_\varepsilon $ that satisfies for all $\varepsilon>0$ small (cf.\ Lemma \ref{lemma: variance without assum. g}): $\int H^2_{\varepsilon} (z) g (z) dz \leq g(x + \delta^*_n s \wedge t) \left( \varepsilon \delta^*_n\right)^2 \log{\frac{1}{(\varepsilon \delta^*_n)^2}}$ (as proved in Lemma \ref{lemma: conv weak argmax} this class of functions is also VC).
Therefore, by using maximal inequality 3.1 (ii) from 
\cite{12} we obtain:
\begin{align*}
    &\limsup_{n \rightarrow \infty}\pb \left\{ \sup_{\mathscr{F}_{n, \varepsilon}^M} |Z_n(s) - Z_n(t)| > \eta \right\} \leq \limsup_{n \rightarrow \infty} \frac{\ex \sup_{\mathscr{F}_{n, \varepsilon}^M} |Z_n(s) - Z_n(t)|^2 }{\eta^2} \\
    &\leq \limsup_{n \rightarrow \infty} \, \frac{J(1)^2 g(x + \delta^*_n s \wedge t)}{(\delta^*_n)^2 \log{n}} (\varepsilon \delta^*_n)^2 \log{(\varepsilon \delta^*_n)^{-2}} =  J(1)^2 \gamma^{-1} \eta^{-2} g(x) \varepsilon^2 \stackrel{\varepsilon \rightarrow 0}{\longrightarrow} 0.
\end{align*}
where the reason why the uniform entropy integral is finite is as in the proof of Lemma \ref{lemma: conv weak argmax}.
The stochastic equicontinuity condition, together with the convergence in distribution of the finite-dimensional distributions, implies convergence in distribution of $Z_n$ to $Z$ in $(C(I_x),d)$ where $d$ is the metric of uniform convergence on compacta. Next, we need the stochastic boundedness of the sequence $T_n$ given below in \eqref{eq: stochastic boundedness seq argmax} (the \hyperlink{proof of lemma 3}{proof} of which is given in Appendix A). For $a_0 = V(x)$ and $n \rightarrow \infty$:
\begin{align}\label{eq: stochastic boundedness seq argmax}
T_n\left(a_0 + \delta_n a \right)=V^{-1}\left(a_0\right)+O_p\left(\delta^*_n\right).
\end{align}
One can now apply Theorem 2.7 from \cite{12} to the process $\tilde{Z}_n$. 
\begin{align}\label{eq: convergence argmax}
(\delta^*_n)^{-1}\left(T_n\left(a_0+\delta_n a\right)-x\right) & \rightsquigarrow \underset{t \in I_x}{\arg \max }\left\{t X-\frac{|t|^{\gamma + 1} \pi}{2 m_0} K  -a t \right\}. 
\end{align}
But it is immediate to see that, almost surely:
\begin{align}\label{eq: iff relation argmax}
    \underset{t \in I_x}{\arg \max }\left\{t X-\frac{|t|^{\gamma + 1} \pi}{2 m_0} K  -a t \right\} \leq 0 \: \: \: \: \iff \: \: \: \: X \leq a.
\end{align}
Hence, using the two relations above:
\begin{align*}
&\mathbb{P}\left\{\delta_n^{-1}\left(\hat{V}_n(x)-V(x)\right) \leq a\right\}=\mathbb{P}\left\{(\delta^*_n)^{-1}\left(T_n\left(a_0+\delta_n a\right)-x\right) \leq 0\right\} \stackrel{n \rightarrow \infty}{\longrightarrow} \\
& \quad \stackrel{\eqref{eq: convergence argmax}}{\longrightarrow} \mathbb{P} \left\{ \underset{t \in I_x}{\arg \max }\left\{t X-\frac{|t|^{\gamma + 1} \pi}{2 m_0} K  -a t \right\} \leq 0 \right\} \stackrel{\eqref{eq: iff relation argmax} }{=} \mathbb{P}\{ X \leq a\}.
\end{align*}
\end{proof}

\begin{example}\label{example: holder + local behavior}
As in \cite{10}, it is possible to conduct the analysis in accordance with the global smoothness properties of the underlying distribution function $F$. Consider, for arbitrary $\gamma=\alpha+\beta, \alpha=0,1, \ldots, 0<\beta \leq 1$, the Hölder class $C^{\gamma}$ of functions $f(y), y \geq 0$ having finite norm:
\begin{align}\label{eq: holder norm definition}
\|f\|_\gamma=\sup _{y \geq 0}|f(y)|+\sup _{x, y \geq 0} \frac{\left|f^{(\alpha)}(x)-f^{(\alpha)}(y)\right|}{|x-y|^\beta}.
\end{align}
Then assume $F \in C^{\gamma}$ for $\gamma > \frac{1}{2}$. Moreover, assume:
\begin{align}\label{eq: condition roughness 2}
    \lim_{t \rightarrow 0} \frac{|F(x+t)-F(x)|}{|t|^{\gamma}} = K > 0.
\end{align}
Note that it is possible to prove that if $F \in C^{\gamma}$ for $\gamma > 1/2$, then the density of the observations $g$ is bounded and continuous (see Lemma \ref{lemma: bounded density}). We show that the combination of \eqref{eq: holder norm definition} and \eqref{eq: condition roughness 2} implies \eqref{eq: condition roughness}.
Using the Dominated Convergence Theorem in combination with the assumption of the Hölder condition, and \eqref{eq: condition roughness 2}, we obtain as $\delta \rightarrow 0$:
\begin{align}\label{eq: reasoning asympt holder 2}
 \int_{0}^{1} \frac{F(x + u \delta) - F(x)}{ \operatorname{sgn}(\delta) |\delta|^{\gamma} u^{\gamma}} u^{\gamma} \, du \rightarrow K. 
\end{align}
\end{example}

Also for the case of Theorem \ref{thm: main thm iso 2}, we look first at the behavior of $\hat{V}_n(x)$. 

\begin{theorem}\label{thm: distrib flat F at x}
    Let $x \geq 0$, and $K_x :=[\ubar{x},\bar{x}] \subset [0,\infty)$ the biggest interval that contains $x$ on which $F$ is constant. Assume that \eqref{eq: finite first moment} holds true. Then:
    \begin{align*}
        \sqrt{n} \left( \hat{V}_n(x) - V(x) \right) \rightsquigarrow L_x,
    \end{align*}
    where, for any $a \in \mathbb{R}$:
    \begin{align*}
        \pb \left( L_x \leq a \right) = \pb \left( \argmax_{s \in K_x} \left\{ \mathbb{Z}_x(s) - a s\right\} \leq x \right). 
    \end{align*}
    for $\mathbb{Z}_x$ a centered continuous Gaussian Process with covariance structure given by:
    \begin{align}
        \mathrm{Cov} \left( \mathbb{Z}_x(t),\mathbb{Z}_x(s) \right) = \mathrm{Cov} \left( \sqrt{(Z_x)_{+}} - \sqrt{(Z_t)_{+}}, \sqrt{(Z_x)_{+}} - \sqrt{(Z_s)_{+}} \right).
    \end{align}
    where $Z_x := Z-x$ and analogously $Z_t$ and $Z_s$.
\end{theorem}
\begin{proof}
Using \eqref{eq: switch relation} we obtain for each fixed $a \in \mathbb{R}$, for all $n$ sufficiently large:
\begin{align*}
    &\sqrt{n} \left( \hat{V}_n (x) - V(x) \right) \leq a \quad \\
     \quad \Leftrightarrow \quad  &\inf \big\{s \geq 0: U_n(s) - \left(V(x)+ a/\sqrt{n} \right) s \: \text { is maximal} \big\} \leq x \quad \quad \\
    \quad \Leftrightarrow \quad  &\inf \big\{s \geq 0: (U_n (s) - U_n(x) - U (s) + U(x) ) + \\
    & \quad \quad \quad \quad + (U(s) - U(x) - V(x)(s-x) ) -  (a s)/\sqrt{n}  \text { is maximal} \big\} \leq x. 
\end{align*}
where we use that the location of a maximum of a function is equivariant under translations. Now using the invariance under multiplication by a positive number and addition of a constant, and defining:
\begin{align}
    &Z_n(s) := \sqrt{n} (U_n (s) - U_n(x) - U (s) + U(x) ) \quad \text{and} 
    \\ \quad &h_x(s) := U(x) - U(s) - V(x)(x-s).
\end{align}
we obtain that the above event can be rewritten as:
\begin{align*}
    \inf \big\{s \geq 0: Z_n(s) - a s - \sqrt{n} h_x(s)  \text { is maximal} \big\} \leq x,
\end{align*}
\begin{lemma}\label{lemma: conv weak argmax}
Let $\ell^{\infty}$ be endowed with uniform norm. Define:
$$ Z^*_n(s) := \sqrt{n} \left( U_n(s) - U(s) \right) $$
then:
\begin{align}\label{eq: convergence ell^infty Z 1}
    Z^*_n \rightsquigarrow \mathbb{Z} \quad \text{in} \: \: \ell^{\infty}[0,\infty),
\end{align}    
where $\mathbb{Z}$ is a Gaussian Process with a.s. continuous sample paths relative to the Euclidean metric. In particular:
\begin{align}\label{eq: norm infinity convergence}
    \|U_n - U \|_{\infty} = \sup_{s\geq 0} |U_n(s)- U(s)| \stackrel{\pb}{\rightarrow} 0.
\end{align}
\end{lemma}
The \hyperlink{proof of lemma weak conv}{proof} of this lemma is given in Appendix A. Consequently:
\begin{align}\label{eq: convergence ell^infty Z}
    Z_n = Z^*_n - Z^*_n(x) \rightsquigarrow \mathbb{Z} - \mathbb{Z}(x) =: \mathbb{Z}_x \quad \text{in} \: \: \ell^{\infty}[0,\infty),
\end{align}    
where also $\mathbb{Z}_x$ is a Gaussian Process with a.s. continuous sample paths relative to the Euclidean metric. By the above computations, we conclude that the limiting behavior of: 
$$\widetilde{Z}_n (s) :=  Z_n(s) - a s - \sqrt{n} h_x(s), $$ 
is determined by the process $\mathbb{W}$ which is defined by:
$$
\mathbb{W}(s) = \begin{cases}
        \mathbb{Z}_x(s) -as, & \text{for } s \in [\ubar{x},\bar{x}]\\
        -\infty, & \text{for } s \notin [\ubar{x},\bar{x}] 
        \end{cases} 
$$
This is because: $h_{x}(s) \geq 0$ for all $s \in (0 ,\ubar{x}) \cup (\bar{x},\infty)$, and therefore multiplied by $-\sqrt{n}$ will go to $-\infty$ in case the inequality is strict; however $h_x(s) = 0, \: \: \forall \: s \in [\ubar{x},\bar{x}]$ and thus on that interval the limiting behavior is completely determined by $\mathbb{Z}_x(s) -as$ (using the limiting behavior of $Z_n$ derived above). Let:
\begin{align*}
&\hat{s}_n := \inf \left\{ s \geq 0 \: : \:  Z_n(s) - a s - \sqrt{n} h_x(s)  \: \: \text{is maximal} \right\} \numberthis \label{eq: def hat{t}_n },\\
&\hat{s} := \argmax_{s \geq 0} \left\{\mathbb{W}(s) \right\} = \argmax_{s \in [\ubar{x},\bar{x}]} \left\{\mathbb{Z}_x(s) - as \right\}. \numberthis \label{eq: def hat{t} }
\end{align*}
To prove convergence in distribution we use the Portmanteau Lemma and show for every closed subset $F$ of $[0,\infty)$: $\limsup_{n \rightarrow \infty} \pb \left( \hat{s}_n \in F \right) \leq \pb \left( \hat{s} \in F \right)$. 

Let us start by choosing a sequence $\varepsilon_n \downarrow 0$ such that $\sqrt{n} h_x(\ubar{x}-\varepsilon_n) \rightarrow \infty$ and $\sqrt{n} h_x(\bar{x}+\varepsilon_n) \rightarrow \infty$. This exists since $\forall \: \varepsilon > 0$, we have $\sqrt{n} h_x(\ubar{x}-\varepsilon) \wedge h_x(\bar{x}+\varepsilon) \rightarrow \infty$. Moreover, $h_x(\bar{x}+\varepsilon)$ is an increasing function in $\varepsilon$ that has a zero in $\varepsilon =0$. The behavior is analogous for $h_x(\ubar{x}-\varepsilon)$. Therefore it is sufficient to construct a sequence of $\varepsilon_n$ such that $h_x(\ubar{x}-\varepsilon_n) \wedge h_x(\bar{x}+\varepsilon_n) = O(1/n^\alpha)$ for some $\alpha < 1/2$. 
For $K_n := [\ubar{x}-\varepsilon_n,\bar{x}+\varepsilon_n]$ and for $K := [\ubar{x},\bar{x}]$ (denoted by $K_x$ in the statement):
\begin{align*}
    &\limsup_{n \rightarrow \infty} \pb \left(\hat{s}_n \in F  \right) \leq \limsup_{n \rightarrow \infty} \pb \left(\hat{s}_n \in F \cap K_n \right) + \limsup_{n \rightarrow \infty} \pb \left(\hat{s}_n \in K^{c}_n  \right)  \\
    & \leq \underbrace{\limsup_{n \rightarrow \infty} \pb \left( \sup_{s \in F \cap K_n} \widetilde{Z}_n(s) \geq \sup_{s \in K_n} \widetilde{Z}_n(s) \right)}_{(1)} + \underbrace{\limsup_{n \rightarrow \infty} \pb \left(\hat{s}_n \in K^{c}_n  \right).}_{(2)} \numberthis \label{eq: two terms to show Portmanteau}
\end{align*}
We look at (2) in \eqref{eq: two terms to show Portmanteau} and show it goes to zero. To do this, we need relation \eqref{eq: conv K^c_n} and its \hyperlink{proof of lemma conv K^c_n}{proof} is given Appendix A. Thus:
\begin{align}\label{eq: conv K^c_n}
    \sup_{s \in K^c_n} \widetilde{Z}_n(s) = \sup_{s \in K^c_n}  \big\{ Z_n(s) - a s - \sqrt{n} h_x(s) \big\} \stackrel{\pb}{\rightarrow} - \infty,
\end{align}
implies that:
\begin{align*}
    \pb \left(\hat{s}_n \in K^{c}_n  \right) &\leq \pb \bigg( \sup_{s \in K^c_n} \widetilde{Z}_n(s) \geq \sup_{s \in K_n} \widetilde{Z}_n(s) \bigg) \\
    & \leq \pb \bigg( \sup_{s \in K^c_n} \widetilde{Z}_n(s) \geq \sup_{s \in K} \big\{ Z_n(s) -as \big\} \bigg) \rightarrow 0. \numberthis \label{eq: probabilities hat{t}_n in K^c}
\end{align*}
The last convergence is a consequence of \eqref{eq: conv K^c_n} and that $Z_n \rightsquigarrow \mathbb{Z}$ on $\ell^{\infty} (K)$, combined with the continuous mapping theorem (i.e.\ the term on the left-hand side of the inequality is $O_p(1)$).

Now we argue the behavior of (1). For that, we need an additional convergence in probability given in \eqref{eq: conv prob argmax}, and its \hyperlink{proof of lemma conv prob}{proof} is provided Appendix A.
\begin{align}\label{eq: conv prob argmax}
    \sup_{s \in F \cap K_n} \left\{Z_n(s) - \sqrt{n} h_x(s) -as\right\} - \sup_{s \in F \cap K} \left\{Z_n(s) -as\right\} \stackrel{\pb}{\longrightarrow} 0.
\end{align}
Using that:
\begin{align*}
    \sup_{s \in F \cap K} \left\{Z_n(s) -as\right\} \rightsquigarrow \sup_{s \in F \cap K} \left\{\mathbb{Z}_x(s) -as\right\},
\end{align*}
we conclude using Theorem 2.7 (iv) from \cite{6} and \eqref{eq: conv prob argmax}:
\begin{align*}
    \sup_{s \in F \cap K_n} \left\{Z_n(s) - \sqrt{n} h_x(s) -as\right\} \rightsquigarrow  \sup_{s \in F \cap K} \left\{\mathbb{Z}_x(s) -as\right\}.
\end{align*}
We derive the behavior of (1) in \eqref{eq: two terms to show Portmanteau}. This term can be upper bounded by:
\begin{align*}
    &\limsup \pb \left( \sup_{s \in F \cap K_n} \left\{Z_n(s) - \sqrt{n} h_x(s) -as\right\} \geq \sup_{s \in K} \left\{ Z_n(s) -as \right\} \right) \\
    &\leq \pb \left( \sup_{F \cap K} \left\{\mathbb{Z}_x(s) -as\right\} \geq \sup_{s \in K} \left\{\mathbb{Z}_x(s) -as\right\} \right) \leq \underbrace{\pb (\hat{s} \in K^c)}_{=0} + \pb \left( \hat{s} \in F \right) \numberthis \label{eq: first term portmanteau}
\end{align*}
where we used the above derivations, together with Lemma \ref{lemma: conv weak argmax} and Theorem 2.7 (v) from \cite{6}. Therefore by \eqref{eq: two terms to show Portmanteau} combined with \eqref{eq: probabilities hat{t}_n in K^c}, \eqref{eq: first term portmanteau} and Portmanteau Lemma we obtain the desired convergence in distribution, which can be restated:
\begin{align*}
    & \pb \left(\sqrt{n} \big( \hat{V}_n (x) - V(x) \big) \leq a \right) =  \pb \Big( \hspace{-0.05cm}  \inf_{s \geq 0} \big\{ Z_n(s) - a s - \sqrt{n} h_x(s)  \text { is maximal} \big\} \hspace{-0.05cm} \leq x \Big) \\
    & \stackrel{n \rightarrow \infty}{\longrightarrow}  \pb \bigg( \argmax_{s \in K} \left\{ \mathbb{Z}_x(s) - a s\right\} \leq x \bigg).
\end{align*}
Note that the last expression defines the cdf of a random variable $L_x$, which is the right derivative of the least concave majorant of the process $\mathbb{Z}_x$ evaluated at $x$. Indeed, by the switch relation, we have 
$\pb \left( \argmax_{s \in K} \left\{ \mathbb{Z}_x(s) - a s\right\} \leq x \right) = \pb \left( L_x \leq a \right).$
\end{proof}

Note that the \hyperlink{proof of main thm}{proof} of Theorem \ref{thm: main thm iso 1} (included in Appendix A) is a consequence of Slutsky's Lemma, Theorem \ref{thm: normality Holder} and \ref{thm: distrib flat F at x} and the fact that $1 - F(x) = V(x)/V(0)$. The proof of Theorem \ref{thm: main thm iso 2} follows the same steps as the proof of Theorem \ref{thm: distrib flat F at x}.

\begin{remark}
Let $x>0$, and suppose that $F$ satisfies condition \eqref{eq: condition roughness} both at $x$ and at 0 for a common $\gamma > \frac{1}{2}$ and \eqref{eq: finite first moment} holds true. Then, by Theorem \ref{thm: main thm iso 1}:
\begin{align}
& \sqrt{\frac{n}{\log n}}\left( \hat{F}_n(x) -F(x)\right) \rightsquigarrow N\left(0, \frac{g(x) V(0)^2+g(0) V(x)^2}{2 \gamma V(0)^4}\right).
\end{align}
The asymptotic variance obtained here coincides with the one obtained in obtained in \cite{10} where it is given a unique $\gamma$ and the one in \cite{1} with $\gamma = 1$. 
\end{remark}

\subsubsection*{Simulation study limiting process Theorem \ref{thm: distrib flat F at x}}

\begin{figure}[ht]
  \centering
  \includegraphics[width=6cm]{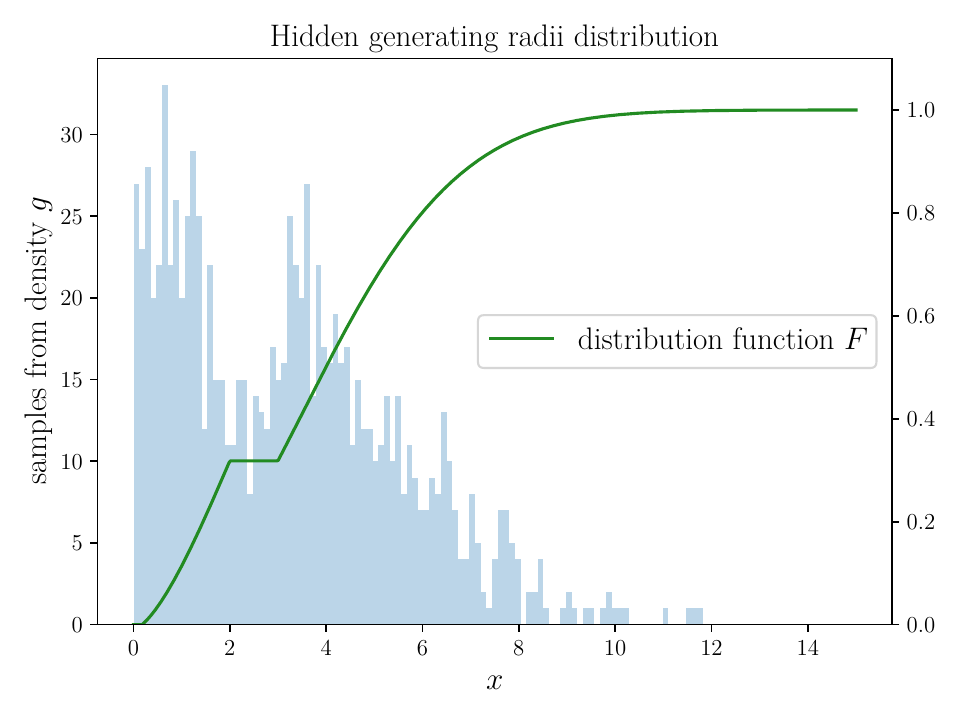}
  \includegraphics[width=5.75cm]{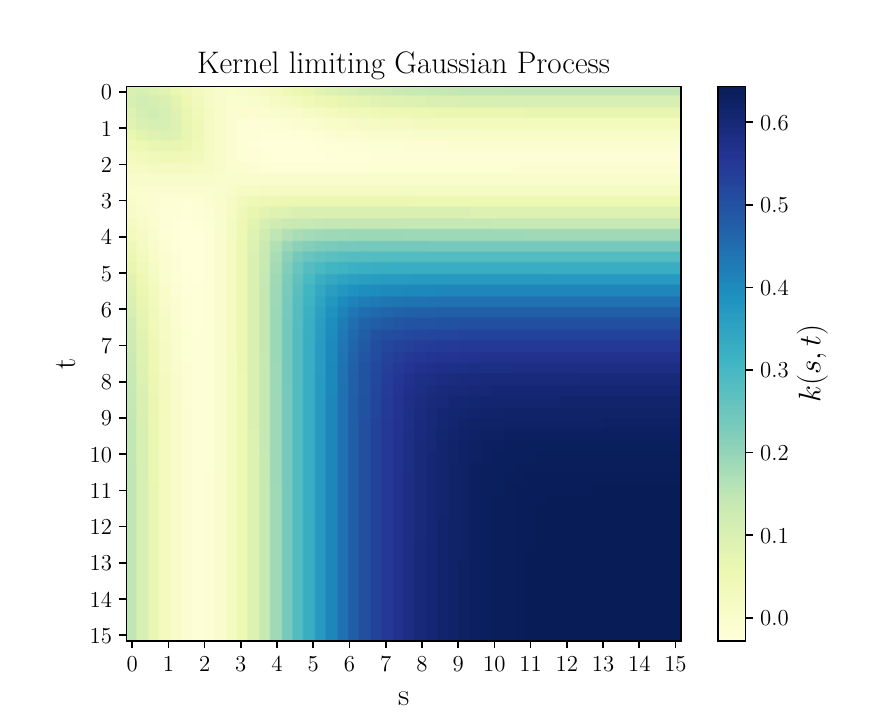}
  \setlength{\belowcaptionskip}{-4pt}
  \setlength{\abovecaptionskip}{0pt}
  \caption{Left panel: underlying cdf $F$ with samples from $g$. Right panel: representation of covariance structure of the limiting Gaussian Process $\mathbb{Z}_x$.}\label{fig: 1}
\end{figure}

In this subsection, we present some simulations about the limiting process of Theorem \ref{thm: distrib flat F at x}, as we believe it is insightful to observe its behavior.
First, we consider an underlying distribution function $F$ of the squared radii as displayed in the left panel of Figure \ref{fig: 1}, which is constant on $[2,3]$; observe in light blue a histogram of $1000$ gathered samples from the associated density of the observations $g$ (note that we obtain observations also on the interval where $F$ is constant).

\begin{figure}[ht]
  \centering
  \includegraphics[width=6.2cm]{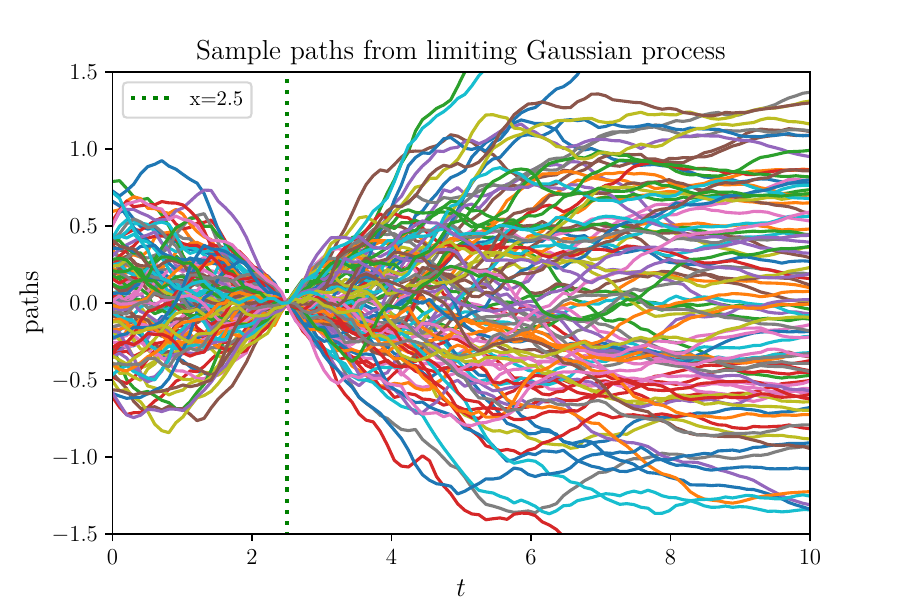}\hspace{-0.5cm}
  \includegraphics[width=6.2cm]{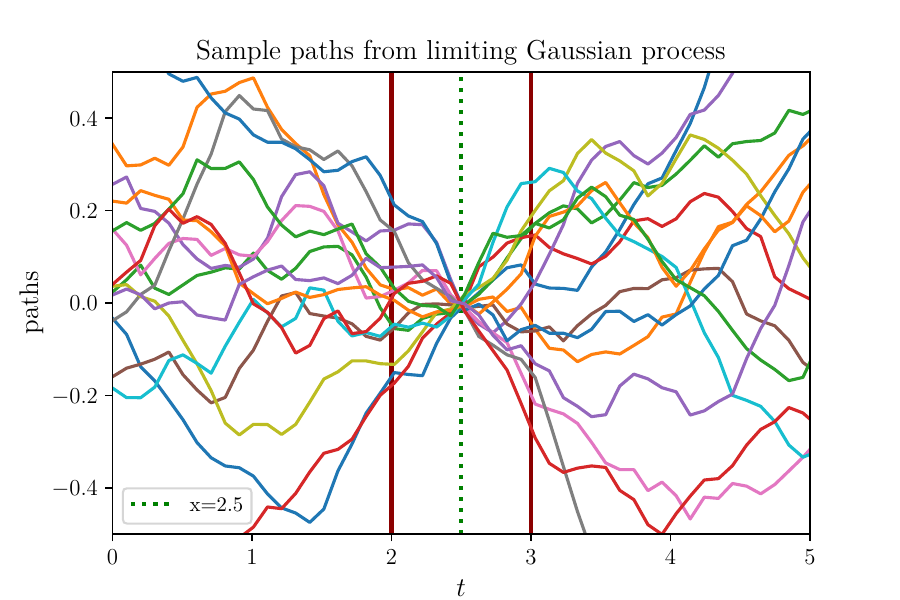}
  \setlength{\belowcaptionskip}{-4pt}
  \setlength{\abovecaptionskip}{-4pt}
  \caption{Left panel: $200$ samples from limiting Gaussian Process $\mathbb{Z}_x$. Right panel: zoom in on $15$ samples paths, putting into evidence the interval $[2,3]$.}\label{fig: 2}
\end{figure}

In this setting, we show some simulations of the distribution of the limiting random variable of Theorem \ref{thm: distrib flat F at x}, above denoted by $L_x$. To do that, we compute the kernel given in \eqref{eq: covariance gauss process} of the limiting Gaussian process $\mathbb{Z}_x$. For the $F$ given in the left panel of Figure \ref{fig: 1}, the corresponding kernel is given in the right panel of Figure \ref{fig: 1}.
Subsequently, we simulate $200$ sample paths from the limiting Gaussian Process $\mathbb{Z}_x$ with the given kernel. These sample paths are shown in Figure \ref{fig: 2}. In the right panel of Figure \ref{fig: 2} we zoom in on the interval $[2,3]$ of interest and look at the behavior of $15$ selected sample paths. As we would expect, all sample paths pass through $(2.5,0)$.

\begin{figure}
   \vspace{-0.1cm} 
  \centering
  \includegraphics[width=6cm]{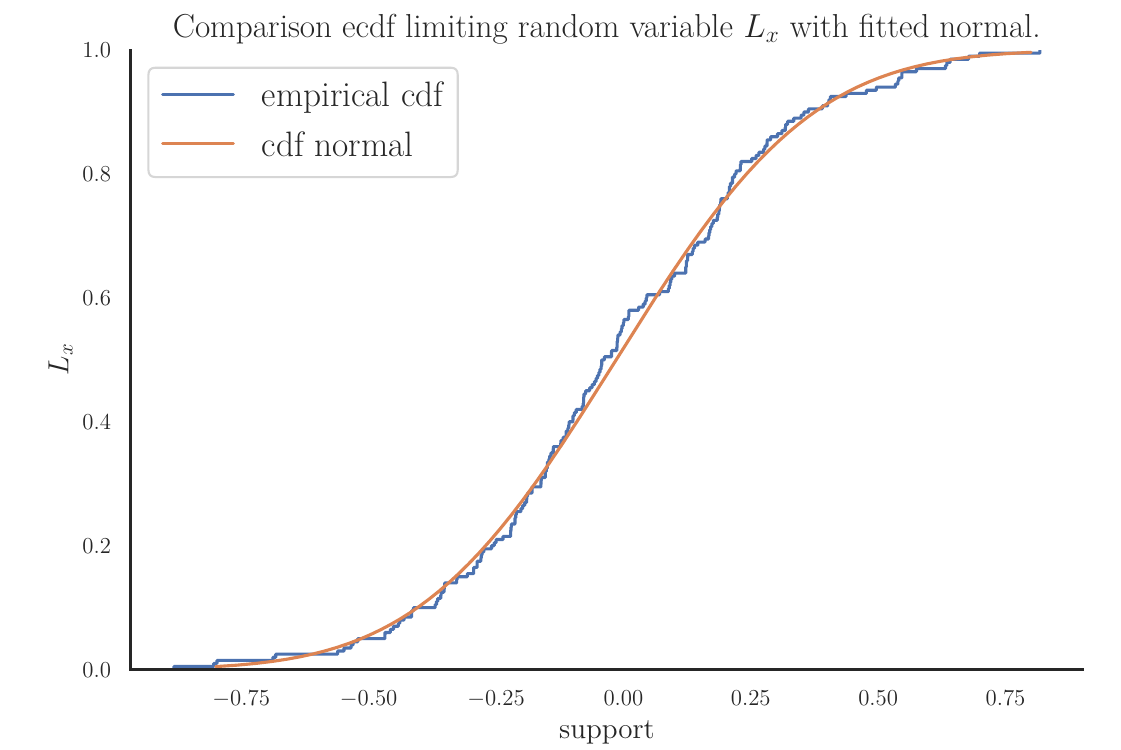}\hspace{-0.5cm}
  \includegraphics[width=6.5cm]{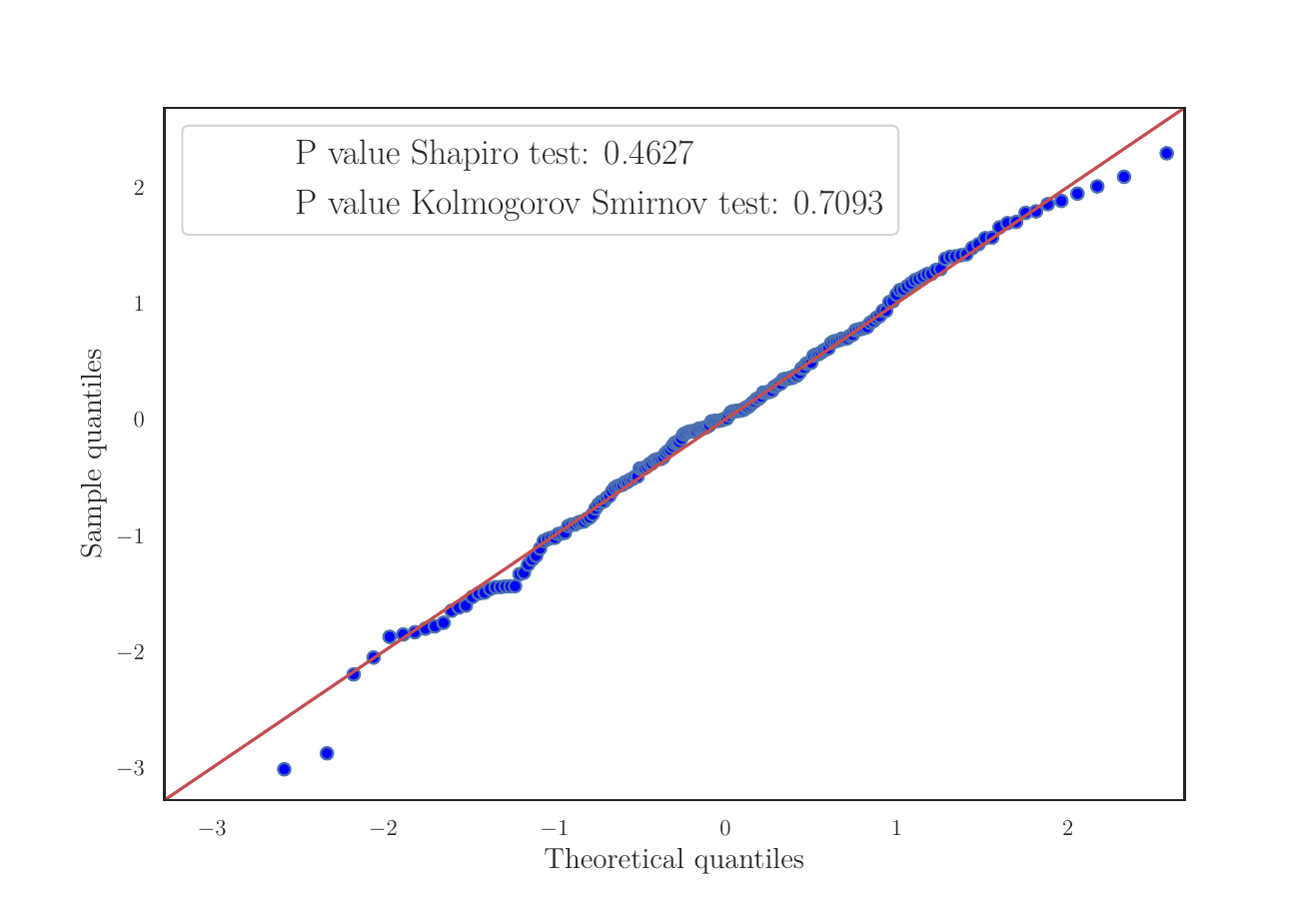}
  \setlength{\belowcaptionskip}{-6pt}
  \setlength{\abovecaptionskip}{-4pt}
  \caption{Left panel: ecdf of limiting random variable $L_x$ compared to a suitably fitted normal. Right panel: qq-plot of $L_x$ compared to the fitted normal and p-values of the standard tests.}\label{fig: 3}
\end{figure}

In order to simulate the distribution of $L_x$, we use the fact that such random variable, by the proof of the theorem, is the right derivative of the least concave majorant of the Gaussian Process $\mathbb{Z}_x$ evaluated at $x = 2.5$.  With all the gathered samples from $L_x$ we looked at its estimated distribution and below we show some tests about the normality of the distribution of $L_x$. With the gathered samples, we see in Figure \ref{fig: 3} that we could not reject the null hypothesis that $L_x$ is normally distributed.

\section{Local Asymptotic Minimax Theorem}

In this section, we show that the isotonic estimator (IIE) is asymptotically efficient, in the sense that its asymptotic variance is the smallest attainable. For this reason, we establish a statement about a lower bound of the minimax risk for any estimator sequence. We do this by constructing a parametric perturbation along the least favorable direction in the proximity of the true distribution function $F$. 

Let $x > 0$ be the point at which we want to estimate $F$. The existence of a least favorable path is a non-trivial matter. We will be guided by intuition in the construction of this path. We want to estimate is: $F(x) = 1 - V(x)/V(0)$. The local behavior at $x$ and $0$ of the kernels in the definitions of $V(x)$ and $V(0)$ is what essentially characterizes the estimation in Wicksell's problem and leads to the unusual rate $n^{-1/2}\sqrt{\log{n}}$. Thus, it should be enough to introduce a local perturbation of $F$ around $x$ and $0$. However, such perturbation must be tailored according to the local level of smoothness of the true $F$ near $x$ and $0$. Therefore, we need to identify such a level, and this is done in the following condition \eqref{eq: condition roughness infimum}.

The level of locality of the condition must be tied to the smoothness itself, because if the true $F$ is very rough around $x$ or $0$, then we expect a more stringent level of locality (and vice versa). Therefore we define the following intervals close to $p \in \{0,x\}$ for constants $\gamma_0 > \frac{1}{2}, \gamma_x > \frac{1}{2}$ to be determined later:
\begin{gather*}
\prescript{r}{}{\mathrm{K}}^{0}_n := \left[ n^{-\frac{1}{2\gamma_{\scaleto{0}{2.5pt}}}}, \eta_n  \right],\\
\prescript{l}{}{\mathrm{K}}^{x}_n := \left[ x - \eta_n, x - n^{-\frac{1}{2\gamma_{\scaleto{x}{2.5pt}}}} \right], \: \: \: \prescript{r}{}{\mathrm{K}}^{x}_n := \left[ x + n^{-\frac{1}{2\gamma_{\scaleto{x}{2.2pt}}}},x + \eta_n  \right],
\end{gather*}
where $\eta_n \downarrow 0$, $0 < \eta_n \leq 1, \: \forall \, n$ and $\exists \: \alpha \in (0,1)$ such that $\eta_n \gtrsim \frac{1}{\log^{\alpha}(n)}$. Now, let $\gamma_0, \gamma_x$ be the smallest constants greater than $1/2$, such that the true $F$, for $\mathrm{K}^p_n \in \{ \prescript{r}{}{\mathrm{K}}^{x}_n, \prescript{l}{}{\mathrm{K}}^{x}_n, \prescript{r}{}{\mathrm{K}}^{0}_n\}$ and therefore $p \in \{0,x\}$ satisfies:
\begin{align*}
    &\liminf_{n \rightarrow \infty } \inf_{\substack{x_1, x_2 \in \mathrm{K}^p_n }}\sqrt{n\log{n}}  \frac{|F(x_2)-F(x_1)|}{\Big|\log{\left\{\frac{p-x_1}{p-x_2} \right\}}\Big|} = \infty. \numberthis \label{eq: condition roughness infimum} 
\end{align*}

With condition \eqref{eq: condition roughness infimum}, our path can be shown to be well-defined and it is the least restrictive we can require. As we naturally would expect, the $\gamma_0, \gamma_x$ introduced above contain information about the local roughness of $F$ near $x$ and $0$. We illustrate this with our previous example \ref{example: cdf holder} and an additional example.
\begin{example}
Let us consider the same class of functions as in Example \ref{example: cdf holder}. These functions can be thought to be tightly Hölder smooth at $x$ and satisfy \eqref{eq: condition roughness infimum}. \\
First notice that in this case, we have for all $n$ big enough so that $\prescript{r}{}{\mathrm{K}}^{x}_n \subseteq U_x$:
\begin{align*}
    \inf_{\substack{x_1, x_2 \in \prescript{r}{}{\mathrm{K}}^{x}_n}} &\frac{F(x_2) - F(x_1)}{\log{\{\frac{x-x_1}{x-x_2} \}}} = \inf_{\substack{x_1, x_2 \in \prescript{r}{}{\mathrm{K}}^{x}_n }} K \frac{|x - x_1|^{\gamma_x} - |x - x_2|^{\gamma_x}}{\log{(x-x_1)} - \log{(x-x_2)}} \\
    &=\inf_{\substack{x_1, x_2 \in \prescript{r}{}{\mathrm{K}}^{x}_n }} K (x_2 - x)^{\gamma_x} \frac{\left( (1 - \frac{x_2 - x_1}{x_2 - x})^{\gamma_x} -1\right)}{\log{\left(1 -  \frac{x_2 - x_1}{x_2 - x} \right)}}.
\end{align*}
where, $x_1 < x_2$, both $x_1, x_2 \in \prescript{r}{}{\mathrm{K}}^{x}_n$ and thus $0<\frac{x_2 - x_1}{x_2 - x} <1$ going to zero as $n \rightarrow \infty$. The above expression, for all $n$ big enough behaves like $(1/n^{1/(2\gamma_x)})^{\gamma_x} \sim 1/\sqrt{n}$. Multiplied by $\sqrt{n\log{n}}$ \eqref{eq: condition roughness infimum} diverges to $\infty$.
\end{example}
\begin{example}
Instead of condition \eqref{eq: condition roughness infimum}, a stronger one, which is similar to the conditions of \cite{10} (cf. (24) and (25)), suffices: 
$\liminf_{n \rightarrow \infty } \inf \frac{|F(x_2)-F(x_1)|}{|x_2 - x_1|^{\gamma_{\scaleto{x}{2.2pt}}} } > 0$, where the infimum is as in \eqref{eq: condition roughness infimum}. We show that it implies \eqref{eq: condition roughness infimum} on $\prescript{r}{}{\mathrm{K}}^{x}_n$:
\begin{align*}
    \inf_{\substack{x_1, x_2 \in \prescript{r}{}{\mathrm{K}}^{x}_n}} &\frac{F(x_2) - F(x_1)}{\log{\{\frac{x-x_1}{x-x_2} \}}} = \inf_{\substack{x_1, x_2 \in \prescript{r}{}{\mathrm{K}}^{x}_n }} \frac{|F(x_2)-F(x_1)|}{|x_2 - x_1|^{\gamma_{\scaleto{x}{2.2pt}}} } \frac{|x_2 - x_1|^{\gamma_{\scaleto{x}{2.2pt}}}}{\log{\left( 1 - \frac{x_2 - x_1}{x_2-x} \right)}}
\end{align*}
where, $x_1 < x_2$, both $x_1, x_2 \in \prescript{r}{}{\mathrm{K}}^{x}_n$ and thus $0<\frac{x_2 - x_1}{x_2 - x} <1$ going to zero as $n \rightarrow \infty$. The above expression, for all $n$ big enough behaves like in the previous example. Multiplied by $\sqrt{n\log{n}}$ \eqref{eq: condition roughness infimum} diverges to $\infty$.
\end{example}

The second fundamental assumption to obtain our minimax lower bound is that $g$ is continuous in a neighborhood of $x$ and in a right neighborhood of $0$. We will assume this all along. A sufficient condition to obtain this is that $F$ is locally Hölder continuous around $x$ and $0$ of degrees greater than $\frac{1}{2}$ (see Lemma \ref{lemma: bounded density} in Appendix A). We also note that from \eqref{eq: condition roughness infimum} we have $F(x)<1$.

We will now proceed with the construction of the path on the space of distribution functions and verify that it is well defined. The construction of our path is based on ideas from \cite{10}. Condition \eqref{eq: condition roughness infimum} gives us the appropriate $\gamma_0, \gamma_x$ that we use for the definition of the perturbation.
 For $u \geq 0$ define: 
\begin{align}\label{eq: path}
    F_{h_n}(u) = \frac{1}{D_{h_n}} \left( F(u) + h_{1,n} \int_{0}^{u} \chi_{{\scaleto{1,n}{4.5pt}}}(v) \, dv + h_{2,n} \int_{0}^{u} \chi_{{\scaleto{2,n}{4.5pt}}}(v-x) \, dv \right),
\end{align}
where for $\gamma_0 > \frac{1}{2}, \: \gamma_x > \frac{1}{2}$, using that $\int_{0}^{\infty} \chi_{{\scaleto{2,n}{4.5pt}}}(v-x) \, dv = 0$:
\begin{align*}
    &\chi_{{\scaleto{1,n}{4.5pt}}}(v)= (v)^{-1} \ind_{\left\{ n^{-1/(2\gamma_{\scaleto{0}{2.2pt}})} \leq v \leq \eta_n \right\}},  \\
    &\chi_{{\scaleto{2,n}{4.5pt}}}(v-x)= (v-x)^{-1} \ind_{\left\{ n^{-1/(2\gamma_{\scaleto{x}{2.2pt}})} \leq|v-x| \leq \eta_n \right\}}, \\
    &D_{h_n} = 1 + h_{1,n} \int_{0}^{\infty} \chi_{{\scaleto{1,n}{4.5pt}}}(v) \, dv + h_{2,n} \int_{0}^{\infty} \chi_{{\scaleto{2,n}{4.5pt}}}(v-x) \, dv \\
    &\quad \quad = 1 + h_{1,n} \int_{0}^{\infty} \chi_{{\scaleto{1,n}{4.5pt}}}(v) \, dv. 
\end{align*}
We will see that the truncation levels of the function $1/(v-x)$ proposed here above work, but they may not be unique. The function $\chi_{{\scaleto{i,n}{4.5pt}}}$, with similar truncation levels, was introduced for the first time in \cite{10}, however there it was not used directly for the perturbation as done here and the assumption $\gamma_0 = \gamma_x$ was used.

We will check that the perturbation $ F_{h_n}$ is well defined for: 
\begin{align}\label{eq: h_i}
    h_{i,n}= \frac{h_i}{\sqrt{n \log{n}}}, \quad h_i \in \mathbb{R}, \quad i=1,2. 
\end{align}

\begin{figure}
    \vspace{-0.1cm}
  \centering
  \includegraphics[width = 8cm]{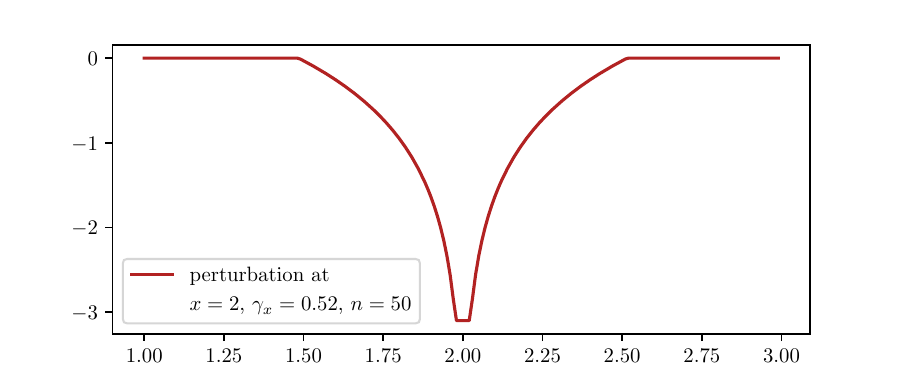}
  \setlength{\belowcaptionskip}{-4pt}
  \caption{Visualization of $\int_{0}^{u} \chi_{{\scaleto{2,n}{4.5pt}}}(v-x)dv$ for $x=2$, $n=50$ and $\gamma_x = 0.52$, $\eta_n = 1/\sqrt{\log{n}}$.} \label{fig: perturbation}
\end{figure}

First, we compute the case of $\chi_{{\scaleto{2,n}{4.5pt}}}$, to show how the perturbation looks like:
\begin{align*}
    &\int_{0}^{u} \chi_{{\scaleto{2,n}{4.5pt}}} (v-x) \, dv = \log{\left(\eta^{-1}_n (x-u) \right)}\ind_{\big\{ x - \eta_n \leq u \leq x - n^{-\frac{1}{2 \gamma_{\scaleto{x}{2.2pt}}}} \big\}} +\\
    &\quad + \log{\left( \eta^{-1}_n \, n^{-\frac{1}{2 \gamma_{\scaleto{x}{2.2pt}}}} \right)} \bigg( \ind_{\big\{ u > x - n^{-\frac{1}{2 \gamma_{\scaleto{x}{2.2pt}}}} \big\}} - \ind_{\big\{ u > x + \eta_n \big\}} \bigg)+\\
    &\quad + \log{\left(\eta^{-1}_n(u-x) \right)} \ind_{\big\{x + n^{-\frac{1}{2 \gamma_{\scaleto{x}{2.2pt}}}} \leq u \leq x+\eta_n  \big\}}. \numberthis \label{eq: written out perturbation}
\end{align*}
Figure \ref{fig: perturbation} visualizes the perturbation $\int_{0}^{u} \chi_{{\scaleto{2,n}{4.5pt}}}(v-x)dv$ in \eqref{eq: written out perturbation}. It is immediate from \eqref{eq: written out perturbation} that this perturbation is continuous and therefore the summation with $F$ results in a function $F_{h_n}$ that is right-continuous. Moreover, thanks to the normalization $D_{h_n}$, $\lim_{u \rightarrow + \infty} F_{h_n} (u) =1$ and, from the definition of $F_{h_n}$, note that $\lim_{u \rightarrow 0} F_{h_n} (u) =0$.
What is not evident is if $F_{h_n}$ is non-decreasing. Intuitively, for this to be true, we need $F$ to be "increasing enough" to compensate for the negative drop introduced by the perturbation. Condition \eqref{eq: condition roughness infimum} ensures that this is the case.  
\begin{lemma}\label{lemma: non decreasing F_h}
    For all $n$ sufficiently large, $F_{h_n}$ is non decreasing.
\end{lemma}
The \hyperlink{proof of lemma 5}{proof} of this lemma is provided in Appendix B.  We start by stating the most interesting result of this section which is the so-called LAM (locally asymptotically minimax) theorem and it entails the minimax risk for any estimator sequence. The proof of this theorem requires a few additional technical statements (the LAN expansion and the expression of the derivative of the functional of interest), which will be provided right after.
\begin{theorem}[LAM]\label{thm: LAM}
    Let $\gamma_0, \gamma_x$ be the smallest constants above $\frac{1}{2}$ such that $F$ satisfies \eqref{eq: condition roughness infimum}. Let then $J$ and $g$ as in Theorem \ref{thm: LAN}, $\dot{\psi}$ as in Lemma \ref{lemma: hadamard derivative}, and $\ell$ be a bowl-shaped loss function, i.e., $\ell$ is non-negative, symmetric and subconvex. Then for every estimator sequence $\left( F_n(x)\right)_{n \in \mathbb{N}}$:
    \begin{align}\label{eq: LAM}
    \sup_{\substack{I \subset \mathbb{R}^2 \\ I \: \textup{finite}}} \liminf_{n \rightarrow \infty} \sup _{h \in I} \ex_{h} \ell\left(\sqrt{\frac{n}{\log{n}}}\left(F_n(x)- F_{h_n}(x)\right)\right) \geq \ex \ell(L) .
    \end{align}
The first supremum is taken over all finite subsets $I$ of $\mathbb{R}^2$, the expectation $\ex_h$ is under $G_{F_{h_n}}$, the probability measures associated to the perturbed densities $g_n$ in \eqref{eq: g_n}, and $L \sim N (0, \dot{\psi}^{\top} J^{-1} \dot{\psi})$;
\begin{align*}
    {\dot{\psi}}^{\top} J^{-1} \dot{\psi} = \frac{4m^2_0}{\pi^2} \left(\frac{g(x)}{2\gamma_x} + (1-F(x))^2 \frac{g(0)}{2\gamma_0} \right).
\end{align*}
\end{theorem}
\begin{proof}
The proof is an immediate application of Theorem 3.11.5 from \cite{5}, using Theorem \ref{thm: LAN} and Lemma \ref{lemma: hadamard derivative} given below.
\end{proof}

Now we proceed by proving the technical results needed for the proof of Theorem \ref{thm: LAM}. First, we show that given the path \eqref{eq: path}, we can obtain a LAN expansion (see \cite{6}, Chapter 7).
In the statement of the LAN Theorem, the perturbed density of the observations $Z_i$ appears. Therefore, using \eqref{eq: path}, we derive it before proceeding:
\begin{align*}
    &g_n (z) = \frac{1}{2m_{h_n}} \int_{z}^{\infty} \frac{dF_{h_n} (v)}{\sqrt{v-z}}  \numberthis \label{eq: g_n} \\
    &= \frac{1}{2m_{h_n} D_{h_n}} \left\{ \int_{z}^{\infty} \frac{dF(v)}{\sqrt{v-z}} + h_{1,n} \int_{z}^{\infty} \frac{\chi_{{\scaleto{1,n}{4.5pt}}}(v)}{\sqrt{v-z}} \, dv + h_{2,n} \int_{z}^{\infty} \frac{\chi_{{\scaleto{2,n}{4.5pt}}}(v-x)}{\sqrt{v-z}} \, dv \right\}.
\end{align*}
where $m_{h_n} = \int_{0}^{\infty} \sqrt{s} \, dF_{h_n} (s)$, analogously to $m_0$.

\begin{theorem}[LAN expansion]\label{thm: LAN}
Let $\gamma_0, \gamma_x$ be constants above $\frac{1}{2}$ such that $F$ satisfies \eqref{eq: condition roughness infimum}. Let $g$ be continuous in a neighborhood of $x$ and in a right-neighborhood of $0$. Then for $h^{\top} = [h_1,h_2] \in \mathbb{R}^2$, $h_n = \frac{h}{\sqrt{n\log{n}}}$, $F_{h_n}$ as in \eqref{eq: path} and accordingly $g_n$ in \eqref{eq: g_n}:
\begin{align*}
    \sum_{i=1}^{n} \log{\frac{g_n (Z_i)}{g(Z_i)}} = h^{\top} \Delta_n - \frac{1}{2} h^{\top} J h + o_p (1),
\end{align*}
where the $o_p(1)$ is under $Z_1,\cdots,Z_n \stackrel{\scaleto{i.i.d.}{3.2pt}}{\sim} g$, for $J =  \frac{\pi^2 }{8 m^2_0 } \operatorname{diag}\left(\frac{1}{\gamma_{\scaleto{0}{2.5pt}} g(0)}, \frac{1}{ \gamma_{\scaleto{x}{2.2pt}} g(x)}\right)$ and $\Delta_n \rightsquigarrow N\left( 0, J \right)$, where:
\begin{align*}
    \Delta_n = \begin{bmatrix} \sum_{i=1}^n \frac{1}{ m_{h_n} D_{h_n} \sqrt{n\log{n}} } \left( \frac{1}{2g(Z_i)} \int_{Z_i}^{\infty} \frac{\chi_{{\scaleto{1,n}{3.5pt}}}(v)}{\sqrt{v-Z_i}} \, dv - \int_{0}^{\infty} \sqrt{v} \chi_{{\scaleto{1,n}{4.5pt}}}(v) \, dv \right)\\ \sum_{i=1}^n \frac{1}{ m_{h_n} D_{h_n} \sqrt{n\log{n}}} \left( \frac{1}{2g(Z_i)} \int_{Z_i}^{\infty} \frac{\chi_{{\scaleto{2,n}{3.5pt}}} (v-x)}{\sqrt{v-Z_i}} \, dv - \int_{0}^{\infty} \sqrt{v} \chi_{{\scaleto{2,n}{4.5pt}}}(v - x) \, dv \right) \end{bmatrix}.
\end{align*}
\end{theorem}
\begin{proof}

The structure of the proof of this result is as follows. First, we give the expansions of the functions $\int_z^{\infty} \frac{\chi_{{\scaleto{2,n}{3.5pt}}}(v-x) }{\sqrt{v-z}} \, d v, \: \int_{z}^{\infty} \frac{\chi_{{\scaleto{1,n}{3.5pt}}}(v)}{\sqrt{v-z}} \, dv$ that will highlight the leading terms that contribute to the limiting behavior. 

Next, we use the obtained approximations to write out the Taylor approximation of the logarithm and prove some additional relations that will eventually give the LAN. Finally, we show that the remaining terms are negligible. First of all, using the definition of $F_{h_n}$ in \eqref{eq: path}, we obtain:
\begin{align*}
    \frac{g_{n} (z)}{g (z)} = \frac{1}{m_{h_n} D_{h_n}} \left\{ m_0+ \frac{h_{1,n}}{2g(z)} \int_{z}^{\infty} \frac{\chi_{{\scaleto{1,n}{4.5pt}}}(v)}{\sqrt{v-z}} \, dv + \frac{h_{2,n}}{2g(z)} \int_{z}^{\infty} \frac{\chi_{{\scaleto{2,n}{4.5pt}}} (v-x)}{\sqrt{v-z}} \, dv \right\}.
\end{align*}
Then using \eqref{eq: approximation xi} (see Lemma \ref{lemma: expansion GL} proved in Appendix B): 
\begin{align*}
& \int_z^{\infty} \frac{\chi_{{\scaleto{2,n}{4.5pt}}}(v-x) }{\sqrt{v-z}} \, d v = \left\{ \frac{\pi}{\sqrt{z-x}} + O\left(\eta_n^{-\frac{1}{2}}\right)  \right\} \ind_{\left\{ n^{-\frac{1}{2 \gamma_{\scaleto{x}{2.2pt}}}} \leq z-x \leq \eta_n \right\}} + \numberthis \label{eq: only term that matters in expansion}\\
& +O\left(n^{\frac{1}{4\gamma_{\scaleto{x}{2.2pt}}}}\right) \ind_{\left\{ -n^{-\frac{1}{2 \gamma_{\scaleto{x}{2.2pt}}}}  \leq z-x \leq n^{-\frac{1}{2 \gamma_{\scaleto{x}{2.2pt}}}} \right\}} +  O\left(\eta_n^{-\frac{1}{2}}\right) \ind_{\left\{ z-x \leq - \eta_n \right\}} + \\
& + \left\{ O\left(n^{-\frac{1}{2\gamma_{\scaleto{x}{2.2pt}}}}(x-z)^{-\frac{3}{2} }\right)+O\left(\eta_n^{-\frac{1}{2}}\right) \right\} \ind_{\left\{ - \eta_n \leq z-x \leq - n^{-\frac{1}{2 \gamma_{\scaleto{x}{2.2pt}}}}\right\}}. \numberthis \label{eq: expansion zeta}
\end{align*}
The proof of Lemma \ref{lemma: expansion GL} has the following relation as an important ingredient:
\begin{align}\label{eq: abel integral pi}
 \int_u^{\infty} \frac{1}{t \sqrt{t-u}} \, d t=\pi \ind_{\{u \geq 0\}} u^{-1 / 2}, 
\end{align}
which follows from: $\int_{1}^{\infty} \frac{1}{t \sqrt{t-1}} dt = \pi$.
We will see that asymptotically the only term that matters is the first one in \eqref{eq: only term that matters in expansion}. It is immediate to see that $\sup_u | \int_{0}^{u} \chi_{{\scaleto{2,n}{4.5pt}}}(v-x)dv | = O( \log{\left( \eta_n^{-1} n^{1/2\gamma_x} \right)})$, once multiplied by $h_2/\sqrt{n\log{n}}$ goes to zero as $n \rightarrow \infty$. Note also that using the same computations as in equation \eqref{eq: written out perturbation} also for $\chi_{{\scaleto{1,n}{4.5pt}}}$, it follows that $D_{h_n} =1+o(1)$. 

Denote by: $A^{i,x}_n := \left\{n^{-\frac{1}{2\gamma_{\scaleto{x}{2.2pt}}}} \leq Z_i^x \leq \eta_n \right\} $, where $Z_i^x := Z_i - x$. As will be seen, the only contributing term is $\pi \ind_{A^{i,x}_n} / \sqrt{Z_i-x}$, we denote all the remainder terms by $R_{i}^{n,x}$ and it is immediate to see that:
\begin{align*}
& R_{i}^{n,x} = O\left(\eta_n^{-\frac{1}{2}}\right) \Big\{ \ind_{\big\{ - \eta_n \leq Z_i^x \leq - n^{-\frac{1}{2\gamma_{\scaleto{x}{2.2pt}}}} \big\}} + \ind_{\left\{ Z_i^x \leq - \eta_n  \right\}} + \ind_{A^{i,x}_n} \Big\} + \\
& O\left(n^{-\frac{1}{2\gamma_{\scaleto{x}{2.2pt}}}}(x-Z_i)^{-\frac{3}{ 2}}\right) \ind_{\big\{ - \eta_n \leq Z_i^x \leq - n^{-\frac{1}{2\gamma_{\scaleto{x}{2.2pt}}}} \big\}}  +O\left(n^{\frac{1}{4\gamma_{\scaleto{x}{2.2pt}}}}\right)\ind_{\big\{ -n^{-\frac{1}{2\gamma_{\scaleto{x}{2.2pt}}}}  \leq Z_i^x \leq n^{-\frac{1}{2\gamma_{\scaleto{x}{2.2pt}}}} \big\}}
\end{align*}
Analogously, we define $R^{n,0}_i$. Using the Taylor expansion of the logarithm with Peano residual, and \eqref{eq: expansion zeta}:
\begin{align*}
    &\sum_{i=1}^{n} \log{\frac{g_n (Z_i)}{g(Z_i)}} = \sum_{i=1}^{n} \Bigg\{ \frac{1}{m_{h_n} D_{h_n}} \bigg[ m_0 - m_{h_n} D_{h_n} + \frac{h_{1,n}}{2g(Z_i)} \int_{Z_i}^{\infty} \frac{\chi_{{\scaleto{1,n}{4.5pt}}}(v)}{\sqrt{v-Z_i}} \, dv + \\ 
    &+ \frac{h_{2,n}}{2g(Z_i)} \int_{Z_i}^{\infty} \frac{\chi_{{\scaleto{2,n}{4.5pt}}}(v-x)}{\sqrt{v-Z_i}} \, dv \bigg] - \frac{1}{2} \bigg[ \frac{m_0 - m_{h_n} D_{h_n}}{m_{h_n} D_{h_n}} + \frac{h_{1,n}}{2 m_{h_n} D_{h_n} g(Z_i)} \bigg(\frac{\pi \ind_{A^{i,0}_n}}{\sqrt{Z_i}}   \\
    & + R_n^{i,0} \bigg) + \frac{h_{2,n}}{2 m_{h_n} D_{h_n} g(Z_i)} \bigg(\frac{\pi \ind_{A^{i,x}_n}}{\sqrt{Z_i -x}}  + R_n^{i,x} \bigg) \bigg]^2 \hspace{-0.2cm} + o_p\Bigg( \bigg[ \frac{m_0 - m_{h_n} D_{h_n}}{m_{h_n} D_{h_n}} + \numberthis \label{eq: LAN exp} \\
    & \hspace{-0.05cm} + \frac{h_{1,n}}{2 m_{h_n} D_{h_n} g(Z_i)} \bigg(\frac{\pi \ind_{A^{i,0}_n}}{\sqrt{Z_i}}  + R_n^{i,0} \bigg) + \frac{h_{2,n}}{2 m_{h_n} D_{h_n} g(Z_i)} \bigg(\frac{\pi \ind_{A^{i,x}_n}}{\sqrt{Z_i -x}} + R_n^{i,x} \bigg) \bigg]^2 \Bigg) \Bigg\}.
\end{align*}
The Taylor expansion of the logarithm is granted by relation \eqref{eq: converge to zero for Taylor} (see Appendix B for its \hyperlink{proof of lemma 6}{proof}).
For $\zeta_{1,n} (z) := \int_{z}^{\infty} \frac{\chi_{{\scaleto{1,n}{3.5pt}}} (v)}{\sqrt{v-z}} \, dv, $ and for $\zeta_{2,n} (z - x) := \int_{z}^{\infty} \frac{\chi_{{\scaleto{2,n}{3.5pt}}} (v-x)}{\sqrt{v-z}} \, dv$ and $F$ as in \eqref{eq: condition roughness infimum}; then in probability: 
\begin{align}\label{eq: converge to zero for Taylor}
    \max_{1 \leq i \leq n} &\bigg|\frac{h_{1,n} \zeta_{1,n}(Z_i )}{ g(Z_i)} + \frac{h_{2,n}\zeta_{2,n} (Z_i - x)}{ g(Z_i)} \bigg| \rightarrow 0.
\end{align}
Now we provide a lemma on the normal behavior of the term inside the first square parenthesis in \eqref{eq: LAN exp}.
\begin{lemma}\label{lemma: lindeberg convergence}
For $F$ as in \eqref{eq: condition roughness infimum} and $h_{i,n}, \: i = 1,2$ as in \eqref{eq: h_i}:
\begin{align*}
    \frac{1}{m_{h_n} D_{h_n}}& \sum_{i=1}^{n} \Bigg\{ \frac{h_{1,n}}{2 g(Z_i)} \int_{Z_i}^{\infty} \frac{\chi_{{\scaleto{1,n}{4.5pt}}} (v)}{\sqrt{v-Z_i}} \, dv + \frac{h_{2,n}}{2 g(Z_i)} \int_{Z_i}^{\infty} \frac{\chi_{{\scaleto{2,n}{4.5pt}}} (v-x)}{\sqrt{v-Z_i}} \, dv + \\
    &\quad - m_{h_n} D_{h_n} + m_0 \Bigg\} \rightsquigarrow  N\left(0, \frac{\pi^2}{8  m^2_0} \left( \frac{h_1^2}{ \gamma_0 g(0) } + \frac{h_2^2}{ \gamma_x g(x) } \right) \right). \numberthis \label{eq: expressed variance}
\end{align*}
\end{lemma}
The \hyperlink{proof of lemma 7}{proof} of this lemma is provided in Appendix B. Now we can partition $m_0 - m_{h_n}D_{h_n}$ into the means showed in the definition of $\Delta_n$ (use \eqref{eq: m_{h_n}D_h}). The terms appearing in Lemma \ref{lemma: lindeberg convergence} are a linear combination of the components of $\Delta_n$. By Lemma \ref{lemma: lindeberg convergence} all such linear combinations have normal behavior. But since we can show the normal convergence for any pair of $h_i$, $i=1,2$, the Cramér-Wold device allows us to obtain the joint asymptotic normality claim for the components of $\Delta_n$. 

To determine the behavior of the leading terms inside the second square parenthesis in \eqref{eq: LAN exp}, i.e.\ :
\begin{align}\label{eq: squared terms converge}
    \sum_{i=1}^n \frac{h^2_{1,n}}{4 m^2_{h_n} D^2_{h_n} g^2(Z_i)} \frac{\pi^2 \ind_{A^{i,0}_n}}{Z_i}, \quad \quad \sum_{i=1}^n \frac{h^2_{2,n}}{4 m^2_{h_n} D^2_{h_n} g^2(Z_i)} \frac{\pi^2 \ind_{A^{i,x}_n}}{Z_i -x},
\end{align}
it is enough to use that (see Appendix B for its \hyperlink{proof of lemma 8}{proof}):
\begin{align}\label{eq: convergence empirical variance}
\frac{1}{n\log{n}} \sum_{i=1}^{n} \frac{\pi^2}{4 m^2_{h_n} g^2(Z_i)(Z_i - x)} \ind_{\left\{ \eta_n^{-1} \leq \frac{1}{Z_i -x} \leq n^{\frac{1}{2\gamma_{\scaleto{x}{2.2pt}}}} \right\}} \stackrel{\pb}{\rightarrow} \frac{ \pi^2 }{8 \gamma_x m^2_0 g(x) }.
\end{align}
The same result also applies for $A^{i,0}_n$ around $0$, and this gives convergence in probability of the terms in \eqref{eq: squared terms converge}. Now by grouping out properly the $h^{\top} = [h_1,h_2]$, we get the claim, as convergence in probability of the matrix is implied by the convergence of its components. 

We now show that all the other terms converge in probability to zero. It is enough to look at multiples of the squared terms and show that these go to zero. Relation \eqref{eq: remaining terms converge to 0} gives that these remaining terms are negligible and its \hyperlink{proof of lemma 9}{proof} is given in Appendix B. For $h_{j,n}, \: \: j = 1,2$ as in \eqref{eq: h_i} and accordingly $s \in \{0,x\}$:
\begin{align}\label{eq: remaining terms converge to 0}
     \left( \frac{m_0 - m_{h_n} D_{h_n}}{m_{h_n} D_{h_n}} \right)^2 = o\left(\frac{1}{n} \right), \quad  \sum_{i=1}^n \left( \frac{h_{j,n} R_n^{i,s}}{2 m_{h_n} D_{h_n} g(Z_i)} \right)^2 = o_p (1).   
\end{align}
We obtain also:
\begin{align*}
      o_p\Bigg(\sum_{i=1}^n &\bigg[ \frac{m_0 - m_{h_n} D_{h_n}}{m_{h_n} D_{h_n}} + \frac{h_{1,n}}{2 m_{h_n} D_{h_n} g(Z_i)} \bigg(\frac{\pi \ind_{A^{i,0}_n}}{\sqrt{Z_i}}  + R_n^{i,0} \bigg) + \\
    & + \frac{h_{2,n}}{2 m_{h_n} D_{h_n} g(Z_i)} \bigg(\frac{\pi \ind_{A^{i,x}_n}}{\sqrt{Z_i -x}} + R_n^{i,x} \bigg) \bigg]^2 \Bigg) = o_p(1),
\end{align*}
by combining \eqref{eq: convergence empirical variance} and \eqref{eq: remaining terms converge to 0}.
\end{proof}
\begin{remark}
If we take $\gamma = \gamma_0 = \gamma_x$, the obtained asymptotic covariance looks the same as the one obtained in \cite{10} (cf.\ page 2418), but our $\gamma_0, \gamma_x$ have a different meaning. 
\end{remark}

For the proof of the LAM (locally asymptotically minimax) theorem, the derivative of the map we are interested in estimating (cf.\ definition 1.10 in \cite{13}) is needed. In particular, let $\psi \: : \: \mathcal{G} \mapsto \mathbb{R}$, where $\mathcal{G}$ is the space of the probability measures $G_{F}$ associated to the random variable $Z$ (the observations). In our case, for our fixed estimation point $x>0$, $\psi(G_F) = F(x) \in \mathbb{R}$. 
\begin{lemma}\label{lemma: hadamard derivative}
    For $\gamma_0, \gamma_x$ as in Theorem \ref{thm: LAN} greater than $\frac{1}{2}$, $h^{\top} = [h_1,h_2] \in \mathbb{R}^2$, $h_n = \frac{h}{\sqrt{n\log{n}}}$, $F_{h_n}$ as in \eqref{eq: path},  $\dot{\psi} = \begin{bmatrix} \frac{1}{2\gamma_{\scaleto{0}{2.5pt}}} (1-F(x)) & - \frac{1}{2\gamma_{\scaleto{x}{2.2pt}}} \end{bmatrix}^{\top}$, as $n \rightarrow \infty$:
    \begin{align*}
        \sqrt{\frac{n}{\log{n}}} \left( \psi(G_{F_{h_n}}) - \psi(G_F) \right) = \sqrt{\frac{n}{\log{n}}} \left( F_{h_n}(x) - F(x) \right) \rightarrow h^{\top} \dot{\psi},
    \end{align*}
\end{lemma}
The \hyperlink{proof of lemma 10}{proof} of this lemma is provided in Appendix B. The combination of Theorem \ref{thm: LAN} and Lemma \ref{lemma: hadamard derivative} yields the main result Theorem \ref{thm: LAM}.

We conclude by stressing the fact that the local asymptotically minimax variance in Theorem \ref{thm: LAM} coincides with the one obtained in Theorem \ref{thm: main thm iso 1}. Moreover, as we naturally expect, both the estimator's behavior and the minimax theorem are obtained under some local conditions at $x$ and $0$ of the smoothness of $F$. The class of functions that are tightly Hölder continuous at these points as in Example \ref{example: cdf holder} satisfy the key conditions needed for the estimator and for the minimax theorem. In general, whenever the hidden $F$ satisfies both \eqref{eq: condition roughness} and \eqref{eq: condition roughness infimum}, then the isotonic estimator (IIE) will attain the minimax asymptotic variance in an adaptive way, and without the need of external tuning. 

\newpage
\begin{appendix}
\section{Proofs of complemental statements for Theorem \ref{thm: normality Holder}.}\label{subsec: proofs lemmas thm normality holder}
\vspace{-0.2cm}
Now we provide some proofs about the fundamental behavior of the covariance structure of the process $Z_n(t)$. We start by showing the case of the variance and we will proceed by giving a proof along the same lines of \eqref{eq: covariance structure final}. We will use several times how samples of $Z$ are obtained in Wicksell's problem. First, note that the observable radii in the cross sections cannot exceed the radii of the associated sphere; the spheres with bigger radius are more likely to be sampled. Thus the distribution of the squared radii of the spheres actually cut is not $F$ but a biased version of it; thus we write $X^b$ with distribution function $F^b$ for the squared radius of the spheres actually cut. Since we assume the samples are gathered uniformly, the observed squared circle radius $Z$ is related to $X^b$ via:
\begin{align}\label{eq: sampling of Z}
    Z = (1-U^2) X^b = V X^b, \quad U \sim U[0,1], \quad V \sim Be(1,1/2)
\end{align}

\begin{lemma}\label{lemma: variance without assum. g}
Fix $x >0$. If $g(x) < \infty$, then:
\begin{align*}
    \mathbb{E}\left(\sqrt{Z-x} \ind_{\{Z>x\}}-\sqrt{Z-x-\varepsilon} \ind_{\{Z>x+\varepsilon\}}\right)^2 \sim \varepsilon^2 \log \left( \frac{1}{\varepsilon} \right) g (x), \quad \text{as} \quad \varepsilon \downarrow 0.
\end{align*}
\end{lemma}
\begin{proof}
In the proof, to make the notation lighter, we will denote $X^b$ by $X$ (see \eqref{eq: sampling of Z}). Define $Z_x := Z-x$. The left side is: 
\begin{align*}
\underbrace{\mathbb{E}Z_x \ind_{\{ Z_x < \varepsilon \}}}_{(1)} + \underbrace{\mathbb{E}\left(\sqrt{Z_x} - \sqrt{Z_x - \varepsilon} \right)^2 \ind_{\left\{ Z_x > \varepsilon \right\}}}_{(2)}. 
\end{align*}
We start with $(1)$. By making a substitution below with $vX - x = y(X-x)$, which means $1-v = (1-y)\frac{X-x}{X}$, we obtain:
\begin{align*}
   \mathbb{E} Z_x \ind_{Z_x < \varepsilon}  &\stackrel{\eqref{eq: sampling of Z}}{=} \mathbb{E} \left\{ \ind_{X > x} \int_{0 < vX - x < \varepsilon} (vX -x) \frac{\ind_{0<v<1}}{2 \sqrt{1 - v}} \, dv \right\}  \\
    &\: \: = \mathbb{E} \left\{ \ind_{X> x} \int_{0}^{\frac{\varepsilon}{X - x} \wedge 1} \frac{y}{2 \sqrt{1-y}} \frac{(X-x)^{3/2}}{\sqrt{X}} \, dy\right\}.
\end{align*}
If we define $y(z) = \int_{0}^{z} \frac{y}{2 \sqrt{1-y}} \, dy$ then the above equals to:
\begin{align*}
    \mathbb{E} \left\{  y\left(\frac{\varepsilon}{X - x} \wedge 1 \right) \ind_{X> x} \frac{(X-x)^{3/2}}{\sqrt{X}} \right\}.
\end{align*}
Using: $\int \frac{y}{2 \sqrt{1-y}} \, dy = \frac{(1-y)^{3/2}}{3} - \sqrt{1-y} + c$, we can see that:
\begin{align*}
    y\left(\frac{\varepsilon}{X - x} \wedge 1\right) &= \int_{0}^{\frac{\varepsilon}{X - x} \wedge 1} \frac{y}{2 \sqrt{1-y}} \, dy =  \frac{1}{4} \left( \frac{\varepsilon}{X-x} \right)^2 + O\left( \left( \frac{\varepsilon}{X-x} \right)^3 \right).
\end{align*}
Moreover it is also true $y(z) \leq z^2, \: \: \forall \: z \in (0,1)$. The above proves:
\begin{align}\label{eq: EZa less equal}
\mathbb{E}Z_x \ind_{Z_x < \varepsilon} = O(\varepsilon^2).     
\end{align}
Therefore $(1)$ is negligible. Now we proceed with $(2)$. Using Taylor's theorem with integral form of the remainder we obtain:
\begin{align*}
    &\mathbb{E}\left( \sqrt{Z_x} - \sqrt{Z_x - \varepsilon} \right)^2 \ind_{\{ Z_x > \varepsilon \}} =  \mathbb{E} \int_{0}^{1} \int_{0}^{1} \frac{\varepsilon}{2 \sqrt{Z_x - \varepsilon s}} \frac{\varepsilon}{2 \sqrt{Z_x - \varepsilon t}} \, ds \, dt \, \ind_{\{Z_x > \varepsilon \}} \\
    &\stackrel{\eqref{eq: sampling of Z}}{=} \mathbb{E} \int_{0}^{1} \int_{0}^{1} \int_{0}^{1} \frac{\varepsilon}{2 \sqrt{vX-x - \varepsilon s}} \frac{\varepsilon}{2 \sqrt{v X - x - \varepsilon t}} \frac{1}{2 \sqrt{1-v}} \ind_{\{vX - x > \varepsilon \}} \, dv \, ds \, dt.
\end{align*}
Now we use again the substitution $vX - x = y(X-x)$ which implies $1-v = (1-y) \frac{X-x}{X}$ and $y(\varepsilon,s,t) = \int_{\varepsilon}^{1} \frac{1}{\sqrt{y - \varepsilon s}} \frac{1}{\sqrt{y - \varepsilon t}} \frac{1}{\sqrt{1-y}} \, dy$ to obtain:
\begin{align*}
    &= \frac{\varepsilon^2}{8} \mathbb{E} \int_{0}^{1} \int_{0}^{1} \int_{0}^{1} \frac{1}{\sqrt{y - \frac{\varepsilon s}{X-x}}} \frac{1}{\sqrt{y - \frac{\varepsilon t}{X-x}}} \frac{1}{ \sqrt{1-y}} \frac{\ind_{\{y(X - x) > \varepsilon \}}}{\sqrt{X-x} \sqrt{X}} \, dy \, ds \, dt \\
    &= \frac{\varepsilon^2}{8} \mathbb{E} \int_{0}^{1} \int_{0}^{1} y\left( \frac{\varepsilon}{X-x},s,t \right) \frac{\ind_{\{X - x > \varepsilon \}}}{\sqrt{X-x} \sqrt{X}} \, ds \, dt.
\end{align*}
Now we show the two following relations $\forall \: 0 < \varepsilon,s,t < 1$:
\begin{align}
    &y(\varepsilon,s,t) \sim \log{\frac{1}{\varepsilon}} - \frac{1}{2}\log{\log{\frac{1}{\varepsilon}}} + O(1) \quad \text{as} \: \varepsilon \downarrow 0 \label{eq: relations on y} \\
    & y(\varepsilon,s,t) \leq \left( \log{\frac{1}{\varepsilon}} + \log{\frac{e^2}{1 - s \vee t}}\right) \frac{1}{\sqrt{1 - s \vee t}} + 1. \label{eq: relations on y 2}
\end{align}
From this it will follow that:
\begin{align}\label{eq: essential second moment rel variance}
    \mathbb{E}\left( \sqrt{Z_x} - \sqrt{Z_x - \varepsilon} \right)^2 \ind_{\{ Z_x > \varepsilon \}} \sim \frac{\varepsilon^2}{4} g(x) \log{\frac{1}{\varepsilon}}. 
\end{align}
Fix $\varepsilon < b < 1$, then we can split the integral in the definition of $y(\varepsilon,s,t)$ into two integrals: one from $\varepsilon$ to $b$ and another one from $b$ to $1$. Let us look at the second one:
\begin{align*}
    \int_{b}^1 \frac{(1-y)^{-\frac{1}{2}}}{y \sqrt{y - \varepsilon s}} \frac{y}{\sqrt{y - \varepsilon t}} \, dy \leq& \frac{b}{b - \varepsilon} \int_{b}^1 \frac{1}{y \sqrt{1-y}} \, dy = \frac{b}{b - \varepsilon} 2 \tanh^{-1} (\sqrt{1-b}). 
\end{align*}
Here we used the fact that the mapping $y \mapsto \frac{y}{y-\varepsilon s}$ is decreasing on $[\varepsilon s,1]$, that $\frac{b}{b - \varepsilon s} \leq \frac{b}{b - \varepsilon}$ and $\int \frac{1}{y} \frac{1}{\sqrt{1-y}} \, dy = -2 \tanh^{-1}(\sqrt{1-y}) + c$. This shows that the asymptotic behavior will be determined by the integral between $\varepsilon$ and $b$;
\begin{align*}
    &\bigg| \int_{\varepsilon}^{b} \frac{1}{\sqrt{y - \varepsilon s}} \frac{1}{\sqrt{y - \varepsilon t}} \frac{1}{\sqrt{1 - y}} \, dy - \int_{\varepsilon}^{b} \frac{1}{\sqrt{y - \varepsilon s}} \frac{1}{\sqrt{y - \varepsilon t}} \, dy \bigg| \\
    &\quad \quad \quad \leq \int_{\varepsilon}^{b} \frac{1}{\sqrt{y - \varepsilon s}} \frac{1}{\sqrt{y - \varepsilon t}} \, dy \left( \frac{1}{\sqrt{1-b}} -1 \right).
\end{align*}
On the other hand, since:
\begin{align*}
\int_{\varepsilon}^b \frac{1}{y-\varepsilon s \wedge t} \, d y \leqslant \int_{\varepsilon}^b \frac{1}{\sqrt{y-\varepsilon s}} \frac{1}{\sqrt{y-\varepsilon t}} \, d y \leqslant  \int_{\varepsilon}^b \frac{1}{y-\varepsilon s \vee t} \, d y.
\end{align*}
Therefore we have:
\begin{align*}
    \int_{\varepsilon}^{b} \frac{1}{\sqrt{y - \varepsilon s}} \frac{1}{\sqrt{y - \varepsilon t}} \frac{1}{\sqrt{1 - y}} \, dy 
    \begin{cases}
        \leq \log{\left\{ \frac{b - \varepsilon s \vee t}{\varepsilon (1 - s \vee t)} \right\}} \frac{1}{\sqrt{1-b}}, \\ 
        \geq \log{\left\{ \frac{b - \varepsilon s \wedge t}{\varepsilon (1 -  s \wedge t)} \right\}}\left(2 - \frac{1}{\sqrt{1-b}} \right).
    \end{cases}
\end{align*}
Using $b = \frac{1}{\sqrt{\log{\frac{1}{\varepsilon}}}}$ we obtain for fixed $c$: 
\begin{align*}
    \log{\left\{ \frac{1/\sqrt{\log{\frac{1}{\varepsilon}}} - \varepsilon c}{\varepsilon (1-c)} \right\}} &= \log{\left\{ 1 - \varepsilon \sqrt{\log{\frac{1}{\varepsilon}}} c\right\}} + \log{\left\{ \frac{1}{\sqrt{\log{\frac{1}{\varepsilon}}}} \right\}} + \log{\frac{1}{\varepsilon}} + \\ 
    &\quad + \log{\frac{1}{1-c}} \sim \log{\frac{1}{\varepsilon}} - \frac{1}{2} \log{\log{\frac{1}{\varepsilon}}} + O(1), \quad \quad \text{as} \: \: \varepsilon \downarrow 0.
\end{align*}
This proves the first relation \eqref{eq: relations on y}. Now we proceed with relation \eqref{eq: relations on y 2}. Assume $s \leq t$. By performing a substitution $y - \varepsilon t = z (1 - \varepsilon t)$ which means $1 - y = (1 - \varepsilon t)(1 - z)$ we obtain:
\begin{align*}
    y(\varepsilon,s,t) \leq \int_{\varepsilon}^1 \frac{1}{y - \varepsilon t} \frac{1}{\sqrt{1 - y}} \, dy = \frac{1}{\sqrt{1 - \varepsilon t}} \int_{\frac{\varepsilon (1-t)}{1-\varepsilon t}}^{1} \frac{1}{z} \frac{1}{\sqrt{1-z}} \, dz. 
\end{align*}
We used the fact that we can compute the following primitive: $\int \frac{1}{y} \frac{1}{\sqrt{1-y}} \, dy = -2 \tanh^{-1}(\sqrt{1-y}) + c$. Therefore:
\begin{align*}
     &\frac{1}{\sqrt{1 - \varepsilon t}}  \int_{\frac{\varepsilon (1-t)}{1-\varepsilon t}}^{1} \frac{1}{z} \frac{1}{\sqrt{1-z}} \, dz = 2 \tanh^{-1}\left( \sqrt{1 - \frac{\varepsilon (1-t)}{1 - \varepsilon t}}\right) \frac{1}{\sqrt{1 - \varepsilon t}}   \\
     & = \frac{1}{\sqrt{1 - \varepsilon t}} \log{\left\{ \frac{1 +  \sqrt{\frac{1 - \varepsilon}{1 - \varepsilon t}} }{ 1 - \sqrt{\frac{1 - \varepsilon}{1 - \varepsilon t}}} \right\}}  = \log{\left\{ \frac{\left(\sqrt{1 - \varepsilon t} + \sqrt{1 - \varepsilon}\right)^2}{\varepsilon - \varepsilon t} \right\}} \frac{1}{\sqrt{1 - \varepsilon t}}  \\
     & \leq \bigg( \log{\frac{1}{\varepsilon}} + \log{\frac{1 - \varepsilon t}{1-t}} + 2\bigg) \frac{1}{\sqrt{1 - \varepsilon t}} \leq \left( \log{\frac{1}{\varepsilon}} + \log{\frac{e^2}{1-t}} \right) \frac{1}{\sqrt{1 - t}} +1.
\end{align*}
And this proves \eqref{eq: relations on y 2}. Now we proceed by using the two relations we just proved. First of all,  \eqref{eq: relations on y} as $\varepsilon \downarrow 0$:
\begin{align*}
    &\frac{y\left( \frac{\varepsilon}{X_x}, s, t \right) \ind_{\{X_x > \varepsilon \}}}{\log{\frac{1}{\varepsilon}}} = \frac{y\left( \frac{\varepsilon}{X_x}, s, t \right) \ind_{\{X_x > \varepsilon \}}}{\log{\frac{X_x}{\varepsilon}}} \cdot \frac{\log{X_x}  + \log{\frac{1}{\varepsilon}}}{\log{\frac{1}{\varepsilon}}} \rightarrow \ind_{\{ X > x\}}
\end{align*}
Now using $\log{xy} \leq \log{x} \ind_{x> 1} + \log{y}$ and \eqref{eq: relations on y 2} we obtain:
\begin{align*}
    \frac{y\left( \frac{\varepsilon}{X-x}, s,t \right) \mathbf{1}_{\{ X-x > \varepsilon\}}}{ \log{\frac{1}{\varepsilon}}} &\leq \frac{\left( \log{\frac{X_x}{\varepsilon}} + \log{\frac{e^2}{1 - s \vee t}} \right) \frac{1}{\sqrt{1 - s \vee t}} + 1}{\log{\frac{1}{\varepsilon}}} \\ 
    &\lesssim \left( \log{X_x} \ind_{X_x > 1} + \log{\frac{e^2}{\sqrt{1 - s \vee t}}}\right) \frac{1}{\sqrt{1 - s \vee t}} + 1.
\end{align*}
Since the idea is to apply the Dominated Convergence Theorem to:
\begin{align*}
    \frac{1}{\log{\frac{1}{\varepsilon}}} \mathbb{E} \int_{0}^{1} \int_{0}^{1} y\left( \frac{\varepsilon}{X-x},s,t \right) \frac{1}{2\sqrt{X-x} \sqrt{X}} \ind_{\{X - x > \varepsilon \}} \: ds \: dt,
\end{align*}
we need to verify that all the upper bounds have finite expectations. First $\frac{\log{X_x}\ind_{\{ X_x > 1\}}}{\sqrt{X}} \lesssim 1$ which means:
\begin{align*}
    \int_{0}^{1} \int_{0}^{1} \frac{\log{X_x} \ind_{\{X_x > 1\}} }{\sqrt{X}} \frac{1}{\sqrt{1 - s \vee t}} \, ds \, dt \lesssim 2 \int_{0}^{1} \int_{0}^{t} \frac{1}{\sqrt{1-t}} \, ds \, dt < \infty.
\end{align*}
On the other hand:
\begin{align*}
    &\int_{0}^{1} \int_{0}^{1} \left( \log{\frac{1}{1 - s \vee t}} \right) \frac{1}{\sqrt{1 - s \vee t}} \, ds \, dt \lesssim 2 \int_{0}^{1} \int_{0}^{t} \left( \log{\frac{1}{1 - t}}\right) \frac{1}{\sqrt{1-t}} \, ds \, dt \\
    &\quad \quad = 2 \int_{0}^{1} t \left( \log{\frac{1}{1-t}} \right) \frac{1}{\sqrt{1-t}} \, dt < \infty.
\end{align*}
So by the Dominated Convergence Theorem, if $g(x)= \mathbb{E} \frac{1}{2 \sqrt{X-x}} \frac{1}{\sqrt{X}} \ind_{\{X>x \}} < \infty$ we obtain the claim.
\end{proof}

\begin{proof}[\hypertarget{proof of lemma 1}{\textbf{Proof of}} \eqref{eq: covariance structure final}]
In the proof, to make the notation lighter, we will denote $X^b$ by $X$ (see \eqref{eq: sampling of Z}). Fix $s, t \in I_x$. Take the case $0 \leq s \leq t$. Fix $\mu$; we define the function $\phi(\cdot, \mu)$ as follows:
$$
\phi(z, \mu)=\sqrt{z-x} \, \ind_{[x, \infty)}(z)-\sqrt{z-x-\mu} \, \ind_{[x+\mu, \infty)}(z) .
$$
Using the i.i.d.\ property, we can write:
\begin{align}\label{eq: equation covariance}
\operatorname{Cov}\left(Z_n(s), Z_n(t)\right)= 4 \delta_n^{-2} ( \delta^*_n)^{-2}  n^{-1} \operatorname{Cov}\left(\phi\left(Z, \delta^*_n t\right), \phi\left(Z, \delta^*_n s\right)\right),
\end{align}
where $Z$ has density $g$. 
Since $g(x) < \infty$, we will prove:
\begin{align}\label{eq: covariance phi function new}
\mathbb{E} \phi\left(Z, \delta^*_n s\right) \phi\left(Z, \delta^*_n t\right)=- \frac{1}{4}st g(x) (\delta^*_n)^2 \log{\delta^*_n} +O\left((\delta^*_n)^2 \right),
\end{align}
which implies \eqref{eq: covariance structure final} for $\delta_n=n^{-1 / 2} \sqrt{\log n}$ and $\delta^*_n = (\delta_n)^{1/\gamma}$. Define $Z_x := Z-x$:
\begin{align*}
    &\mathbb{E}\left[ \sqrt{Z_x} \ind_{\{Z_x > 0 \}} - \sqrt{Z_x - \delta^*_n t} \ind_{\{Z_x > \delta^*_n t \}} \right] \hspace{-0.1cm} \left[ \sqrt{Z_x} \ind_{\{Z_x > 0 \}} - \sqrt{Z_x - \delta^*_n s} \ind_{\{Z_x > \delta^*_n s \}} \right]\\
    &= \underbrace{\mathbb{E} Z_x \ind_{\{ 0 < Z_x < \delta^*_n s\}}}_{(1)} + \underbrace{\mathbb{E} \left\{ \sqrt{Z_x} \left(\sqrt{Z_x} - \sqrt{Z_x - \delta^*_n s} \right) \ind_{\{ \delta^*_n s \leq Z_x \leq \delta^*_n t\}}\right\}}_{(2)} + \\
    &\quad + \underbrace{\mathbb{E}\left\{\left( \sqrt{Z_x} - \sqrt{Z_x - \delta^*_n t} \right) \left( \sqrt{Z_x} - \sqrt{Z_x - \delta^*_n s} \right) \ind_{\{Z_x >  \delta^*_n t \}} \right\}}_{(3)}.
\end{align*}
Let us start with $(1)$. We can use \eqref{eq: EZa less equal} to see that for $\delta^*_n$ going to zero: 
\begin{align*}
    \mathbb{E} Z_x \ind_{\{0 < Z_x < \delta^*_n s \}} \leq (\delta^*_n s)^2 \mathbb{E} \frac{1}{\sqrt{X-x}} \frac{1}{\sqrt{X}} \ind_{\{ X > x\}} = O((\delta^*_n s)^2).
\end{align*}
This shows that this term is negligible as for the covariance it will be multiplied by $4 \delta_n^{-2} ( \delta^*_n)^{-2}  n^{-1} = 4 \frac{1}{\log{n}} ( \delta^*_n)^{-2}$, will be of order $\frac{1}{\log{n}}$ and it will go to zero. \\
Let us continue with $(2)$. First note one can upper bound $(2)$ using Hölder's inequality and notice that by \eqref{eq: essential second moment rel variance} we have:
\begin{align}\label{eq: symmetric case moment squared}
    \mathbb{E} \left( \sqrt{Z_x} - \sqrt{Z_x - \delta^*_n s}\right)^2 \ind_{\{ Z_x > \delta^*_n s \}} \sim \frac{(\delta^*_n s)^2}{4} \log{\frac{1}{\delta^*_n s}} g(x).
\end{align}
Which implies that the square root of this term is of the order $O(\delta^*_n \sqrt{\log{\frac{1}{\delta^*_n}}})$. Moreover, notice that by \eqref{eq: EZa less equal}:
\begin{align*}
    \left\{ \mathbb{E} Z_x \ind_{\{ \delta^*_n s < Z_x < \delta^*_n t \}} \right\}^{1/2} \leq \left\{ \mathbb{E} Z_x \ind_{\{  Z_x > \delta^*_n s \}} \right\}^{1/2} = O(\delta^*_n ).
\end{align*}
All together we obtain that $(2)$ is of the order $O((\delta^*_n)^2 \sqrt{\log{(\delta^*_n})^{-1}})$. This quantity for the covariance structure will be multiplied by $4 \delta_n^{-2} ( \delta^*_n)^{-2}  n^{-1}$, and will be of order $O(\frac{1}{\sqrt{\log{n}}})$ thus it will go to zero because $\sqrt{\log{\frac{1}{\delta^*_n}}} \sim \sqrt{\log{n}}$. \\
Now we argue the order of $(3)$. We argue like in the proof of \eqref{eq: essential second moment rel variance}:
\begin{align*}
    (3) &= (\delta^*_n)^2 st \mathbb{E} \int_{0}^{1} \int_{0}^{1} \frac{1}{2 \sqrt{Z_x - \delta^*_n s p}} \frac{1}{2 \sqrt{Z_x - \delta^*_n t q}} \, dp \,dq \: \ind_{\{Z_x > \delta^*_n t\}} \\
    &=  \mathbb{E} \int_{0}^{1} \int_{0}^{1} \int_{0}^{1} \frac{(\delta^*_n)^2 st}{2 \sqrt{vX-x - \delta^*_n s p}} \frac{\ind_{\{vX - x > \delta^*_n t\}}}{2 \sqrt{v X - x - \delta^*_n t q}} \frac{1}{2 \sqrt{1-v}} \, dv \, dp \, dq.
\end{align*}
Now we use again the substitution $vX - x = y(X-x)$ which implies $1-v = (1-y) \frac{X-x}{X}$ and $y(\delta^*_n,s,t) = \int_{\delta^*_n}^{1} \frac{1}{\sqrt{y - \delta^*_n s}} \frac{1}{\sqrt{y - \delta^*_n t}} \frac{1}{\sqrt{1-y}} dy$ to obtain:
\begin{align*}
    &=\mathbb{E} \int_{0}^{1} \int_{0}^{1} \int_{0}^{1} \frac{(\delta^*_n)^2 st}{8 \sqrt{y - \frac{\delta^*_n t}{X-x} q}} \frac{\ind_{\{y(X - x) > \delta^*_n t\}}}{ \sqrt{y - \frac{\delta^*_n s}{X-x} p}} \frac{1}{ \sqrt{1-y}} \frac{1}{\sqrt{X-x} \sqrt{X}} \, dy \, dp \, dq \\
    &= \frac{(\delta^*_n)^2}{8} st \mathbb{E} \int_{0}^{1} \int_{0}^{1} y\left( \frac{\delta^*_n}{X-x},sp,tq \right) \frac{1}{\sqrt{X-x} \sqrt{X}} \ind_{\{X - x > \delta^*_n t \}} \, dp \, dq.
\end{align*}
Then one can obtain the claim by using \eqref{eq: relations on y} and \eqref{eq: relations on y 2}. The reason why we can apply the Dominated Convergence theorem is as in the proof of Lemma \ref{lemma: variance without assum. g}. The cases $s \leq 0 \leq t$ and $s \leq t \leq 0$ can be argued in a similar way.
\end{proof}

\begin{proof}[\hypertarget{proof of lemma 2}{\textbf{Proof of}} \eqref{eq: marginal convergence}]
We verify the Lindeberg condition for $Z_n(1)$. Recall from \eqref{eq: Z_n(1)} the definition of $Z_n(1)$. Then set $Z_n(1) = \sum_{k=1}^{n} (X^n_k - \mathbb{E}X^n_k)$ and thus $s^2_n = \sum_{k=1}^{n}\mathbb{V}(X^n_k)$, which on the other hand means:
\begin{align*}
X_k^n = \frac{2\delta_n^{-1}}{n} \frac{\sqrt{Z_k-x} \ind_{[x, \infty)}\left(Z_k\right)-\sqrt{Z_k -x-\delta_n^*} \ind_{\left[x+\delta_n^*, \infty\right)} (Z_k)}{\delta_n^*}.  
\end{align*}
Therefore by \eqref{eq: covariance structure final}: $s_n^2 = \mathbb{V}\left( Z_n(1)\right)=\frac{1}{2 \gamma} g(x) \left(1-\frac{\log \log n}{\log n}\right)+O\left(\frac{1}{\log n}\right)$. We show $\mathbb{E} X^n_k = \frac{(\delta^*_n)^{-1} \delta_n^{-1}}{n} O(\delta^*_n) = O(\frac{1}{\sqrt{n \log{n}}})$. It is enough to prove for $\delta^*_n \downarrow 0$:
\begin{align}\label{eq: enough to be shown}
    \int \left\{ \sqrt{z-x} \ind_{[x, \infty)}(z)-\sqrt{z-x-\delta_n^* } \ind_{\left[x+\delta_n^* , \infty\right)}(z) \right\} g (z) d z = O(\delta^*_n).
\end{align}
Let $Z_x := Z -x$ therefore the right-hand side of the above can be written as:
\begin{align*}
& \underbrace{\mathbb{E}\left[ \sqrt{Z_x} \ind_{\{0 < Z_x < \delta^*_n \}}  \right]}_{=(a)} +\underbrace{\mathbb{E}\left( \sqrt{Z_x} - \sqrt{Z_x - \delta^*_n} \right) \ind_{\{Z_x > \delta^*_n \}}}_{=(b)}.
\end{align*}
We show $(a) = O(\delta^*_n)$. By using \eqref{eq: sampling of Z} (we denote $X^b$ by $X$ for brevity) and making a substitution below with $vX - x = y(X-x)$, which means $1-v = (1-y)\frac{X-x}{X}$, we obtain:
\begin{align*}
    \mathbb{E} \sqrt{Z_x} \ind_{Z_x < \delta^*_n} &= \mathbb{E} \left\{ \ind_{X > x} \int_{0 < vX - x < \delta^*_n} \sqrt{vX -x} \frac{\ind_{0<v<1}}{2 \sqrt{1 - v}} \, dv \right\}  \\
    &= \mathbb{E} \left\{ \ind_{X> x} \int_{0}^{\frac{\delta^*_n}{X - x} \wedge 1} \frac{\sqrt{y}}{2 \sqrt{1-y}} \frac{(X-x)}{\sqrt{X}} \, dy \right\}.
\end{align*}
For $y(z) = \int_{0}^{z} \frac{\sqrt{y}}{2 \sqrt{1-y}} dy$, the above can be upper bounded for all $\delta^*_n$ small:
\begin{align}\label{eq: equation to upper bound Z_x}
    \mathbb{E} \left\{  y\left(\frac{\delta^*_n}{X - x} \wedge 1 \right) \ind_{X> x} \frac{X-x}{\sqrt{X}} \right\} \lesssim \mathbb{E} \delta^*_n \frac{1}{\sqrt{X}} \ind_{\{X > x \}} = O(\delta^*_n).
\end{align}
That is because: $\int_{0}^{z} \frac{\sqrt{y}}{2 \sqrt{1-y}} dy = \frac{1}{2} \sin^{-1} (\sqrt{z}) - \frac{1}{2} \sqrt{z(1-z)} \leq z$ for $z \in (0,1)$. Now we proceed with $(b)$. By using again \eqref{eq: sampling of Z} and by making a substitution below with $vX - x = y(X-x)$, which means $1-v = (1-y)\frac{X-x}{X}$, we obtain:
\begin{align*}
& \mathbb{E}\left(\sqrt{Z_x}-\sqrt{Z_x-\delta^*_n}\right) \ind_{\{Z_x>\delta^*_n\}}=\mathbb{E} \int_0^1 \frac{\delta^*_n}{2 \sqrt{Z_x-\delta^*_n s}} \: d s \: \ind_{\{Z_x>\delta^*_n\}} \\
& \quad =  \mathbb{E} \int_0^1 \int_0^1 \frac{\delta^*_n}{2 \sqrt{v X-x-\delta^*_n s}} \frac{1}{2 \sqrt{1-v}} \ind_{\{Z_x>\delta^*_n\}} \: d v \: d s\\
& \quad = \frac{1}{4} \mathbb{E} \int_0^1 \int_0^1 \frac{\delta^*_n}{\sqrt{y(X-x)-\delta^*_n s}} \frac{1}{\sqrt{1-y}} \frac{\sqrt{X}}{\sqrt{X-x}} \frac{X-x}{X} \: d y \: d s \\
& \quad = \frac{\delta^*_n}{4} \mathbb{E} \int_0^1  y\left(\frac{\delta^*_n}{X-x}, s\right) \frac{1}{\sqrt{X}} \ind_{\{X>x+\delta^*_n\}} \: d s,
\end{align*}
where $y(\varepsilon,s) = \int_{\varepsilon}^{1} \frac{1}{\sqrt{y-\varepsilon s}} \frac{1}{\sqrt{1-y}} dy$. For $b \in (\varepsilon,1)$, as done previously, we split the integral between $\varepsilon$ and $b$ and between $b$ and $1$. For the first integral we have:
\begin{align*}
    \bigg| \int_{\varepsilon}^{b} \frac{1}{\sqrt{y - \varepsilon s}} \frac{1}{\sqrt{1 - y}} \, dy - \int_{\varepsilon}^{b} \frac{1}{\sqrt{y - \varepsilon s}} \, dy\bigg| \leq \int_{\varepsilon}^{b} \frac{1}{\sqrt{y - \varepsilon s}} \, dy \left( \frac{1}{\sqrt{1-b}} -1 \right).
\end{align*}
From which we conclude:
\begin{align*}
\int_{\varepsilon}^b \frac{1}{\sqrt{y-\varepsilon s}} \frac{1}{\sqrt{1-y}} \, d y = \left\{\begin{array}{l}
\leq (2 \sqrt{b-\varepsilon s}-2 \sqrt{\varepsilon(1-s)})\left(\frac{1}{\sqrt{1-b}}-1\right), \\
\geq (2 \sqrt{b-\varepsilon s}-2 \sqrt{\varepsilon(1-s)})\left(2-\frac{1}{\sqrt{1-b}}\right).
\end{array}\right.
\end{align*}
For the second integral we have:
\begin{align*}
    \int_{b}^{1} \frac{1}{\sqrt{y-\varepsilon s}} \frac{1}{\sqrt{1-y}} \, dy &\leq \frac{1}{\sqrt{b^2 - b\varepsilon}} \left\{ \sqrt{b(1-b)} - \sin^{-1}(\sqrt{1-b}) \right\}.
\end{align*}
Putting together the above computations, for all $\delta^*_n$ small and $b \sim \sqrt{\frac{1}{\log{\frac{1}{\delta^*_n}}}}$:
\begin{align*}
    &\mathbb{E}\left(\sqrt{Z_x}-\sqrt{Z_x-\delta^*_n}\right) \ind_{\{Z_x>\delta^*_n\}} \lesssim \delta^*_n \mathbb{E} \frac{\ind_{\{X > x + \delta^*_n \}}}{\sqrt{X}} = O(\delta^*_n).
\end{align*}
This proves that \eqref{eq: enough to be shown} is true. 

Let us now look at the indicator $\ind_{\left\{\left|X^n_k-\mathbb{E}X^n_k\right|>\varepsilon s_n\right\}}$. We look separately at the cases: $ X^n_k-\mathbb{E}X^n_k > \varepsilon s_n$ and $X^n_k-\mathbb{E}X^n_k < - \varepsilon s_n$. Using the following facts:
\begin{align*}
    &\sqrt{Z_x} \ind_{\{Z_x > 0 \}} - \sqrt{Z_x - \delta^*_n} \ind_{\{Z_x > \delta^*_n \}} = \sqrt{Z_x} \ind_{\{ 0 < Z_x < \delta^*_n\}} + \Big( \sqrt{Z_x} \: +\\
    &   - \sqrt{Z_x - \delta^*_n} \Big) \ind_{\{Z_x > \delta^*_n \}} \leq \sqrt{\delta^*_n} + \int_{0}^{1} \frac{\delta^*_n \ind_{\{Z_x > \delta^*_n \}}}{2 \sqrt{Z_x - \delta^*_n s}} \, ds \leq \sqrt{\delta^*_n} + \delta^*_n \frac{ \ind_{\{Z_x > \delta^*_n \}}}{2 \sqrt{Z_x - \delta^*_n }}.
\end{align*}
So we can rewrite the set $ \{X^n_k-\mathbb{E}X^n_k > \varepsilon s_n \}$, using the fact that $\mathbb{E}X_k^n = O\left( \frac{1}{\sqrt{n \log{n}}} \right)$, and noticing that the last inclusion is valid for all $n$ big enough:
\begin{align*}
    &\left\{ \frac{\sqrt{Z_x} \ind_{\{Z_x > 0 \}} - \sqrt{Z_x - \delta^*_n} \ind_{\{Z_x > \delta^*_n \}}}{ \delta^*_n} > \varepsilon s_n \sqrt{n \log{n}}  + O(1) \right\} \\
    & \quad \subseteq \left\{ \frac{\ind_{\{Z_x > \delta^*_n \}}}{2 \sqrt{Z_x - \delta^*_n}} > \varepsilon s_n \sqrt{n \log{n}} - \frac{1}{\sqrt{\delta^*_n}} + O(1) \right\} \\
    & \quad \subseteq \left\{ Z_x < \frac{\delta^*_n}{4\left( \sqrt{\delta^*_n n \log{n}} \varepsilon s_n - 1 + O(\sqrt{\delta^*_n}) \right)^2} + \delta^*_n \right\}.
\end{align*}
Since $\frac{\delta^*_n}{4\left( \sqrt{\delta^*_n n \log{n}} \varepsilon s_n - 1 + O(\sqrt{\delta^*_n}) \right)^2} \sim \frac{1}{n \log{n}}$, for all $n$ big enough we conclude:
\begin{align*}
    \{X^n_k-\mathbb{E}X^n_k > \varepsilon s_n \} \subseteq \left\{ Z_x < \frac{1}{n\log{n}} + \delta^*_n \right\} =: A_n.
\end{align*}
Now we see that the set $\{X^n_k-\mathbb{E}X^n_k < - \varepsilon s_n\}$ is negligible, because the object on the left-hand side in the following expression is nonnegative:
\begin{align*}
   \left\{ \frac{\sqrt{Z_x} \ind_{\{Z_x > 0 \}} - \sqrt{Z_x - \delta^*_n} \ind_{\{Z_x > \delta^*_n \}}}{ \delta^*_n} < - \varepsilon s_n \sqrt{n \log{n}}  + O(1) \right\}.
\end{align*}
Therefore we see that the integral over the set $\left|X^n_k-\mathbb{E}X^n_k\right|>\varepsilon s_n$ for all $n$ big enough has the same behavior as integrating over the set: $A_n$. Using i.i.d.:
\begin{align*}
\sum_{k=1}^{n} \mathbb{E}[(X^n_k - \mathbb{E}X^n_k)^2 \ind_{A_n}] \leq \underbrace{2 n \mathbb{E}[(X^n_k)^2 \ind_{A_n}]}_{(1)} + \underbrace{2 n (\mathbb{E} X^n_k)^2 \mathbb{P}(A_n)}_{(2)} \rightarrow 0.   
\end{align*}
Since $\mathbb{E} X_k^n = O(\frac{1}{\sqrt{n \log{n}}})$ we see immediately that $(2)$ goes to 0. Now we prove that $(1)$ also does. We need to look at the following upper bound:
\begin{align*}
\frac{(\delta^*_n)^{-2}}{\log{n}} \mathbb{E} \left[ \left\{ Z_x \ind_{\{Z_x > 0 \}} + (Z_x - \delta^*_n) \ind_{\{Z_x > \delta^*_n \}} \right\} \ind_{A_n} \right] \hspace{-0.05cm} \leq \frac{2(\delta^*_n)^{-2}}{\log{n}} \mathbb{E} \left[  Z_x \ind_{\{Z_x > 0 \} \cap A_n} \right]
\end{align*}
From \eqref{eq: EZa less equal}, we recall that: $\mathbb{E}Z_x \ind_{0 < Z_x < \varepsilon} = O(\varepsilon^2)$ for $\varepsilon \downarrow 0$. We conclude that:
\begin{align*}
    &\frac{(\delta^*_n)^{-2}}{\log{n}} \mathbb{E} \left[  Z_x \ind_{\{0 < Z_x < \delta^*_n + \frac{1}{n \log{n}} \}} \right] \sim \frac{1}{\log{n}} \left( \frac{n}{\log{n}} \right)^{\frac{1}{\gamma}} \left( \left(\frac{\log{n}}{n} \right)^{\frac{1}{2\gamma}} + \frac{1}{n \log{n}} \right)^2 \sim \\
    & \quad \sim \frac{1}{\log{n}} + \frac{1}{\log{n}} \left(\frac{n}{\log{n}} \right)^{1/2\gamma} \frac{1}{n \log{n}} + \frac{1}{\log{n}}
 \left( \frac{n}{\log{n}} \right)^{1/\gamma} \frac{1}{n^2 \log^2{n}} \rightarrow 0.
 \end{align*}
 For $\gamma > \frac{1}{2}$. Therefore the Lindeberg condition is satisfied for $Z_n(1)$.
\end{proof}

\begin{proof}[\hypertarget{proof of lemma 3}{\textbf{Proof of}} \eqref{eq: stochastic boundedness seq argmax}]
We start with the class of functions $h_\theta:[0, \infty) \rightarrow \mathbb{R}$:
\begin{align*}
\mathscr{H}= \bigg\{h_\theta \: : \: h_\theta(z) =& \: 2 \sqrt{z-x} \ind_{[x, \infty)}(z)-2 \sqrt{z-\theta} \ind_{[\theta, \infty)}(z)-a_0(x-\theta), \theta \geq 0 \bigg\}
\end{align*}
and, for $\eta>0$, its subclass $ \mathscr{H}_\eta=\left\{h_\theta \in \mathscr{H}:|\theta-x| \leq \eta\right\}$.
In \cite{1}, they show that the envelope of $\mathscr{H}_\eta$, $H_\eta(z)=2\left(\sqrt{z-x} \ind_{[x, \infty)}(z)-\sqrt{z-x-\eta} \ind_{[x-\eta, \infty)}(z)\right)+a_0 \eta$, satisfies for all $\eta$ small (cf.\ (11) in \cite{1}):
\begin{align}\label{eq: bound on envelope}
\int H_\eta^2(z) g(z) \, d z \leq-2 g(x) \eta^2 \log \eta.    
\end{align}
The first thing that needs to be proved is that for $|\theta -x|$ small, we have:
\begin{align*}
    \int  h_\theta(z) - h_x (z) \, d G(z) \lesssim - |\theta - x|^{\gamma + 1}.
\end{align*}
Using a Taylor expansion for $\int h_\theta(z) d G(z)$ around $x$ (note $h_x(z) = 0$), we get, using the definition of $U$: $U(x)=\frac{\pi}{2 m_0} \int_0^x(1-F(y)) \, d y =\int_0^x V(y) \, d y$, that there exists $\eta_0 > 0$ such that:
\begin{align*}
&\int h_\theta(z) \, d G(z) = 2 \int \sqrt{z-x} \ind_{[x,\infty)}(z) -  \sqrt{z-\theta} \ind_{[\theta,\infty)}(z) \, dG(z)  - V(x) (\theta-x) \\ 
&= U(x) - U(\theta) - V(x)(x-\theta) \stackrel{(T)}{=} \frac{(\theta-x)\pi}{2 m_0} \int_{0}^{1} F(x + u(\theta-x)) -F(x) \, du  \\
&\stackrel{\text{\eqref{eq: condition roughness}}}{=}-\frac{\pi}{2 m_0} |\theta - x|^{\gamma + 1}(K + o(1)),
\end{align*}
for all $\theta$ with $|\theta-x|<\eta_0$ and for some constant $K>0$. Note that here we used the same reasoning as in \eqref{eq: reasoning asympt holder} and we also used assumption \eqref{eq: condition roughness}. 

Now we show the weak consistency of $T_n\left(a_0 + \delta_n a\right)$ for $x=V^{-1}\left(a_0\right)$. This will follow from Theorem 3, Theorem 4 in \cite{1} and the switch relation \eqref{eq: switch relation}. In other words we prove: $\forall \: \varepsilon >0$ then $\mathbb{P}\left(\left|T_n\left(a_0 + \delta_n a \right)-x\right| > \varepsilon  \right) \rightarrow 0$. Indeed using the switch relation and consistency of $\hat{V}_n$ we get:
\begin{align*}
    &\mathbb{P}\left(\left|T_n\left(a_0 + \delta_n a \right)-x\right| > \varepsilon \right) \\
    & \quad \leq \mathbb{P} \left( T_n (V(x) + \delta_n a) > \varepsilon + x \right) + \mathbb{P} \left( T_n (V(x) + \delta_n a) < x - \varepsilon \right) \\
    & \quad = \mathbb{P} \left( \hat{V}_n (x + \varepsilon) > V(x) + \delta_n a \right) + \mathbb{P} \left( \hat{V}_n (x - \varepsilon) < V(x) + \delta_n a \right) \longrightarrow 0.
\end{align*}
As by consistency $\hat{V}_n(x + \varepsilon) \rightarrow V(x+\varepsilon)$, $\delta_n \downarrow 0$ and $V$ is strictly decreasing at $x$, it will never be the case that $V (x + \varepsilon) \geq V(x)$, nor $V (x - \varepsilon) \leq V(x)$. Note that $V$ is already decreasing by definition, and it cannot be constant in a neighborhood of $x$, as we would get a contradiction with the fact that $F$ is not constant in a neighborhood of $x$. \\
Now we will proceed by deriving the rate of convergence of $T_n(a_0 + \delta_n a)$ to $x = V^{-1}(a_0)$. We will do it by a standard technique using a peeling argument, combined later with what we showed above. 

First it needs to be proved that $\forall \: \varepsilon > 0$ there exists $ \{M_n\}$ sequence of stochastically bounded random variables and a positive constant $\eta_1$ such that:
\begin{align}\label{eq: relation to be proved}
2&\left|\int\left\{\sqrt{z-x} \ind_{[x, \infty)}(z) - \sqrt{z-\theta} \ind_{[\theta, \infty)}(z)\right\} \, d\left(\mathbb{G}_n-G\right)(z)\right| \\
&\quad \quad \quad \quad \quad \quad \quad \quad \leq \varepsilon|x-\theta|^{\gamma + 1}+(\delta^*_n)^{\gamma +1} M_n,
\end{align}
for all $\theta$ with $|\theta-x|<\eta_1$, and $\gamma$ as in assumption \eqref{eq: condition roughness}. Remember that $x$ has been fixed and that $\delta_n=n^{-1 / 2} \sqrt{\log n}$, whereas $\delta^*_n = (\delta_n)^{1/\gamma}$. Now we proceed with the proof of \eqref{eq: relation to be proved}. Fix $\varepsilon > 0$, and let $\{M_n (\omega) \}$ be the random variable defined omega-wise as the infimum of all positive numbers $\eta$ so that:
\begin{align}\label{eq: equation with epsilon}
\left|\int h_\theta(z) \, d\left(\mathbb{G}_n-G\right)(z ; \omega)\right| \leq \varepsilon|x-\theta|^{\gamma + 1}+(\delta^*_n)^{\gamma + 1} \eta,
\end{align}
holds for all $\theta \in\left[x-\eta_1, x+\eta_1\right]$. If we define, for $n, j \geq 1$,
$$
A(n, j)=\left\{\theta \geq 0:(j-1) \delta^*_n \leq|\theta-x| \leq j \delta^*_n\right\},
$$
we can write, for all $n \geq 1$ and $\nu$ constant:
\begin{align*}
& \mathbb{P}\left\{M_n>\nu\right\} \leq \\
& \quad \leq \sum_{j=1}^{\eta_1 / \delta^*_n} \mathbb{P}\left\{\sup _{\theta \in A(n, j)} (\delta^*_n)^{-(\gamma + 1)}\left|\int h_\theta(z) \, d\left(\mathbb{G}_n-G\right)(z)\right|>\varepsilon(j-1)^{\gamma + 1}+\nu\right\} \\
& \quad \leq \sum_{j=1}^{\eta_1 / \delta^*_n} \frac{\mathbb{E}\left\{\sup _{|\theta-x| \leq j \delta^*_n} (\delta^*_n)^{-2(\gamma + 1)}\left[\int h_\theta(z) \, d\left(\mathbb{G}_n-G\right)(z)\right]^2\right\}}{\left[\varepsilon(j-1)^{\gamma + 1}+\nu\right]^2}. \numberthis \label{eq: equation to max bound}
\end{align*}
By maximal inequality 3.1 (ii) in \cite{12} and relation \eqref{eq: bound on envelope} (with $\eta = \delta^*_n j$), $\exists$ a $C>0$ such that, for each $j \leq \eta_1 / \delta^*_n$ and using that $(\delta^*_n)^{-2(\gamma + 1)} = (\delta^*_n)^{-2} \delta^{-2}_n$:
\begin{align*}
&\mathbb{E}\bigg\{\sup _{|\theta-x| \leq j \delta^*_n} (\delta^*_n)^{-2(\gamma + 1)}\left[\int h_\theta(z) \, d\left(\mathbb{G}_n-G\right)(z)\right]^2\bigg\}  \\
& \quad = (\delta^*_n)^{-2} \frac{n}{\log{n}} \mathbb{E}\bigg\{\sup _{|\theta-x| \leq j \delta^*_n} \left[\int h_\theta(z) \, d\left(\mathbb{G}_n-G\right)(z)\right]^2\bigg\}  \\
& \quad \leq (\delta^*_n)^{-2} \frac{J(1)^2}{\log{n}} \left( -2 g(x) j^2 (\delta^*_n)^2\log{(\delta^*_n \cdot j)}\right) =  C g(x) j^2 + O\left( (\log{n})^{-1} \right).
\end{align*}
where $J(1)$ is the uniform entropy integral appearing from the application of the maximal inequality 3.1 (ii) in \cite{12}. We argue that $J(1)$ is finite in the proof of Lemma \ref{lemma: conv weak argmax}. Applying this inequality to \eqref{eq: equation to max bound}, it follows that we can upper bound it by: $\sum_{j=1}^{\eta_1 / \delta^*_n} \frac{C g(x) j^2 + o(1)}{\left[\varepsilon(j-1)^{\gamma + 1}+\nu\right]^2}$,
by choosing $\nu$ sufficiently large, $\mathbb{P}\left\{M_n>\nu\right\}$ can be made arbitrarily small uniformly in $n$, which proves \eqref{eq: relation to be proved}. 

We are now ready to conclude that $T_n\left(a_0 + \delta_n a \right)$ is a $\delta^*_n$-consistent estimator for $x=V^{-1}\left(a_0\right)$. From the above we know with probability tending to 1, $\left|T_n\left(a_0 + \delta_n a \right)-x\right|<\eta_0 \wedge \eta_1,$
for $n \rightarrow \infty$. If $\left|T_n\left(a_0 + \delta_n a \right)-x\right|<\eta_0 \wedge \eta_1$, take $\varepsilon$ in relation \eqref{eq: equation with epsilon} equal to $\frac{\pi K}{4 m_0}$ and obtain that:
\begin{align*}
0 & =\int h_x(z) \, d \mathbb{G}_n(z) \leq \int h_{T_n\left(a_0 + \delta_n a\right)}(z) \, d \mathbb{G}_n(z) \\
& \leq \int h_{T_n\left(a_0 + \delta_n a \right)}(z) \, d G(z)+\frac{\pi K}{4 m_0}\left|x-T_n\left(a_0 + \delta_n a\right)\right|^{\gamma +1
} +(\delta^*_n)^{\gamma + 1} M_n \\
& = \left( \frac{\pi K}{4 m_0} - \frac{\pi K}{2 m_0} + o(1) \right)\left|x-T_n\left(a_0 + \delta_n a \right)\right|^{\gamma +1
}+(\delta^*_n)^{\gamma + 1} M_n,
\end{align*}
where the sequence $\{M_n\}$ is stochastically bounded by the above part. This proves that $T_n\left(a_0 + \delta_n a \right)-x=O_p\left(\delta^*_n\right)$ for $n \rightarrow \infty$. 
\end{proof}

\begin{lemma}\label{lemma: finite g if rough 2}
    Condition \eqref{eq: condition roughness 2} implies that $F$ is not constant in a neighborhood of  $x$ and that: $ \int_{x}^{\infty} \frac{dF(s)}{\sqrt{s - x}}  < \infty.$
\end{lemma}
\begin{proof}
Assume \eqref{eq: condition roughness 2}. Suppose $F$ is constant in a neighborhood of  $x$ then the limit in \eqref{eq: condition roughness 2} is zero, a contradiction. From condition \eqref{eq: condition roughness 2}, it follows that we can fix $K > \epsilon > 0$ knowing that there exists $\delta >0$ such that $\forall \: t \in \mathbb{R}^{+} \: : \: t \leq \delta$ then: $(K - \epsilon) t^{\gamma} \leq F(x + t) - F(x) \leq (K+\epsilon) t^{\gamma}$. Let us set a grid that goes from $x$ to $x + \delta$ as: $x_0 = x + \delta \geq x_1 \geq ... \geq x$, $x_j = x + 2^{-j} \delta$. It is enough to look at:
\begin{align*}
    &\int_{x}^{x+\delta} \frac{dF(s)}{\sqrt{s-x}} \leq \sum_{j=0}^{\infty} \int_{x_{j+1}}^{x_j} \frac{dF(s)}{\sqrt{s-x}} \leq \sum_{j=0}^{\infty} \frac{(K+\epsilon)2^{-\gamma (j+1)}}{\sqrt{2^{-(j+1)}\delta}} - \sum_{j=0}^{\infty} \frac{(K-\epsilon)2^{-\gamma j}}{\sqrt{2^{-(j+1)}\delta}},
\end{align*}
where the last expression is finite as $\gamma > \frac{1}{2}$.
\end{proof}
\begin{lemma}\label{lemma: bounded density}
If $F$ is locally Hölder continuous of degree $\gamma > 1/2$ around $x_0$, then the density of the observations $g$ is bounded and continuous around $x_0$.
\end{lemma}
\begin{proof}
Fix any $z$ in a neighborhood of $x_0$, so that $ \exists \: \eta >0 \, : \, F$ is Hölder continuous on $(z- \eta, z+ \eta)$. Set $z_j = z + \eta 2^{-j}$ for $j=0,1,2...$ so that $z_0 = z+\eta \geq z_1 \geq ... \geq z$ and therefore $z_j - z_{j+1} = \eta 2^{-j-1}$. We have (note $m_0 < \infty$ by \eqref{eq: finite first moment}):
\begin{align*}
    &g (z)m_0 = \int_{z}^{\infty} \frac{1}{\sqrt{x-z}} dF(x) = \sum_{j=0}^{\infty} \int_{z_{j+1}}^{z_j} \frac{1}{\sqrt{x-z}} dF(x) + \int_{z+\eta}^{\infty} \frac{1}{\sqrt{x-z}} dF(x)  \\
    &\leq\sum_{j=0}^{\infty} \frac{(F(z_j)-F(z_{j+1}))}{\sqrt{z_{j+1}-z}}  + \frac{1}{\sqrt{\eta}} \leq \frac{1}{\sqrt{\eta}} \left( \sum_{j=0}^{\infty} \frac{2^{-\gamma(j+1)}}{2^{-\frac{1}{2}(j+1)}} + 1 \right), 
\end{align*}
where the last expression is $< \infty $ if and only if $\gamma > \frac{1}{2}$. Now we know that the function $g$ is bounded around $x_0$. We have to prove that it is continuous. Using the same reasoning as in the previous part we have that:
\begin{align*}
     g (z) & =\frac{1}{m_0} \underbrace{\sum_{j=0}^{\infty} \int_{z_{j+1}}^{z_j} \frac{1}{\sqrt{x-z}} dF(x)}_{(1)} + \frac{1}{m_0} \underbrace{\int_{z+\eta}^{\infty} \frac{1}{\sqrt{x-z}} dF(x)}_{(2)}.
\end{align*}
Since in the previous part of the proof we have been able to bound $(1)$ with the convergent series independent of $z$: $\sum_{j=0}^{\infty} \frac{2^{-\gamma(j+1)}}{2^{-\frac{1}{2}(j+1)}}$, then the Weierstrass criterion for uniform convergence of series of functions (M-test) yields that the series $(1)$ converges uniformly. Now we need to show that $z \mapsto \int_{z_{j+1}}^{z_j} \frac{1}{\sqrt{x-z}} dF(x)$ are continuous $\forall \: j$. We also need to show that $(2)$ is continuous. First of all, notice that the following function: $f(s) = \frac{\ind_{[\eta 2^{-j-1}, \eta 2^{-j}]}(s)}{\sqrt{s}}$, is continuous on the compact set $[\eta 2^{-j-1},\eta 2^{-j}]$, thus uniformly continuous.
By taking $s=x-z$, one can retrieve $\frac{\ind_{[z+ \eta 2^{-j-1}, z + \eta 2^{-j}]} (x)}{\sqrt{x-z}}$. Now take $\delta$ so that whenever $|z_1 - z_2| < \delta$, by uniform continuity we have:
\begin{align*}
    \bigg| \frac{\ind_{[z_1 + \eta 2^{-j-1}, z_1 + \eta 2^{-j}]} (x)}{\sqrt{x-z_1}} - \frac{\ind_{[z_2 +\eta 2^{-j-1}, z_2 +\eta 2^{-j}]} (x)}{\sqrt{x-z_2}} \bigg| < \varepsilon.
\end{align*}
Then we have, using the triangle inequality for integrals:
\begin{align*}
    \bigg| \int \frac{\ind_{[z_1 + \eta 2^{-j-1}, z_1 + \eta 2^{-j}]} (x)}{\sqrt{x-z_1}} - \frac{\ind_{[z_2 + \eta 2^{-j-1}, z_2 + \eta 2^{-j}]} (x)}{\sqrt{x-z_2}} dF(x) \bigg| \leq \varepsilon \int dF(x) = \varepsilon.
\end{align*}
Finally, $(2)$ is continuous because $\frac{\ind_{[z+\eta,\infty)} (x)}{\sqrt{x-z}} = \frac{\ind_{[\eta,\infty)} (x-z)}{\sqrt{x-z}}$ is uniformly continuous on $[\eta, \infty)$. The Weierstrass M-test together with the uniform limit theorem yield continuity of $g$.
\end{proof}

\begin{proof}[\hypertarget{proof of lemma weak conv}{\textbf{Proof of Lemma}} \ref{lemma: conv weak argmax}]
Because
\begin{align*}
    U_n(s) := \frac{2}{n} \sum_{i=1}^n \left\{ \sqrt{(Z_i)_{+}} - \sqrt{(Z_i-s)_{+}} \right\}, 
\end{align*}
the stated convergence is equivalent to saying that the class of functions $z \rightarrow \sqrt{(z)_{+}} - \sqrt{(z-s)_{+}}, \: s \geq 0$ is Donsker. From Lemma 2.6.16 in \cite{5}, it follows that the class of functions $\{ z- s \: : \: s \geq 0 \}$ is VC, next by Lemma 2.6.18 (ii) in \cite{5} we see that the class $\{ (z -s)_{+} \: : \: s \geq 0 \}$ is VC. Again by Lemma 2.6.18 (vii), (iv) and (v) in \cite{5} we conclude that the class of functions $\{ \sqrt{(z)_{+}} - \sqrt{(z-s)_{+}} \: : \: s \geq 0 \}$ is VC. Such class has envelope: $0 \leq \sqrt{(z)_{+}} - \sqrt{(z-s)_{+}} \leq \sqrt{(z)_{+}}$. Therefore the class is $P$-Donsker for any $P$ with $\int_{0}^{\infty} z dP(z) < \infty$, which is the case here by assumption \eqref{eq: finite first moment} with $k = \frac{3}{2}$.
Now we prove that $\mathbb{Z}$ is a Gaussian process with a.s.\ continuous sample paths. Consider:
\begin{align*}
    &\ex \left( \mathbb{Z} (t_j) - \mathbb{Z} (t_i) \right)^2 =  \ex  \bigg|\sqrt{(Z - t_j)_{+}} - \sqrt{(Z - t_i)_{+}} \bigg|^2  \leq |t_j - t_i|.
\end{align*}
Where we use the fact that the square root is Hölder continuous of degree $1/2$. 

Since the intrinsic semimetric $d (s,t) := (\ex \left( \mathbb{Z} (s) - \mathbb{Z} (t) \right)^2)^{1/2} \leq \sqrt{|s-t|}$, we obtain sample path continuity of $\mathbb{Z}$. 

By the very same computation as above, we deduce that $\mathbb{V}(\mathbb{Z}(t) - \mathbb{Z}(s)) \neq 0$ for $s \neq t$ and since $\mathbb{Z}$ is indexed by a $\sigma$-compact metric space, we obtain from Lemma 2.6 in \cite{12} that the maxima of the sample paths of $\mathbb{Z}$ are a.s.\ unique. Finally the covariance structure of $\mathbb{Z}$ is given by direct computation:
\begin{align*}
    \mathrm{Cov} \left( \mathbb{Z}(t),\mathbb{Z}(s) \right) = \mathrm{Cov} \Big( \sqrt{\smash[b]{(Z)_{+}}} - \sqrt{\smash[b]{(Z_t)_{+}}}, \sqrt{\smash[b]{(Z)_{+}}} - \sqrt{\smash[b]{(Z_s)_{+}}} \Big).
\end{align*}
\end{proof}

\begin{proof}[\hypertarget{proof of lemma conv K^c_n}{\textbf{Proof of}} \eqref{eq: conv K^c_n}]
    Note that:
    \begin{align*}
        \sup_{s \in K^c_n} \left\{ Z_n(s) - \sqrt{n} h_x(s)-as \right\} \leq \sup_{s \geq 0} Z_n(s) - \inf_{\substack{s < \ubar{x}-\varepsilon_n \\ s > \bar{x}+\varepsilon_n}} \left\{ \sqrt{n} h_x(s) + as\right\},
    \end{align*}
    where the first term on the right-hand side is $O_p(1)$ by the convergence given in Lemma \ref{lemma: conv weak argmax}. The second infimum is the minimum over the infima over the intervals $[0, \ubar{x}-c]$, $[\ubar{x}-c, \ubar{x}-\varepsilon_n]$, $[\bar{x}+\varepsilon_n,\bar{x}+1]$ and $[\bar{x}+1, \infty)$ for a constant $\varepsilon_1<c<\ubar{x}$. Because $F$ is constant in $[\ubar{x},\bar{x}]$, we have $U(s)-U(x)-V(x)(s-x)= 0$ for $s \in [\ubar{x},\bar{x}]$ and from the concavity of $U$:
    \begin{align*}
        U(s)-U(x)-V(x)(s-x)< 0,  \quad \text{for} \quad s \notin [\ubar{x},\bar{x}]. 
    \end{align*} 
    Again from the fact that $U$ is concave and because $[\ubar{x},\bar{x}]$ is maximal:
    \begin{align*}
        &\forall \: \delta >0, \: \exists \: \eta >0 \: : \: U(s) - U(x) - V(x)(s-x)< -\eta |s-\ubar{x}| \quad \text{for} \quad s < \ubar{x}- \delta,\\
        &\forall \: \delta >0, \: \exists \: \eta >0 \: : \: U(s) - U(x) - V(x)(s-x)< -\eta |s-\bar{x}| \quad \text{for} \quad s > \bar{x}+ \delta.
    \end{align*} 
    Therefore using these three facts stated above we have: on $[\bar{x}+1,\infty)$,
    \begin{align*}
        \sqrt{n} h_x(s) + as \geq \sqrt{n} \eta |s-\bar{x}|- |a|s \gtrsim \sqrt{n} \rightarrow \infty
    \end{align*}
    and on $[\bar{x}+ \varepsilon_n, \bar{x}+1]$,
    \begin{align*}
        \sqrt{n} h_x(s) + as \geq \sqrt{n} h_x(\bar{x}+\varepsilon_n) - |a| \rightarrow \infty
    \end{align*}
    by the construction of $\varepsilon_n$. The other two intervals can be handled similarly.
\end{proof}

\begin{proof}[\hypertarget{proof of lemma conv prob}{\textbf{Proof of}} \eqref{eq: conv prob argmax}]
First note:
\begin{align*}
    0 &\leq \sup_{s \in F \cap K_n} \left\{Z_n(s) - \sqrt{n} h_x(s) -as\right\} - \sup_{s \in F \cap K} \left\{Z_n(s) -as\right\} \\
    & \leq \sup_{s \in F \cap K_n} \left\{Z_n(s) -as\right\} - \sup_{s \in F \cap K} \left\{Z_n(s) -as\right\}
\end{align*}
Define the process:
\begin{align*}
    \widetilde{Z}_n(s) := \begin{cases}{Z_n(s)-as,} & \text { if } s \in [\ubar{x},\bar{x}] \\ Z_n(\ubar{x})-a\ubar{x}, & \text { if } s \in [0,\ubar{x}] \\ Z_n(\bar{x})-a\bar{x}, & \text { if } s \in [\bar{x},\infty) \end{cases}
\end{align*}
because: $\sup_{s \in F \cap K_n} \widetilde{Z}_n(s) = \sup_{s \in F \cap K}  \widetilde{Z}_n(s) = \sup_{s \in F \cap K} \left\{Z_n(s) -as\right\}$, then
\begin{align*}
    &\quad \bigg| \sup_{s \in F \cap K_n} \left\{Z_n(s) -as\right\} - \sup_{s \in F \cap K} \left\{Z_n(s) -as\right\} \bigg| \\
    &= \bigg| \sup_{s \in F \cap K_n} \left\{Z_n(s) -as\right\} - \sup_{s \in F \cap K_n} \widetilde{Z}_n(s) \bigg| \leq \sup_{s \in F \cap K_n} \big| Z_n(s) - as - \widetilde{Z}_n(s) \big| \\
    & \leq \sup_{s \in F \cap K_n \cap [0,\ubar{x}]} \big| Z_n(s) - as - \widetilde{Z}_n(\ubar{x}) \big| \vee \sup_{s \in F \cap K_n \cap [\bar{x},\infty)} \big| Z_n(s) - as - \widetilde{Z}_n(\bar{x}) \big| \rightarrow 0
\end{align*}
where we used the asymptotic equicontinuity of the $Z_n$.
\end{proof}

\begin{proof}[\hypertarget{proof of main thm}{\textbf{Proof of Theorem}} \ref{thm: main thm iso}]
The claim is a consequence of $1 - F(x) = V(x)/V(0)$ and by analogy, we construct the estimator using the adaptive isotonic estimator. Apply Slutsky's Lemma together with Theorem \ref{thm: normality Holder} (cf.\ Corollary 1 in \cite{1}) to get Theorem \ref{thm: main thm iso 1}. To prove Theorem \ref{thm: main thm iso 2}, it is possible to follow the same steps as in the proof of Theorem \ref{thm: distrib flat F at x} and show the joint convergence, for all $a,b \in \mathbb{R}$:
\begin{align*}
        &\pb \bigg( \sqrt{n} \left(\hat{V}_n(x) - V(x) \right) \leq a, \sqrt{n} \left(\hat{V}_n(0) - V(0) \right) \leq b \bigg) \rightarrow\\
        &\rightarrow \pb \bigg( \argmax_{s \in K_x} \left\{ \mathbb{Z}_x(s) - a s\right\} \leq x, \argmax_{s \in K_0} \left\{ \mathbb{Z}_0(s) - b s\right\} \leq 0 \bigg) \\ 
        &= \pb  \left( L_x \leq a, L_0 \leq b \right),
\end{align*}
where in the last step we used the switch relation. Using again Slutsky's Lemma with the above convergence, we get claim. 
\end{proof}

\section{Proofs complemental statements for the LAN expansion. }\label{appendix: statements LAN}

\begin{proof}[\hypertarget{proof of lemma 5}{\textbf{Proof of Lemma}} \ref{lemma: non decreasing F_h}]
It is enough to consider the case $h_1 > 0$ and $h_2 > 0$, (depicted in Figure \ref{fig: perturbation} around $x$) and we prove in this case that $F_{h_n}$ is non-decreasing for all $n$ big enough. We check only around $x$ on $\prescript{l}{}{\mathrm{K}^x_n}$, as around $0$ is the same. From \eqref{eq: condition roughness infimum}, for any $h_2 \geq 0$, $\exists \: N^*$ such that $\forall \: n \geq N^*$, all $x_1, x_2 \in K^x_n$, $x_1 < x_2$: $F(x_2)-F(x_1) \geq \frac{h_2}{\sqrt{n\log{n}}} \log{\left\{\frac{x-x_2}{x-x_1} \right\}}$,
and therefore eventually:
\begin{align*}
    D_{h_n} \left( F_{h_n}(x_2) - F_{h_n} (x_1) \right) = F(x_2)-F(x_1) - \frac{h_2}{\sqrt{n\log{n}}} \log{\left\{\frac{x-x_2}{x-x_1} \right\}} \geq 0. 
\end{align*}
\end{proof}

\begin{lemma}\label{lemma: expansion GL}
    Let $\zeta_n(x):= \int_{x}^{\infty} \frac{\ind_{\{ \delta_n \leq |v| \leq \eta_n \}}}{v \sqrt{v-x}} \, dv$, then:
    \begin{align*}
         &\zeta_n(x) = 0 \quad  \text { if } x>\eta_n, \quad \quad \bigg|\zeta_n(x) \bigg| \lesssim \frac{\delta_n}{(-x)^{3/2}} + \frac{1}{\sqrt{\eta_n}} \quad  \text { if }-\eta_n \leq x \leq-\delta_n,\\
         &\sup_{[\delta_n, \eta_n]} \bigg|\zeta_n(x) - \frac{\pi}{\sqrt{x}} \bigg| \leq \frac{\pi}{\sqrt{\eta_n}}, \quad  \sup_{[-\delta_n, \delta_n]} \bigg|\zeta_n(x)  \bigg| \leq \frac{\pi}{\sqrt{\delta_n}}, \quad  \sup_{x \leq -\eta_n}\bigg|\zeta_n(x) \bigg| \lesssim \frac{1}{\sqrt{\eta_n}}.
    \end{align*}
\end{lemma}
\begin{remark}
    Another more compact way of stating the above lemma is (cf.\ \cite{10}):
    \begin{align}
         \zeta_n(x) =& \begin{cases}0, & \text { if } x>\eta_n, \\ 
         \pi x^{-1/2} + O\left(\eta_n^{-1/2} \right), & \text { if } \delta_n \leq x \leq \eta_n, \\  O\left(\delta_n^{-1 / 2}\right), & \text { if }-\delta_n \leq x \leq \delta_n, \\ O\left(\delta_n(-x)^{-3 / 2}\right)+O\left(\eta_n^{-1 / 2}\right), & \text { if }-\eta_n \leq x \leq-\delta_n, \\ O\left(\eta_n^{-1 / 2}\right), & \text { if } x \leq-\eta_n .\end{cases}\label{eq: approximation xi}
    \end{align}
\end{remark}
\begin{proof}
    \begin{enumerate}
        \item If $x > \eta_n$, then clearly $\zeta_n (x) = 0$, due to the indicator function.
        \item If $\delta_n \leq x \leq \eta_n$, using \eqref{eq: abel integral pi} and the triangle inequality:
        \begin{align*}
            \bigg| \zeta_n(x) - \frac{\pi}{\sqrt{x}} \ind_{x \geq 0} \bigg| &= \bigg| \int_{x}^{\eta_n}\frac{1}{v\sqrt{v-x}} \, dv - \int_{x}^{\infty}\frac{1}{v\sqrt{v-x}} \, dv \bigg| \\
            & = \bigg| \int_{\eta_n}^{\infty} \frac{1}{v\sqrt{v-x}} \, dv \bigg| \leq \int_{\eta_n}^{\infty} \frac{1}{v\sqrt{v-\eta_n}} \, dv = \frac{\pi}{\sqrt{\eta_n}}.
        \end{align*}
        \item If $-\delta_n \leq x \leq \delta_n$, using \eqref{eq: abel integral pi}:
        \begin{align*}
            \hspace{-0.25cm} \big| \zeta_n(x) \big| &= \int_{\delta_n}^{\eta_n}\frac{1}{v\sqrt{v-x}} \, dv  \leq \int_{\delta_n}^{\infty} \frac{1}{v\sqrt{v-x}} \, dv \leq \int_{\delta_n}^{\infty} \frac{1}{v\sqrt{v-\delta_n}} \, dv = \frac{\pi}{\sqrt{\delta_n}}.
        \end{align*}
        \item If $-\eta_n \leq x \leq-\delta_n$, $J(y) := \int_{y}^{1} \frac{1}{v} \left( \frac{1}{\sqrt{v+1}} - \frac{1}{\sqrt{1-v}}\right)  dv$, $H(y) := \int_{1}^{y} \frac{1}{v \sqrt{v+1}}  dv$:
        \begin{align*}
            \zeta_{n} (x) &= \int_{x}^{-\delta_n} \frac{1}{v\sqrt{v-x}} \, dv + \int_{\delta_n}^{\eta_n} \frac{1}{v \sqrt{v-x}} \, dv  \\
            & = \int_{\delta_n}^{-x} \frac{1}{v} \left( \frac{1}{\sqrt{v-x}} - \frac{1}{\sqrt{-v-x}}\right) \, dv + \int_{-x}^{\eta_n} \frac{1}{v \sqrt{v-x}} \, dv \\
            &= \frac{1}{\sqrt{-x}} \left( J\left(\frac{\delta_n}{-x}\right) + H\left(\frac{\eta_n}{-x}\right) \right).
        \end{align*}
        Because $H(\infty) = -J(0) = \log{(3+2\sqrt{2})}$ we obtain:
        \begin{align*}
            &\big| \zeta_n(x) \big| = \bigg| \frac{1}{\sqrt{-x}} \left( J\left(\frac{\delta_n}{-x}\right) - J(0) - H(\infty)+ H\left(\frac{\eta_n}{-x}\right) \right) \bigg| =\\
            &= \bigg|  \frac{1}{\sqrt{-x}} \left( \int_{0}^{\delta_n/-x} \frac{1}{v} \left( \frac{1}{\sqrt{v+1}} - \frac{1}{\sqrt{1-v}}\right) \, dv + \int_{\eta_n/-x}^{\infty} \frac{1}{v \sqrt{v+1}} \, dv \right) \bigg|\\
            &\lesssim \frac{\delta_n}{(-x)^{3/2}} + \frac{1}{\sqrt{\eta_n}}
        \end{align*}
        \item If $x \leq-\eta_n$, using the same reasoning as above:
        \begin{align*}
            \big| \zeta_n(x) \big| &= \frac{1}{\sqrt{-x}} \bigg| \int_{\delta_n/-x}^{\eta_n/-x} \frac{1}{v} \left(\frac{1}{\sqrt{v+1}} - \frac{1}{\sqrt{1-v}} \right) \bigg| \, dv \leq \frac{\log{(3+2\sqrt{2})}}{\sqrt{\eta_n}}.
        \end{align*}
    \end{enumerate}
\end{proof}

\begin{proof}[\hypertarget{proof of lemma 6}{\textbf{Proof of}} \eqref{eq: converge to zero for Taylor}]
Fix $\epsilon >0$, using the triangle inequality, we can just look separately at $p \in \{0,x\}$ (where for $j=1$, $p=0$; for $j=2$, $p=x$ and recall $h_{j,n}, \: j = 1,2$ as in \eqref{eq: h_i}):
\begin{align*}
    &\pb \left( \max_{1 \leq i \leq n} \bigg|\frac{h_{j,n} \zeta_{j,n} (Z_i - p)}{ g(Z_i)} \bigg| > \epsilon \right) \leq \sum_{i=1}^n \pb \left( \bigg|\frac{h_{j,n} \zeta_{j,n} (Z_i - p)}{ g(Z_i)} \bigg| > \epsilon \right) \\
    & \leq \sum_{i=1}^n \pb \left( \bigg|\frac{ \zeta_{j,n} (Z_i - p)}{g(Z_i)}  \bigg| > \epsilon \frac{\sqrt{n \log{n}}}{h_{j}}  \right), \quad \quad j=1,2.
\end{align*}
Consider now the expansion in \eqref{eq: expansion zeta}. Define:
\begin{align}\label{eq: def A^i_n}
    A^{i,x}_n := \big\{ n^{-\frac{1}{2\gamma_{\scaleto{x}{2.2pt}}}} \leq Z_i - x \leq \eta_n \big\},
\end{align}
and analogously $A^{i,0}_n$. Now we can use \eqref{eq: expansion zeta} and, in particular, notice that the term in \eqref{eq: only term that matters in expansion} is the one of greatest order. So if we show that for $j=1,2$ and accordingly $p \in \{0,x\}$:
\begin{align*}
    \sum_{i=1}^n \pb \left( \frac{ \pi}{g(Z_i) \sqrt{Z_i-p}} \ind_{A^{i,p}_n} >  \epsilon \frac{\sqrt{n \log{n}}}{h_{j}}  \right) \rightarrow 0,
\end{align*}
where $A^{i,p}_n$ is as in \eqref{eq: def A^i_n}, we obtain the claim. The above summation is zero for all $n$ big enough because the indicator on $A^{i,p}_n$ requires $1/\sqrt{Z_i -p}$ to be bounded by $n^{\frac{1}{4\gamma_p}} \leq \sqrt{n \log{n}} (C+o(1))$, for $\gamma_p > \frac{1}{2}$, $C>0$ and $n$ big enough. Moreover the density of the observations on $A^{i,p}_n$ is bounded away from zero.    
\end{proof}

\begin{proof}[\hypertarget{proof of lemma 7}{\textbf{Proof of Lemma}} \ref{lemma: lindeberg convergence}]
We will structure the proof as follows: first, we show that the expression inside the parenthesis in \eqref{eq: expressed variance} is mean zero, secondly, we show that by using the expansion \eqref{eq: expansion zeta} the term in line \eqref{eq: only term that matters in expansion} is the one that contributes to the normal behavior (both around $x$ and around $0$) and lastly we show that all the other terms converge to zero in probability. 
First note:
\begin{align*}
    &\int_{0}^{\infty} \int_{z}^{\infty} \frac{\chi_{{\scaleto{2,n}{4.5pt}}}(v-x)dv}{2g (z) \sqrt{v-z}}  g(z) dz = \int_{0}^{\infty} \int_{0}^{v} \frac{\chi_{{\scaleto{2,n}{4.5pt}}}(v-x)}{2\sqrt{v-z}} \, dz dv  = \int_{0}^{\infty} \frac{\chi_{{\scaleto{2,n}{4.5pt}}}(v-x)}{v^{-\frac{1}{2}}} \, dv
\end{align*}
and it is immediate to see, using the definition of $m_{h_n}$ and $m_0$, that:
\begin{align}\label{eq: m_{h_n}D_h}
    D_{h_n} m_{h_n} = m_0 + h_{1,n} \int_{0}^{\infty} \sqrt{v} \chi_{{\scaleto{1,n}{4.5pt}}} (v) \, dv + h_{2,n} \int_{0}^{\infty} \sqrt{v} \chi_{{\scaleto{2,n}{4.5pt}}}(v-x) \, dv. 
\end{align}
This proves that the expression inside the parenthesis in the statement of Lemma \ref{lemma: lindeberg convergence} is mean zero. Let $\overline{k}^x_n$ be the mean of the term first term on the right hand side of the equality in \eqref{eq: only term that matters in expansion} multiplied $1/(2m_{h_n} \sqrt{n \log{n}})$ and analogously $\overline{k}^0_n$ the same mean with $0$ instead of $x$. Define (independently of $i$ and recall \eqref{eq: h_i}):
\begin{align*}
    \overline{k}_n := \overline{k}^x_n + \overline{k}^0_n = \ex \left\{ \frac{\pi}{2 m_{h_n}} \frac{h_{2,n}}{ g(Z_i) \sqrt{Z_i - x}} \ind_{A^{i,x}_n} \right\} + \ex \left\{ \frac{\pi}{2 m_{h_n}} \frac{h_{1,n}}{ g(Z_i) \sqrt{Z_i}} \ind_{A^{i,0}_n}  \right\}, 
\end{align*}
where we used again the definition of $A^{i,p}_n$ in \eqref{eq: def A^i_n}.
Note that we just isolate the term $\overline{k}_n$ from the centering $(m_{h_n} D_{h_n} - m_0)/(m_{h_n})$ which already contained it. We need to isolate $\overline{k}_n$ for the next part (note $\overline{k}_n = o(1)$). Finally, denote by $\Sigma_h$ the asymptotic variance in \eqref{eq: expressed variance}. Using the fact that $D_{h_n} = 1 + o(1)$, we can prove that for $h_{i,n }, \: i=1,2$ as in \eqref{eq: h_i}:
\begin{align*}
 \frac{\pi}{2 m_{h_n}} \sum_{i=1}^{n} \left\{ \frac{h_{1,n}}{ g(Z_i) \sqrt{Z_i}} \ind_{A^{i,0}_n} +  \frac{h_{2,n}}{ g(Z_i) \sqrt{Z_i - x}} \ind_{A^{i,x}_n} - \overline{k}_n \right\} \rightsquigarrow N \left(0, \Sigma_h \right).
\end{align*}
We apply the Lindeberg Central Limit Theorem. If we look at:
\begin{align*}
    s^2_n := \frac{\pi^2}{4 m^2_{h_n}} \sum_{i=1}^{n} \mathbb{V} \left\{ \frac{h_{1,n}}{ g(Z_i) \sqrt{Z_i}} \ind_{A^{i,0}_n} +  \frac{h_{2,n}}{ g(Z_i) \sqrt{Z_i - x}} \ind_{A^{i,x}_n} \right\}, 
\end{align*}
the covariance of the two terms above inside the parenthesis goes to zero because for all $n$ sufficiently large $\ind_{A^0_i} \ind_{A^x_i} = 0$. Now we focus on the variance:
\begin{align*}
    \sum_{i=1}^{n} \mathbb{V} \left\{ \frac{h_2}{\sqrt{n \log{n}}} \frac{\pi}{2 m_{h_n} g(Z_i)  \sqrt{Z_i - x}} \ind_{A^{i,x}_n} \right\} = \frac{h_2^2 \pi^2}{8 \gamma_x m^2_0 g(x) } + o(1). \numberthis \label{eq: s_n convergence}
\end{align*}
Indeed using the i.i.d. assumption and leaving out the constants $(h_i^2 \pi^2) / (4m^2_0)$ (recall $m_{h_n}/m_0 \rightarrow 1$), we obtain that the left hand side of \eqref{eq: s_n convergence} is equivalent to:
\begin{align*}
    \frac{1}{\log{n}} \mathbb{E} \left[ \frac{1}{g^2(Z)(Z-x)}  \ind_{A^x_n}  \right] - \left(  \frac{1}{\sqrt{\log{n}}} \mathbb{E} \left[ \frac{1}{g(Z)\sqrt{Z-x}}  \ind_{A^x_n}  \right] \right)^2.
\end{align*}
Now analyze the first term by taking limit in $n$. Using the continuity of $g$ at $x$ and that $g(x) >0$, together with $0< \eta_n < 1$, and $\eta_n \gtrsim \frac{1}{\log^{\alpha}(n)}$:
\begin{align}\label{eq: fundamental term variance LAN}
    &\lim_{n \rightarrow \infty} \frac{1}{\log{n}} \int_{x + \frac{1}{n^{1/2\gamma_{\scaleto{x}{2.2pt}}}}}^{x + \eta_n} \frac{1}{z-x} \frac{1}{g (z)} \, dz = \frac{1}{2 \gamma_x g(x)}.
\end{align}
The second term converges to zero because:
$\lim_{n \rightarrow \infty} \frac{1}{\sqrt{\log{n}}} \int_{x + \frac{1}{n^{1/2\gamma_{\scaleto{x}{2.2pt}}}}}^{x + \eta_n} \frac{1}{\sqrt{z-x}} \, dz =  0$ just by integrating and taking the limit.
This proves the claim about the variance. Using the same relation \eqref{eq: s_n convergence} with $0$ instead of $x$, we obtain:
\begin{align}\label{eq: s_n final convergence}
    s^2_n = \frac{h_1^2 \pi^2}{8 \gamma_0 m^2_0 g(0) } + \frac{h_2^2 \pi^2}{8 \gamma_x m^2_0 g(x) } + o(1).
\end{align}

Now we look at Lindeberg condition; define: 
\begin{align*}
    B^{\epsilon}_n := \left\{\left| \frac{\pi}{2 m_{h_n}\sqrt{n \log{n}}} \left( \frac{\ind_{A^{i,x}_n}}{ g(Z_i)  \sqrt{Z_i - x}}  + \frac{\ind_{A^{i,0}_n}}{ g(Z_i)  \sqrt{Z_i}}  \right) - \overline{k}_n \right| > \epsilon s_n  \right\}.
\end{align*}
Then $\forall \: \epsilon > 0$, as $n \rightarrow \infty$, we show:
\begin{align}\label{eq: lindeberg condition 0}
    \frac{1}{ s^2_n}\sum_{i=1}^n \mathbb{E} \left[ \left( \frac{\pi h_{1,n} \ind_{A^{i,0}_n}}{ 2 g(Z_i) m_{h_n} \sqrt{Z_i}} + \frac{\pi h_{2,n} \ind_{A^{i,x}_n}}{ 2 g(Z_i) m_{h_n}   \sqrt{Z_i - x}}  - \overline{k}_n \right)^2 \ind_{B^{\epsilon}_n} \right] \rightarrow 0.
\end{align}
By using the triangle inequality, the fact that $\overline{k}_n = O(1/\sqrt{n \log{n}})$, that from \eqref{eq: s_n final convergence} $s_n$ converges to a constant, and that $\ind_{A^{i,p}_n} / g(Z_i) = \ind_{A^{i,p}_n}/g(p) +o(1)$, for $p \in \{0,x \}$ and all $n$ big enough, we can work with the superset of $B^{\epsilon}_n$:
\begin{align*}
     C^{\epsilon}_{n,i} \hspace{-0.05cm} := \hspace{-0.05cm} \left\{ \frac{\ind_{A^{i,x}_n}}{ (g(x) + o(1))  \sqrt{Z_i - x}}  + \frac{\ind_{A^{i,0}_n}}{ (g(0) + o(1))  \sqrt{Z_i}} > \frac{\epsilon 2 m_0  s_n }{\pi (n \log{n})^{-1/2} } +O(1) \right\}.  
\end{align*}
Since $(\overline{k}_n)^2 = o(1)$, the term in \eqref{eq: lindeberg condition 0} involving $(\overline{k}_n)^2$ goes to zero. Alike, the cross-product terms go to zero. So the only terms that matter are the squared terms. Thus we look at:  
\begin{align*}
     \frac{\pi^2}{4 m^2_{h_n} s^2_n n \log{n} } \sum_{i=1}^n \mathbb{E} \left[ \frac{1}{g^2(Z_i)  (Z_i - p)} \ind_{A^{i,p}_n} \ind_{C^{\epsilon}_{n,i}} \right], \quad \quad \text{for} \: p \in \{0,x\},
\end{align*}
which are zero for all $n$ sufficiently large as for all such $n$ we have $\ind_{A^{i,p}_n} \ind_{C^{\epsilon}_{n,i}} = 0$.

Now we proceed with the last part of the proof and we show that all the terms that we did not take into account in the computations above converge to zero by using the Bernstein inequality (see for instance \cite{11} page 193). We will do it only for the terms deriving from \eqref{eq: expansion zeta} involving the perturbation around $x$. The analogous terms deriving from the perturbation around $0$ can be treated in the same way.

To apply Bernstein's inequality we need to center the random variables of interest. The means of all the terms below are already contained in $(m_{h_n} D_{h_n} - m_0)/(m_{h_n})$ (as it was for $\overline{k}_n$), thus we only need to isolate them for the next part. Recall also the assumptions on $F$ and $g$ (continuous at $x$ and $0$ and $g(x) >0, g(0)>0$).
\begin{enumerate}
    \item Let $A^i_n := \big\{\frac{1}{n^{1/(2\gamma_{\scaleto{x}{2.2pt}})}} \leq Z_i-x \leq \eta_n \big\}$ and $Y^i_n =  O(\eta_n^{-1 / 2})\big(\ind_{A^i_n}/g(Z_i) - \ex \ind_{A^i_n}/g(Z_i)\big)$ then $\mathbb{V}(Y^i_n) = O(1)$. Define:
    \begin{align*}
         \widetilde{Y}^i_n = \eta_n^{1/2} Y^i_n. 
     \end{align*}
    Then $\sqrt{\mathbb{V}(Y^i_n)} = \eta_n^{-1/2} \sqrt{\mathbb{V}(\widetilde{Y}^i_n)}$ and by the assumptions on $g$ we have $|\widetilde{Y}^{i}_n| \leq L$ bounded for all $n$ big enough, therefore $\forall \: \varepsilon >0$:
    \begin{align*}
        &\mathbb{P} \left( \frac{1}{\sqrt{n \log{n}}} \bigg| \sum_{i=1}^{n} Y^i_n \bigg| > \varepsilon \right) = \mathbb{P} \left( \frac{\sqrt{\mathbb{V} (Y^i_n)}}{\sqrt{\mathbb{V} (Y^i_n)} \sqrt{n \log{n}}} \bigg| \sum_{i=1}^{n} Y^i_n \bigg| > \varepsilon \right) \numberthis \label{eq: steps to be repeated}\\
        &\leq \mathbb{P} \left( \frac{O(1)}{\sqrt{\mathbb{V} (\widetilde{Y}^i_n)} \sqrt{n \log{n}}} \bigg| \sum_{i=1}^{n} \widetilde{Y}^i_n \bigg| > \varepsilon \right) \lesssim 2 e^{ - \frac{\varepsilon^2 \log{n}}{2} \left( 1 + \frac{\varepsilon L \sqrt{\log{n}}}{O(\sqrt{\eta_n}) \sqrt{n}} \right)^{-1} } \rightarrow 0.
    \end{align*}
    \item Let $A^i_n := \left\{-\frac{1}{n^{1/(2\gamma_{\scaleto{x}{2.2pt}})}} \leq Z_i-x \leq \frac{1}{n^{1/(2\gamma_{\scaleto{x}{2.2pt}})}} \right\}$,  $Y^i_n =  O(n^{\frac{1}{4\gamma_{\scaleto{x}{2.2pt}}}})\big(\ind_{A^i_n}/g(Z_i) - \ex \ind_{A^i_n}/g(Z_i)\big)$ then $\mathbb{V}(Y^i_n) = O(1)$. Define:
    \begin{align*}
         \widetilde{Y}^i_n = n^{-\frac{1}{4\gamma_x}} Y^i_n. 
     \end{align*}
    Then $\sqrt{\mathbb{V}(Y^i_n)} = n^{\frac{1}{4\gamma_x}} \sqrt{\mathbb{V}(\widetilde{Y}^i_n)}$ and by the assumptions on $g$ we have $|\widetilde{Y}^{i}_n| \leq L$ bounded for all $n$ big enough, using the same steps as for the previous term in \eqref{eq: steps to be repeated} and because $\gamma_{\scaleto{x}{3.2pt}} > \frac{1}{2}$, $\forall \: \varepsilon >0$:
    \begin{align*}
        &\mathbb{P} \left( \frac{1}{ \sqrt{n \log{n}}} \bigg| \sum_{i=1}^{n} Y^i_n \bigg| > \varepsilon \right) \lesssim 2 e^{ - \frac{\varepsilon^2 \log{n} }{2} \left( 1 + \frac{\varepsilon  L \sqrt{\log{n}}}{ O(n^{-1/4\gamma_{\scaleto{x}{2.2pt}}}) \sqrt{n}} \right)^{-1} } \rightarrow 0.
    \end{align*}
    \item Let $A^i_n := \big\{ -\eta_n \leq Z_i-x \leq -\frac{1}{n^{1/(2\gamma_{\scaleto{x}{2.2pt}})}} \big\}$, $Y^i_n = \frac{O(n^{-1/2\gamma_{\scaleto{x}{2.2pt}}})}{g(Z_i)(x-Z_i)^{3/2}} \ind_{A^i_n} - \ex \frac{O(n^{-1/2\gamma_{\scaleto{x}{2.2pt}}})}{g(Z_i)(x-Z_i)^{3/2}} \ind_{A^i_n}$, then $\mathbb{V}(Y^i_n) = O(1)$. Define:
    \begin{align*}
         \widetilde{Y}^i_n = n^{-\frac{1}{4\gamma_x}} Y^i_n. 
     \end{align*}
    Then $\sqrt{\mathbb{V}(Y^i_n)} = n^{\frac{1}{4\gamma_x}} \sqrt{\mathbb{V}(\widetilde{Y}^i_n)}$ and by the assumptions on $g$ we have $|\widetilde{Y}^{i}_n| \leq L$ bounded for all $n$ big enough, using the same steps as for the previous term in \eqref{eq: steps to be repeated} and because $\gamma_{\scaleto{x}{3.2pt}} > \frac{1}{2}$, $\forall \: \varepsilon >0$:
     \begin{align*}
        &\mathbb{P} \left( \frac{1}{ \sqrt{n \log{n}}} \bigg| \sum_{i=1}^{n} Y^i_n \bigg| > \varepsilon \right) \lesssim 2 e^{ - \frac{\varepsilon^2 \log{n} }{2} \left( 1 + \frac{\varepsilon  L \sqrt{\log{n}}}{ O(n^{-1/4\gamma_{\scaleto{x}{2.2pt}}}) \sqrt{n}} \right)^{-1} } \rightarrow 0.
    \end{align*}
    \item Let $A^i_n := \big\{Z_i-x \leq -\eta_n \big\}$, and $Y^i_n = O(\eta_n^{-\frac{1}{2}})\big(\ind_{A^i_n}/g(Z_i) - \ex \ind_{A^i_n}/g(Z_i)\big)$ then $\mathbb{V}(Y^i_n) = O(\eta_n^{-1}) \left(\int_{0}^x 1/g(z) \, dz - x^2 \right)$ (using the assumptions on $g$ we have $\int_{0}^x 1/g(z) \, dz >0$). Define:
    \begin{align*}
         \widetilde{Y}^i_n = \eta_n^{1/2} Y^i_n. 
     \end{align*}
    Then $\sqrt{\mathbb{V}(Y^i_n)} = \eta_n^{-1/2} \sqrt{\mathbb{V}(\widetilde{Y}^i_n)}$ and by the assumptions on $g$ we have $|\widetilde{Y}^{i}_n| \leq L$ bounded for all $n$ big enough, therefore $\forall \: \varepsilon >0$:
    \begin{align*}
        &\mathbb{P} \left( \frac{1}{ \sqrt{n \log{n}}} \bigg| \sum_{i=1}^{n} Y^i_n \bigg| > \varepsilon \right) \leq \mathbb{P} \left( \frac{ O(\eta_n^{-\frac{1}{2}})}{\sqrt{\mathbb{V} (Y^i_n)} \sqrt{n \log{n}}} \bigg| \sum_{i=1}^{n} Y^i_n \bigg|  > \varepsilon \right) =\\ 
        &\mathbb{P} \left( \frac{ O(1)}{\sqrt{n}\sqrt{\mathbb{V} (\widetilde{Y}^i_n)}}  \bigg| \sum_{i=1}^{n} \widetilde{Y}^i_n \bigg|  > \varepsilon  \sqrt{\eta_n \log{n}}\right) \lesssim 2 e^{ - \frac{\varepsilon^2 \log{n}}{2(\log{n})^{\alpha}} \left( 1 + \frac{\varepsilon L \sqrt{\log{n}} }{\sigma \eta_n^{-1/2} \sqrt{n}} \right)^{-1}} \hspace{-0.3cm} \rightarrow 0.
    \end{align*}
    where $\alpha \in (0,1)$ and $\sigma^2 := \int_{0}^x 1/g(z) \, dz - x^2$.
    \item The last term can be treated identically to 1. \qedhere
\end{enumerate}
\end{proof}

\begin{proof}[\hypertarget{proof of lemma 8}{\textbf{Proof of}} \eqref{eq: convergence empirical variance}]
Let $A^i_n$ as in \eqref{eq: def A^i_n}. We are going to verify the conditions of Theorem 2.2.6. from \cite{7} to prove:
\begin{align*}
     \frac{1}{n \log{n}}  \Bigg\{ \sum_{i=1}^{n} \frac{1}{g^2(Z_i)(Z_i - x)}\ind_{A^i_n} - n \int_{x+\frac{1}{n^{1/2\gamma_{\scaleto{x}{2.2pt}}}}}^{x+\eta_n} \frac{1}{z-x} \frac{1}{g(z)} \, dz \Bigg\} \rightarrow 0. \numberthis \label{eq: convergence second moment 2}
\end{align*}
By \eqref{eq: fundamental term variance LAN} we know how $\frac{1}{\log{n}} \int_{x+\frac{1}{n^{1/2\gamma_{\scaleto{x}{2.2pt}}}}}^{x+\eta_n} \frac{1}{z-x} \frac{1}{g(z)} dz$ behaves; for this reason, if we show \eqref{eq: convergence second moment 2}, we get the claim. We verify condition $(i)$ of Theorem 2.2.6. in \cite{7}:
\begin{align*}
    n \int_{\big\{ \frac{1}{g^2(z) (z-x)} \ind_{\left\{ \eta_n^{-1} \leq \frac{1}{z-x} \leq n^{1/2\gamma_{\scaleto{x}{2.2pt}}}  \right\}} > n\log{n} \big\}} g(z) \, dz= 0,
\end{align*}
because the set over which we are integrating is empty for all $n$ big enough. Now we check condition $(ii)$ of Theorem 2.2.6. from \cite{7}; we have, using Leibniz integral rule together with de L'Hôpital rule and the fact that $\frac{1}{2\gamma_{\scaleto{x}{2.2pt}}} + 1 > \frac{1}{\gamma_{\scaleto{x}{2.2pt}}}$: 
$$\lim_{n \rightarrow \infty} \frac{1}{n^2 \log^2{n}} n \int_{x+\frac{1}{n^{1/2\gamma_{\scaleto{x}{2.2pt}}}}}^{x+\eta_n} \frac{1}{(z-x)^2} \frac{1}{g^3(z)} \, dz = 0.$$ 
The above computations now yield \eqref{eq: convergence second moment 2}. By combining \eqref{eq: convergence second moment 2} and \eqref{eq: fundamental term variance LAN}, together with: $\frac{\pi^2}{4 m^2_{h_n}} \rightarrow \frac{\pi^2}{4 m^2_0} $ we obtain the claim.
\end{proof}

\begin{proof}[\hypertarget{proof of lemma 9}{\textbf{Proof of}} \eqref{eq: remaining terms converge to 0}]
Using \eqref{eq: m_{h_n}D_h}, the definition of $h_{i,n}, \: i=1,2$ in \eqref{eq: h_i},  the fact that $\int_{0}^{\infty} \sqrt{v} \chi_n (v) \, dv = O(\sqrt{\eta_n})$ and that
$\int_{0}^{\infty} \sqrt{v} \chi_n (v-x) \, dv = O(\eta_n)$:
\begin{align*}
   \left( \frac{m_0 - m_{h_n} D_{h_n}}{m_{h_n} D_{h_n}} \right)^2 = O\left( \frac{\eta_n}{n \log{n} }\right) = O\bigg(\frac{1}{n \log^{1+\alpha}{(n)}} \bigg). 
\end{align*}
Now we look only at the terms derived from the perturbation around $x$ (i.e.\ the ones for $p= x$), as the one derived from the perturbation around $0$ behave alike. It is enough to consider all the squared terms in $(R^{i,x}_n)^2$. The ones that are at most $O(\eta_n^{-1})$, go to zero as $\frac{n  O(\log^{\alpha}{n})}{n\log{n}}  \rightarrow 0 $, $\alpha \in (0,1)$. The remaining terms converge in probability to 0 by Bernstein's inequality as in the proof of Lemma \ref{lemma: lindeberg convergence}.
\end{proof}

\begin{proof}[\hypertarget{proof of lemma 10}{\textbf{Proof of Lemma}} \ref{lemma: hadamard derivative}]
\begin{align*}
&\sqrt{\frac{n}{\log{n}}} \left( F_{h_n}(x)-F(x) \right) = \frac{\sqrt{n}}{D_{h_n}\sqrt{\log{n}}} \left( (1- D_{h_n})F(x) \right) + \frac{h_1 \int_{0}^{x} \chi_{{\scaleto{1,n}{4.5pt}}}(v)  \, dv}{D_{h_n} \log{n}}\\
& \: + \frac{h_2 \int_{0}^{x} \chi_{{\scaleto{2,n}{4.5pt}}}(v-x)  \, dv}{D_{h_n} \log{n}} = \left(-\frac{h_1 \int_{0}^{\infty} \chi_{{\scaleto{1,n}{4.5pt}}}(v)  \, dv}{D_{h_n} \log{n}} - \frac{h_2 \int_{0}^{\infty} \chi_{{\scaleto{2,n}{4.5pt}}}(v-x)  \, dv}{D_{h_n} \log{n}} \right) F(x) \\
& \: + \frac{h_1 \int_{0}^{x} \chi_{{\scaleto{1,n}{4.5pt}}}(v)  \, dv}{D_{h_n} \log{n}} + \frac{h_2 \int_{0}^{x} \chi_{{\scaleto{2,n}{4.5pt}}}(v-x)  \, dv}{D_{h_n} \log{n}}. 
\end{align*}
Now by using: $\int_{0}^{\infty} \chi_{{\scaleto{1,n}{4.5pt}}}(v)  \, dv = - (1 + o(1)) \log{n^{-1/(2\gamma_{\scaleto{0}{2.5pt}})}}$, $\int_{0}^{\infty}\chi_{{\scaleto{2,n}{4.5pt}}}(v-x)  \, dv = o(1)$,  $\int_{0}^{x} \chi_{{\scaleto{1,n}{4.5pt}}}(v)  \, dv = -(1 + o(1)) \log{n^{-1/(2\gamma_{\scaleto{0}{2.5pt}})}}$, and that $\int_{0}^{x} \chi_{{\scaleto{2,n}{4.5pt}}}(v-x)  \, dv = (1 + o(1)) \log{n^{-1/(2\gamma_{\scaleto{x}{2.2pt}})}}$, the above for all $n$ big is equal to (using also $1/D_{h_n} \rightarrow 1$):
\begin{align*}
    &=\frac{h_1}{D_{h_n} \log{n}} (1 + o(1)) (F(x)-1) \log{n^{-\frac{1}{2\gamma_{\scaleto{0}{2.5pt}}}}} - \frac{h_2}{D_{h_n} \log{n}} ( o(1)) F(x)  +\\
    &  \frac{h_2}{D_{h_n} \log{n}}(1 + o(1)) \log{n^{-\frac{1}{2\gamma_{\scaleto{x}{2.5pt}}}}} \rightarrow \frac{h_1}{2\gamma_0} (1-F(x)) - \frac{h_2}{2\gamma_x} = h^{\top} \begin{bmatrix} \frac{1}{2\gamma_{\scaleto{0}{2.5pt}}} (1-F(x)) \\ - \frac{1}{2\gamma_{\scaleto{x}{2.2pt}}} \end{bmatrix}.
\end{align*}
\end{proof}

\begin{proof}[\hypertarget{proof of thm 2}{\textbf{Proof of Theorem}} \ref{thm: LAM weaker}]
Because for any neighborhood in supremum norm of $F$ there exists $n \in \mathbb{N}$ such that the "least favorable" perturbation $F_{h_n}$ (cf.\ \eqref{eq: path}) is in such neighbourhood, the expression in \eqref{eq: LAM simplified} is lower-bounded by the expression in \eqref{eq: LAM} (cf.\ pag.\ 117 in \cite{6}).
\end{proof}
\end{appendix}


\newpage


\begin{thebibliography}{}
\bibitem[1]{1} \textsc{Groeneboom, P.} and \textsc{Jongbloed, G.} (1995). Isotonic Estimation and Rates of Convergence in Wicksell's problem. \textit{The Annals of Statistic}. Vol.\ 23, No.\ 5, 1528-1542.  \url{https://doi.org/10.1214/aos/1176324310}

\bibitem[2]{2}  \textsc{Groeneboom, P.} and \textsc{Jongbloed, G.} (2014). \textit{Nonparametric Estimation under Shape Constraints}. Cambridge University Press.

\bibitem[3]{3} \textsc{Jongbloed, G.} (2001). Sieved Maximum Likelihood Estimation in Wicksell's Problem and Related Deconvolution Problems. \textit{Scandinavian Journal of Statistics} Vol.\ 2, No.\ 6, 1152-1174. \url{https://doi.org/10.1111/1467-9469.00230}

\bibitem[4]{4} \textsc{Ghosal, S.} and \textsc{van der Vaart, A.} (2017). \textit{Fundamentals of Nonparametric Bayesian Inference}. Cambridge University Press.

\bibitem[5]{5} \textsc{van der Vaart, A.} and \textsc{Wellner, J.} (1996). \textit{Weak convergence and empirical processes}. Springer, New York. 

\bibitem[6]{6} \textsc{van der Vaart, A.} (1998). \textit{Asymptotic Statistics}. Cambridge University Press. 


\bibitem[7]{7} \textsc{Durrett, R.} (2010). \textit{Probability: Theory and Examples}, Cambridge University Press.

\bibitem[8]{8}  \textsc{Feller, W.} (1971). \textit{An Introduction to Probability Theory and its Applications.} Volume II. John Wiley and Sons, Inc.

\bibitem[9]{9} \textsc{Stoyan, D.} and \textsc{Kendall, W.\ S.} (1987). \textit{Stochastic Geometry and Its Applications}. Wiley, New York.

\bibitem[10]{10} \textsc{Golubev, G.\ K.} and \textsc{Levit, B.\ Y.} (1998). Asymptotically Efficient Estimation in the Wicksell Problem. \textit{The Annals of Statistics}, Vol. 26, No. 6, 2407–2419.  \url{https://doi.org/10.1214/aos/1024691477}

\bibitem[11]{11} \textsc{Pollard, D.} (1984). \textit{Convergence of Stochastic Processes}. Springer, New York.

\bibitem[12]{12} \textsc{Kim, J.} and \textsc{Pollard, D.} (1990). Cube Root Asymptotics. \textit{The Annals of Statistics}, Vol. 18, No. 1, 191-219. \url{https://doi.org/10.1214/aos/1176347498}

\bibitem[13]{13} \textsc{Bolthausen, E.}, \textsc{van der Vaart, A.} and \textsc{Perkins, E.} (1999). \textit{Lectures on Probability Theory and Statistics, Ecole d'Eté de Probabilités de Saint-Flour XXIX} - Part III. Aad van der Vaart. Springer, New York.

\bibitem[14]{14} \textsc{Wicksell, S.\ D.} (1925). The Corpuscle Problem, \textit{Biometrika} \textbf{17} 84-99. \url{https://doi.org/10.1093/biomet/17.1-2.84}

\bibitem[15]{15} \textsc{Watson, G.\ S.} (1971). Estimating functionals of particle size distributions. \textit{Biometrika} \textbf{58}, 483-490. \url{https://doi.org/10.1093/biomet/58.3.483}

\bibitem[16]{16} \textsc{Hall, P.} and \textsc{Smith, R.\ L.} (1988). The kernel method for unfolding sphere size distributions. \textit{J.\ Comput.\ Phys.\ } \textbf{74}, 409-421. \url{https://doi.org/10.1016/0021-9991(88)90085-X}

\bibitem[17]{17} \textsc{Golubev, G.\ K.} and \textsc{Enikeeva,  F.\ N.}  (2001). Asymptotically Efficient Smoothing in the Wicksell Problem under Squared Losses. \textit{Problems of Information Transmission}, Vol. 37, No. 1. \url{https://doi.org/10.1023/A:1010495609901}

\bibitem[18]{18} \textsc{Chan, K.\ C.\ G.} and \textsc{Qin, J.} (2016). Nonparametric maximum likelihood estimation for the multisample Wicksell corpuscle problem. \textit{Biometrika} \textbf{103}, 273–286.  \url{https://doi.org/10.1093/biomet/asw011}

\bibitem[19]{19} \textsc{Sen, B.} and \textsc{Woodroofe,  F.\ N.}  (2012). Bootstrap confidence intervals for isotonic estimators in a stereological problem. \textit{Bernoulli}, \textbf{18(4)}, 1249–1266. \url{https://doi.org/10.3150/12-BEJ378}

\bibitem[20]{20} \textsc{Deng, H.}, \textsc{Han, Q.} and \textsc{Zhang, C.} (2021). Confidence intervals for multiple isotonic regression and other monotone models. \textit{The Annals of Statistics}, Vol. 49, No. 4, 2021–2052. \url{https://doi.org/10.1214/20-AOS2025}

\bibitem[21]{21} \textsc{Lopez-Sanchez, M.\ A.} and \textsc{LLana-Fúnez, s.} (2016). An extension of the Saltykov method to quantify 3D grain size distributions in mylonites. \textit{Journal of Structural Geology}, Vol.\ 93, 149-161. \url{https://doi.org/10.1016/j.jsg.2016.10.008}

\bibitem[22]{22} \textsc{Cuzzi, J.\ } and \textsc{Olson, D.\ } (2017). Recovering 3D particle size distributions from 2D sections. \textit{Meteoritics \& Planetary Science} Vol.\ 52, Nr.\ 3, 532-545. \url{https://doi.org/10.1111/maps.12812}

\bibitem[23]{23} \textsc{van der Jagt, T.\ }, \textsc{Jongbloed, G.\ } and \textsc{Vittorietti, M.\ } (2023). Existence and Approximation of Densities of Chord Length- and Cross Section Area Distributions. \textit{Image Analysis and Stereology}, 42(3), 171–184. https://doi.org/10.5566/ias.2923 \url{https://doi.org/10.5566/ias.2923}

\bibitem[24]{24} \textsc{Koshevnik, Yu.\ A.\ } and \textsc{Levit, B.\ Ya.\ } (2023). On a Non-Parametric Analogue of the Information Matrix. \textit{Theory of Probability and its Applications}, Vol.\ 21, Iss.\ 4.\ \url{https://doi.org/10.1137/1121087}


\end{thebibliography}
\end{document}